\documentclass[a4paper,10pt,twoside]{scrartcl}

\usepackage[utf8]{inputenc}
\usepackage{amsmath,amsthm,amssymb}
\usepackage{moreverb}
\usepackage{array}
\usepackage{dsfont}
\usepackage{bbm}
\usepackage{pst-all}
\usepackage{fancyhdr}
\usepackage{graphicx}
\usepackage{float}
\usepackage[left=3cm,right=3cm,top=4cm,bottom=3cm]{geometry}
\usepackage{color}
\usepackage[T1]{fontenc}
\usepackage{hyperref}

\hypersetup{
pdfauthor={Gerold Alsmeyer, Andrea Winkler},
pdftitle={Metabassins.tex},
colorlinks=false,
pdfdisplaydoctitle=true,
pdfborder={0 0 0},
plainpages = false
}

\setkomafont{sectioning}{\bf\normalcolor\scshape}
\setkomafont{descriptionlabel}{\normalcolor\scshape}

\title{Metabasins - a State Space Aggregation for highly disordered Energy Landscapes}
\date{\today}
\author{Gerold Alsmeyer and Andrea Winkler\\ \small{Institut für Mathematische Statistik, Universität Münster,}\\*[-0.2cm] \small{Einsteinstraße 62, DE-48149 Münster, Germany}}

\newenvironment{mydescription}[1]
  {\begin{list}{}
   {\renewcommand\makelabel[1]{##1\hfill}
   \settowidth\labelwidth{\makelabel{#1}}
   \setlength\leftmargin{\labelwidth}
   \addtolength\leftmargin{\labelsep}}}
  {\end{list}}

\newtheoremstyle{theorem}
  { }
  { }
  {\normalfont\itshape}
  { }
  {\bf}
  {.}
  { }
  { }
\theoremstyle{theorem}
\newtheorem{thm}{Theorem}
\newtheorem{lem}[thm]{Lemma}
\newtheorem{prop}[thm]{Proposition}

\newtheoremstyle{definition}
  { }
  { }
  {\normalfont}
  { }
  {\bf}
  {.}
  { }
  { }
\theoremstyle{definition}
\newtheorem{defn}[thm]{Definition}
\newtheorem{bsp}[thm]{Example}
\newtheorem{bem}[thm]{Remark}
\newenvironment{bew}{\begin{proof}[Proof:]}{\end{proof}}

\numberwithin{thm}{section}

\newcommand{\argmin}{\mathop{\mathrm{argmin}}}
\newcommand{\argmax}{\mathop{\mathrm{argmax}}}
\newcommand{\Prob}{\mathbb{P}}
\newcommand{\wh}[1]{{\widehat #1}}

\begin{document}
\pagestyle{fancy} 
\fancyhf{}
\fancyhead[EL]{\thepage \hfill \leftmark}
\fancyhead[OR]{\leftmark \hfill \thepage}
\maketitle
\thispagestyle{empty}

%
%

\begin{abstract}
Glass-forming systems, which are characterized by a highly disordered energy landscape, have been studied in physics by a simulation-based state space aggregation. The purpose of this article is to develop a path-independent approach within the framework of aperiodic, reversible Markov chains with exponentially small transition probabilities which depend on some energy function. This will lead to the definition of certain metastates, also called metabasins in physics. More precisely, our aggregation procedure will provide a sequence of state space partitions such that on an appropriate aggregation level certain properties (see Properties 1--4 of the Introduction) are fulfilled. Roughly speaking, this will be the case for the finest aggregation such that transitions back to an already visited (meta-)state are very unlikely within a moderate time frame.
\end{abstract}


\vspace{.1cm}

%
%

\section*{Introduction}
\addcontentsline{toc}{section}{Introduction}
\markboth{Introduction}{Introduction}

Supercooled liquids of glass forming systems are typical examples of high-dimensional systems with highly disordered energy landscapes and our main concern behind this work. Simulations have shown that many important characteristics of such a system are better described by a process on the set of so-called \textit{metabasins} (MB) than by the more common process on the set of visited minima of the energy landscape (see \cite{He08}). Those MB are formed in the following way by aggregation of suitable states of the describing process $(X_{n})_{n\ge 0}$ along a simulated trajectory:

\vspace{.1cm}
Fixing a reasonable observation time $T$, define $\chi_0\equiv0$ and then recursively for $n\geq 1$
\begin{equation*}
\chi_n:=\inf\big\{k>\chi_{n-1}\,|\,\{X_{k}{,}...,X_{T}\}\cap\{X_{0}{,}...,X_{k-1}\}=\emptyset\big\}.
\end{equation*}
Then the MB up to $\upsilon:=\sup\{n\geq 0\,|\,\chi_n\leq T\}$ are chosen as
\begin{equation*}
\mathcal{V}_{n}:=\{X_{\chi_n}{,}...,X_{\chi_{n+1}-1}\},\quad 0\leq n\leq\upsilon.
\end{equation*}

\vspace{.1cm}
Simulation studies have shown that local sampling within a MB does not affect typical parameters of the process like the diffusion coefficient or the time to reach equilibrium. Dynamical aspects are therefore fully characterized by the MB-valued process. Furthermore, this model reduction by aggregation, as proposed in \cite{He08} and \cite{HeRu}, offers several advantages (referred to as Properties 1--5 hereafter):
\begin{enumerate}
\item The probability of a transition from one MB to any other one does not depend on the state at which this MB is entered.
\item There are basically no reciprocating jumps between two MB. This is in strong contrast to the unaggregated process where such jumps occur very often: The system falls back to a minimum many times before eventually cresting a high energy barrier and then falling into a new valley, where it will again take many unsuccessful trials to escape. These reciprocating jumps are not only irrelevant for the actual motion on the state space but also complicating the estimation of parameters like the diffusion coefficient or the relaxation time.
\item The expected time spent in a MB is proportional to its depth. Thus there is a strong and explicit relation between dynamics and thermodynamics, not in terms of the absolute but the relative energy.
\item The energy barriers between any two MB are approximately of the same height, that is, there is an energy threshold $E_0$ such that, for a small $\varepsilon$, it requires a crossing of at least $E_0-\varepsilon$ and at most $E_0+\varepsilon$ to make a transition from one MB to another. Such systems with $\varepsilon=0$ are called trap models (see \cite{Bo92}).
\item The sojourn times and the jump distances between successively visited MB (measured in Euclidean distance) form sequences of weakly or even uncorrelated random variables, and are also mutually independent, at least approximately. Therefore, the aggregated process can be well approximated by a continuous time random walk, which in turn simplifies its analysis and thus the analysis of the whole process. 
\end{enumerate}

Despite these advantages, the suggested definition of MB has the obvious blemish that it depends on the realization of the considered process and may thus vary from simulation to simulation. To provide a mathematically stringent definition of a \emph{path-independent aggregation} of the state space, which maintains the above properties and is based on the well-established notion of metastable states, is therefore our principal concern here with the main results being Theorem \ref{thm:MB} and Theorem \ref{thm:pd versus pid MB}. In this endeavor, we will draw on some of the ideas developed by \textsc{Bovier} in \cite{Bo06} and by \textsc{Scoppola} in \cite{Scop}, most notably her definition of metastable states.

Metastability, a phenomenon of ongoing interest for complex physical systems described by finite Markov processes on very large state spaces, can be defined and dealt with in several ways. It has been derived from a renormalization procedure in \cite{Scop2}, by a pathwise approach in \cite{CaGaOlVa84}, and via energy landscapes in \cite{Bo06}, the latter being also our approach hereafter. To characterize a supercooled liquid, i.e.\ a glass forming system at low temperature, via its energy landscape was first done by \textsc{Goldstein} in 1969 \cite{Go69} and has by now become a common method. The general task when studying metastability, as well originally raised in physics (\cite{KiThWo89}, \cite{Sa98} or  \cite{HeRu}), is to provide mathematical tools for an analysis of the property of thermodynamical systems to evolve in state space along a trajectory of unstable or metastable states with very long sojourn times.

Inspired by simulations of glass forming systems at very low temperatures with the Metropolis algorithm, we will study (as in \cite{Scop}) finite Markov chains with exponentially small transition probabilities which are determined by an energy function and a parameter $\beta>0$. This parameter can be understood as the inverse temperature and we are thus interested in the behavior of the process as $\beta\to \infty$. We envisage an energy function of highly complex order and without the hierarchical ordering that is typical in spin glass models. 
A good picture is provided by randomly chosen energies with correlations between neighbors or by an energy landscape that looks like a real mountain landscape. We will show that, towards an aggregation outlined above, the metastable states as defined in \cite{Scop} are quite appropriate because they have an ordering from a kind of ``weak'' to a kind of ``strong'' metastability. Around those states we will define and then study connected valleys [Definition \ref{def:Tal}] characterized by minimal energy barriers. 
In the limit of low temperatures, any such barrier will determine the speed, respectively probability of a transition between the two valleys it separates. More precisely, in the limit $\beta\to\infty$, the process, when starting in a state $x$, will almost surely reach a state with lower barrier earlier than a state with higher barrier [Theorem \ref{thm:DriftZumMinimum}]. In the limit of low temperatures, the bottom (minimum) of an entered valley will therefore almost surely be reached before that valley is left again [Proposition \ref{prop:Übergänge}]. As a consequence, the probability for a transition from one valley to another is asymptotically independent of the state where the valley is entered. This is Property 1 above. 

Furthermore, since valleys as well as metastable states have a hierarchical ordering, we can build valleys of higher order by a successive merger of valleys of lower order [Proposition \ref{prop:TaelerRekursion}]. Given an appropriate energy landscape, this procedure can annihilate (on the macroscopic scale) the accumulation of reciprocating jumps by merging valleys exhibiting such jumps into a single valley [Subsection \ref{subsec:ÜWkeiten}]. Hence, valleys of sufficiently high order will have Property 2. 

Beside the macroscopic process [Section \ref{sec:MakroskopischerProzess}], which describes the transitions between valleys, one can also analyze the microscopic process [Section \ref{sec:MikroskopischerProzess}], that is, the system behavior when moving within a fixed valley. Here we will give a formula for the exit time and connect it with its parameters [Theorem \ref{thm:Verlassenszeit}]. This will confirm Property 3. 

Having thus established Properties 1--4 [Theorem \ref{thm:MB}], we will finally proceed to a comparison of our path-independent definition of MB with the path-dependent one given above. It will be shown [Theorem \ref{thm:pd versus pid MB}] that both coincide with high probability under some reasonable conditions on the connectivity of valleys which, in essence, ensure the existence of reasonable path-dependent MB. We will also briefly touch on the phenomenon of quasi-stationarity [Proposition \ref{prop:Eigenwerte}] which is a large area \cite{Po10} but to our best knowledge less studied in connection with the aggregation of states of large physical systems driven by energy landscapes.

Let us mention two further publications which, despite having a different thrust, provide definitions of valleys, called basins of attraction or metastates there, to deal with related questions. \textsc{Olivieri \& Scoppola} \cite{OlSc96} fully describe the tube of exit from a domain in terms of which basins of attraction of increasing order are visited during a stay in that domain and for how long these basins are visited. In a very recent publication, \textsc{Beltr\'an \& Landim} \cite{BeLa11}, by working with transition rates instead of energies, aim at finding a universal depth (and time scale) for all metastates. However, we rather aim at the \emph{finest aggregation such that transitions back to an already visited metastate are very unlikely within a time frame used in simulations}. This finest aggregation will lead to valleys of very variable depth just as simulations do not exhibit a universal depth or timescale.

%
%

\section{Valleys}\label{sec:valleys}

Let $X$ be a Markov chain on a finite set $\mathcal{S}$ with transition matrix $\mathbf{P}=(p(r,s))_{r,s\in \mathcal{S}}$ and stationary distribution $\pi$, and let $E:\mathcal{S}\to \mathbb{R}$ be an energy function such that the following conditions hold:
\begin{description}
\item[Irreducibility:] $\mathbf{P}$ is irreducible with $p(s,s)>0$ and $p(r,s)>0$ iff $p(s,r)>0$ for all $r,s\in \mathcal{S}$.
\item[Transition Probabilities:] There exist parameters $\beta>0$ and $\gamma_{\beta}>0$ with $\gamma_{\beta}\to 0, \beta \gamma_{\beta}\to \infty$ as $\beta \to \infty $ such that
\[e^{-\beta((E(s)-E(r))^++\gamma_{\beta}/|\mathcal{S}|)}
\leq p(r,s)
\leq e^{-\beta((E(s)-E(r))^+-\gamma_{\beta}/|\mathcal{S}|)}\]
for all distinct $r,s\in \mathcal{S}$ with $p(r,s)>0$. Furthermore,
\begin{align*}
p^{*}(r,s):=\lim_{\beta\to\infty}p(r,s)
\end{align*}
exists for all $r,s\in\mathcal{S}$, is positive if $E(r)\geq E(s)$ and $=0$ otherwise.
\item[Reversibility:] The pair $(\pi,\mathbf{P})$ satisfies the detailed balance equations, i.e.
\[\pi(r)p(r,s)=\pi(s)p(s,r)\]
for all $r,s\in \mathcal{S}$.
\item[Non-Degeneracy:] $E(r)\neq E(s)$ for all $r,s \in \mathcal{S}, r\neq s$.
\end{description}

We are thus dealing with a reversible Markov chain with exponentially small transition probabilities driven by an energy landscape. As an example, which is also the main motivation behind this work, one can think of a Metropolis chain with transition probabilities of the form
\[p(r,s)=\frac{1}{C(r)}e^{-\beta(E(s)-E(r))^+}.\]
Here $\beta$ is the inverse temperature and $C(r), r\in\mathcal{S}$, is a parameter, independent from $\beta$, giving the number of neighbors of $r$. For $\gamma_{\beta}/|\mathcal{S}|:=\max_{r\in \mathcal{S}}\ln(C(r))(\beta+1)^{-1/2}$ the above conditions are fulfilled. Let us start with the following basic result for the ratios of the stationary distribution.

\begin{lem}\label{lem:statVert}
For any two states $r,s\in\mathcal{S}$ with $E(r)>E(s)$, we have
\[e^{-\beta(E(r)-E(s)+2\gamma_{\beta})}
\leq \frac{\pi(r)}{\pi(s)}\leq e^{-\beta(E(r)-E(s)-2\gamma_{\beta})}.\]
\end{lem}

\begin{bew}
To start with, assume $r\sim s$. Reversibility and the assumptions on the transition probabilities imply
\begin{align*}
\frac{\pi(r)}{\pi(s)}\ 
&=\ \frac{p(s,r)}{p(r,s)}\ 
\leq\ \frac{e^{-\beta((E(r)-E(s))^+-\gamma_{\beta}/|\mathcal{S}|)}}{e^{-\beta((E(s)-E(r))^++\gamma_{\beta}/|\mathcal{S}|)}}\ 
=\ e^{-\beta(E(r)-E(s)-2\gamma_{\beta}/|\mathcal{S}|)}.
\end{align*}
and
\begin{align*}
\frac{\pi(r)}{\pi(s)}\ 
&=\ \frac{p(s,r)}{p(r,s)}\ 
\geq\ \frac{e^{-\beta((E(r)-E(s))^++\gamma_{\beta}/|\mathcal{S}|)}}{e^{-\beta((E(s)-E(r))^+-\gamma_{\beta}/|\mathcal{S}|)}}\ 
=\ e^{-\beta(E(r)-E(s)+2\gamma_{\beta}/|\mathcal{S}|)}.
\end{align*}
Now let $r$ and $s$ be arbitrary. By the irreducibility, there is a path $r=r_0,r_1{,}...,r_n=s$ from $r$ to $s$ of neighboring states with $\pi(r_i)/\pi(r_{i+1})\in [e^{-\beta(E(r_i)-E(r_{i+1})+2\gamma_{\beta}/|\mathcal{S}|)},e^{-\beta(E(r_i)-E(r_{i+1})-2\gamma_{\beta}/|\mathcal{S}|)}],\, 0\leq i\leq n-1$. This finishes the proof.
\end{bew}

Under the stated assumptions, \textsc{Scoppola} \cite{Scop} has shown the existence of a successive filtration (aggregation) $\mathcal{S}=M^{(0)}\supset M^{(1)}\supset ... \supset M^{(\mathfrak{n})}=\{s_0\}$, $\mathfrak{n}\in \mathbb{N}$, of the state space such that the elements of each set $M^{(i)}, 1\leq i\leq \mathfrak{n}$, can be called \textit{metastable} in the following sense:
\begin{itemize}
\item They arise from the local minima of the energy function or certain modifications of it.
\item There is a lower bound on the expected time needed for a transition from $m_1$ to $m_2$ for any $m_1,m_2\in M^{(i)}$ which increases very fast with $i$.
\item There exists a constant $C$ such that
\[\mathbb{P}_m(X_t\notin M^{(i+1)})\leq e^{-C\beta}\]
for all $m\in M^{(i)}, 0\leq i\leq \mathfrak{n}-1,$ and sufficiently large $t$.
\end{itemize}
This filtration starts with
\[M^{(1)}:=\{s\in \mathcal{S}| E(s)<E(r) \textrm{ for all $r\sim s$}\},\]
and deletes one local minimum at each step. In fact, the local minimum with minimal activation energy for a transition to another minimum is deleted, see \cite{Scop} and \cite{Scop2} for further details.

\begin{figure}[ht]
\begin{center}
\includegraphics[trim=5.5cm 22.5cm 6cm 3.5cm, clip, width = 8cm]{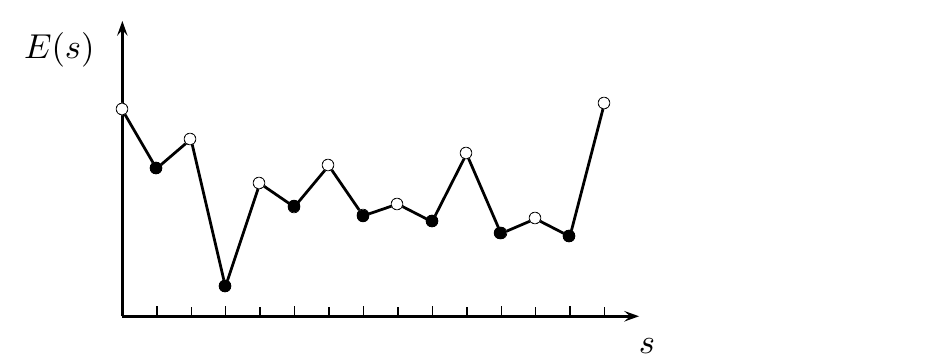}
\end{center}
\caption{Example of an energy landscape with minima shown as black dots ($\bullet$)}
\label{fig:PEL}
\end{figure}

\begin{figure}[htb]
\begin{center}
\includegraphics[trim=2cm 12cm 4.55cm 3.5cm, clip, width = 12cm]{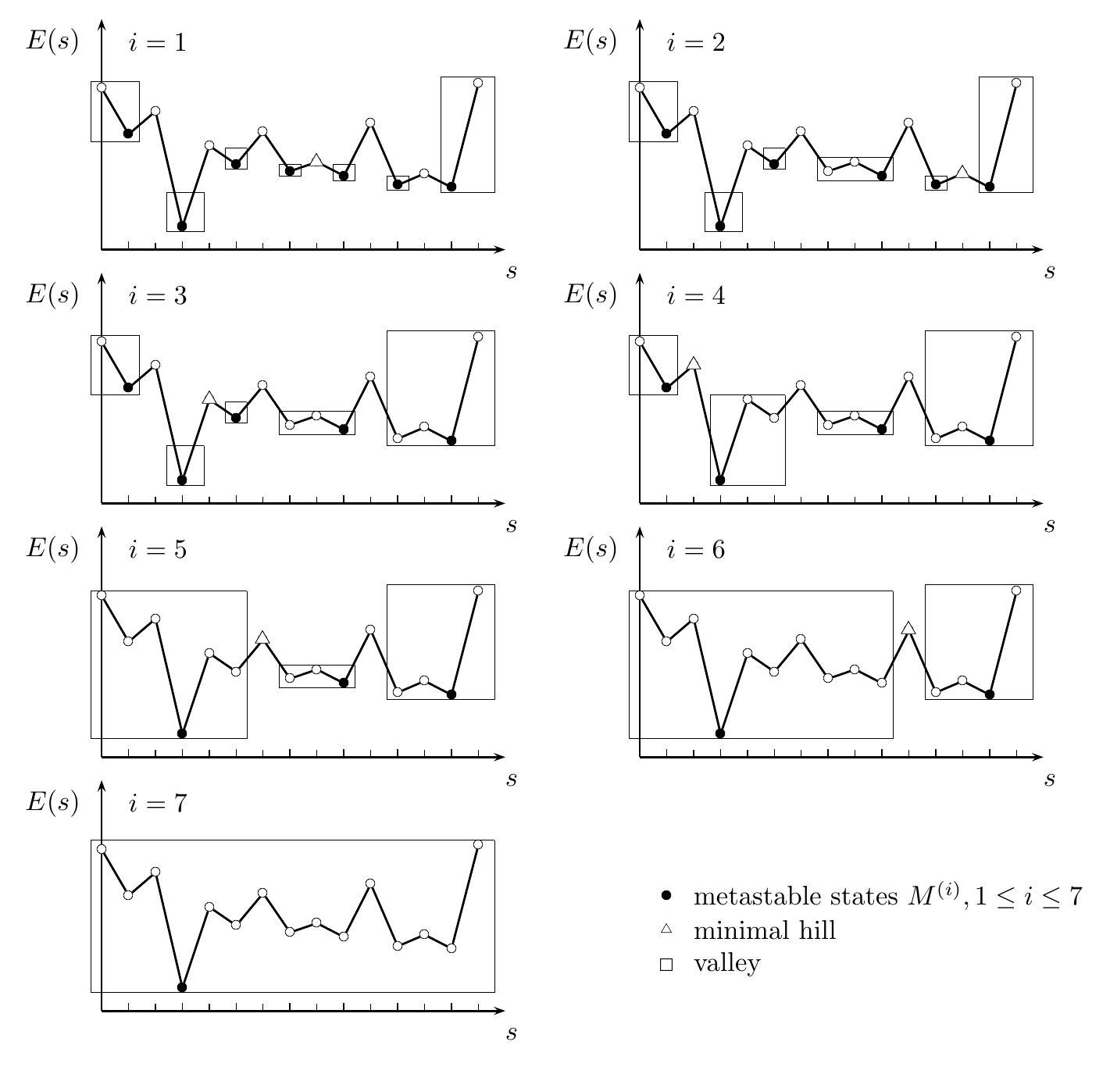}
\end{center}
\caption{Successive application of the algorithm in \protect\cite{Scop} to the energy landscape in Figure \protect\ref{fig:PEL}. For each step $i$, the metastable states as well as the corresponding valleys are shown.}
\label{fig:PELChange}
\end{figure}

\begin{bsp}\label{bsp:Energielandschaft}
For the simple energy function depicted in Figure \ref{fig:PEL}, a successive application of the algorithm from \cite{Scop} as illustrated in Figure \ref{fig:PELChange} leads to the following decomposition into subsets of metastable states:
\begin{align*}
M^{(1)}&=\{2,4,6,8,10,12,14\}, \ &M^{(2)}&=\{2,4,6,10,12,14\}, \\
M^{(3)}&=\{2,4,6,10,14\}, &M^{(4)}&=\{2,4,10,14\},\\
M^{(5)}&=\{4,10,14\}, \ &M^{(6)}&=\{4,14\},\\ 
M^{(7)}&=\{4\}.
\end{align*}
\end{bsp}

Based on the filtration of $\mathcal{S}$ just described, we now proceed to a definition of a sequence of metastable sets associated with the metastable states which will induce the MB. In order to do so, we must study first minimal paths between two states and maximal energies along such paths.

\begin{defn}\label{def:Pfade}
(a) For any two distinct states $r,s\in\mathcal{S}$, let 
\begin{align*}
\Gamma(r,s):=\{(x_0{,}...,x_k)|\,k\in \mathbb{N},\,x_0=r,\,x_k=s,\,
x_i\neq x_j,\, p(x_i,x_{i+1})>0 \textrm{ for } 0\leq i\leq k-1,\,i\neq j\}
\end{align*}
\begin{itemize}
\item[]
be the set of all finite self-avoiding paths from $r$ to $s$ having positive probability. For any such path $\gamma=(\gamma_0{,}...,\gamma_k)\in\Gamma(r,s)$, let $|\gamma|:=k$ be its length. We further write $t\in\gamma$ if $t\in\{\gamma_{1}{,}...,\gamma_{k}\}$.
\item[(b)] A self-avoiding path $\gamma=(\gamma_{1}{,}...,\gamma_{k})$ from $r$ to $s$ is called \textit{minimal} if its maximal energy $\max_{1\leq i\leq k} E(\gamma_i)$ is minimal among all $\gamma'\in\Gamma(r,s)$. The set of these paths is denoted $\Gamma^*(r,s)$.
\item[(c)] The \textit{essential saddle} $z^*(r,s)$ between $r$ and $s$ is then defined as
\[z^*(r,s):=\argmax_{t\in\gamma}E(t)\in\mathcal{S}\]
for any $\gamma\in\Gamma^{*}(r,s)$.
\end{itemize}
\end{defn}

As for (c), it is to be noted that, due to the assumed non-degeneracy of the energy function, the essential saddle is unique, which means that it does not depend on (as it must) which minimal path we choose in the definition of $z^*(r,s)$. There may indeed be several minimal paths, every single one thus crossing the saddle at some time. With the help of these notions the valleys can now be defined in a quite concrete way. Let us label the local minima as $m^{(1)}{,}...,m^{(\mathfrak{n})}$, so that $M^{(i)}=\{m^{(i)}{,}...,m^{(\mathfrak{n})}\}$ for each $i=1{,}...,\mathfrak{n}$.

\begin{defn}\label{def:Tal1}
For each $m \in M^{(i)}$, $1\le i\le\mathfrak{n}$, let
\begin{equation*}
V^{(i)}_<(m):=\left\{s\in\mathcal{S}\Big|E(z^*(s,m))< E(z^*(s,m'))\text{ for all }m'\in M^{(i)}\backslash \{m\}\right\}.
\end{equation*}
We say that \emph{state $s$ is attracted by $m$ at level $i$}, expressed as $s\leadsto m$ at level $i$, if
\[ E(z^*(s,m))=\min_{n\in M^{(i)}}E(z^*(s,n))\]
and every minimal path from $s$ to a state $m'\in M^{(i)}\backslash\{m\}$ with $E(z^*(s,m'))=E(z^*(s,m))$ hits $V^{(i)}_{<}(m)$ at some time. Finally, let
\[l(i):=\inf\big\{i<j\le\mathfrak{n}|m^{(i)}\leadsto m\text{ at level }j\text{ for some }m\in M^{(j)}\big\}\]
denote the minimal level at which the minimal state $m^{(i)}$ becomes attracted by a minimal state of superior level.
\end{defn}

\begin{defn}\label{def:Tal}
(a) \emph{Initialization}: For each $m \in M^{(1)}$, define
\begin{equation*}
V^{(1)}(m):=\left\{s\in\mathcal{S}\,\Big|\,s\leadsto m\text{ at level }1\right\}.
\end{equation*}
as the \textit{valley of order $1$ containing $m$} and let
\begin{equation*}
N^{(1)}:=\left(\bigcup_{j=1}^{\mathfrak{n}}V^{(1)}(m^{(j)})\right)^{c}
\end{equation*}
be the set of non-assigned states at level 1.

\vspace{.2cm}
(b) \emph{Recursion:} For $2\le i\le\mathfrak{n}$ and $m\in M^{(i)}$, define 
\begin{equation*}
V^{(i)}(m):=V^{(i-1)}(m)\,\cup\,\left\{s\in N^{(i-1)}\,\Big|\,s\leadsto m\text{ at level }i\right\}\,\cup\,\bigcup_{j:l(j)=i,m^{(j)}\leadsto m\text{ at level }i}\hspace{-6pt}V^{(j)}(m^{(j)})
\end{equation*}
as the  \textit{valley of order $i$ containing $m$} and let
\begin{equation*}
N^{(i)}:=\left(\bigcup_{j=1}^{\mathfrak{n}}V^{(i\wedge j)}(m^{(j)})\right)^{c}
\end{equation*}
be the set of non-assigned states at level i.
\end{defn}

Here is a more intuitive description of what the previous two definitions render in a formal way: First, we define, for each level $i$ and $m\in M^{(i)}$, the set $V^{(i)}_{<}(m)$ of those states $s$ that are strongly attracted by $m$ in the sense that $E(z^{*}(s,m))$ is strictly smaller than $E(z^{*}(s,m'))$ for any other $m'\in M^{(i)}$. Then, starting at level one, each valley $V^{(1)}(m)$, $m\in M^{(1)}$, is formed from $V^{(1)}_{<}(m)$ by adjoining all further states $s$ attracted by $m$ at this level. This leaves us with a set of non-assigned states, denoted $N^{(1)}$. In the next step (level 2), any $V^{(2)}(m)$ for $m\in M^{(2)}$ is obtained by adjoining to $V^{(1)}(m)$ all those $s\in N^{(1)}$ which are attracted by $m$ at level 2. Observe that this ensures $V^{(2)}_{<}(m)\subset V^{(2)}(m)$. Moreover, if $m^{(1)}$ is attracted by $m$ at level 2, then $V^{(1)}(m^{(1)})$ is merged into $V^{(2)}(m)$ as well. If no such $m$ exists (thus $l(1)>2$), it remains untouched until reaching level $l(1)$ where 
its bottom state $m^{(1)}$ becomes attracted by some $m'\in M^{(l(1))}$ causing its valley to be merged into $V^{(l(1))}(m')$. This procedure continues in the now obvious recursive manner until at level $\mathfrak{n}$ all states have been merged into one valley. Obviously, valleys of the same order are pairwise disjoint. Also, valleys once formed at some level can only be merged as a whole and will thus never be ripped apart during the recursive construction. For the energy function depicted in Figure \ref{fig:PEL}, the successively derived valleys of order $i=1{,}...,7$ are shown in Figure \ref{fig:PELChange}.

\vspace{.1cm}
Before proceeding to results on the general shape of valleys, we collect some basic, mostly technical properties of essential saddles which will be useful thereafter.

\begin{prop}\label{prop:essentiellerSattel}
For any $r,s,u\in \mathcal{S}$, $0\leq i\leq \mathfrak{n}$, $m_1,m_2\in M^{(i)}, m_1\neq m_2,$ and $x_1,x_2\in \mathcal{S}$ with $x_1\in V^{(i)}_<(m_1)$ and $x_2\in V^{(i)}(m_2)$, we have
\begin{itemize}
\item[(a)] $z^*(r,s)=z^*(s,r)$.
\item[(b)] $E(z^*(r,s))\leq E(z^*(r,u))\vee E(z^*(u,s))$.
\item[(c)] $E(z^*(x_2,m_2))\leq E(z^*(x_2,m'))$ for all $m'\in M^{(i)}$.
\item[(d)] $E(z^*(x_1,m_2))=E(z^*(m_1,m_2))$.
\item[(e)] $E(z^*(x_1,x_2))\geq E(z^*(m_1,m_2))$.
\item[(f)] $z^*(x_1,x_2)\neq x_1$.
\end{itemize}
\end{prop}

\begin{bew}
Parts (a) and (b) are obvious. 

For (c) we use an induction over $i$ and note that there is nothing to show when $i=1$. For general $i$, we must only verify that $E(z^{*}(x_{2},m_{2}))\le E(z^{*}(x_{2},m'))$ for all $m'\in M^{(i)}$ if $x_{2}\in V^{(j)}(m^{(j)})$ for some $j<i$ such that $l(j)=i$ and $m^{(j)}\leadsto m$ at level $i$ (due to the recursive definition of $V^{(i)}(m)$). But the latter ensures that $E(z^{*}(x_{2},m^{(j)}))\le E(z^{*}(x_{2},n))$ for all $n\in M^{(j)}\supset M^{(i)}$ (inductive hypothesis) as well as $E(z^{*}(m^{(j)},m_{2}))\le E(z^{*}(m^{(j)},m'))$ for all $m'\in M^{(i)}$. Consequently, for any such $m'$,
\begin{align*}
E(z^{*}(x_2,m_2))&\le E(z^{*}(x_{2},m^{(j)}))\vee E(z^{*}(m^{(j)},m_{2}))\\
&\le E(z^{*}(x_{2},m^{(j)}))\vee E(z^{*}(m^{(j)},m'))\\
&\le E(z^{*}(x_{2},m^{(j)}))\vee E(z^{*}(x_{2},m')))\\
&=E(z^{*}(x_{2},m'))
\end{align*}
as asserted.

For assertion (d), note that $E(z^*(x_1,m_2))>E(z^*(x_1,m_1))$, which in combination with (a) and (b) implies
\begin{equation*}
E(z^*(m_1,m_2))\leq E(z^*(m_1,x_1))\vee E(z^*(x_1,m_2))=E(z^*(x_1,m_2))
\end{equation*}
and then further
\[E(z^*(x_1,m_2))\leq \underbrace{E(z^*(x_1,m_1))}_{<E(z^*(x_1,m_2))}\vee \underbrace{E(z^*(m_1,m_2))}_{\leq E(z^*(x_1,m_2))}\leq E(z^*(x_1,m_2)).\]
So the above must be an identity, i.e.\ $E(z^*(x_1,m_2))=E(z^*(m_1,m_2))$.

Turning to part (e), we first infer with the help of (c) and (d) that
\begin{align}\label{eqn:minBerge0}
E(z^*(x_1,m_1))&< E(z^*(x_1,m_2)) \nonumber\\
&=E(z^*(m_1,m_2)) \nonumber\\
&\leq E(z^*(m_1,x_2))\vee E(z^*(x_2,m_2))\\
&= E(z^*(x_2,m_1)) \nonumber\\
&\leq E(z^*(x_2,x_1))\vee E(z^*(x_1,m_1)), \nonumber
\end{align}
thus
\begin{equation}\label{eqn:minBerge}
E(z^*(x_1,m_1))<E(z^*(x_1,x_2)).
\end{equation}
Together with the just shown inequality $E(z^*(m_1,m_2))\leq E(z^*(x_2,m_1))$ (see \eqref{eqn:minBerge0}) and another use of (c), this yields
\[E(z^*(x_1,x_2))=E(z^*(x_2,x_1))\vee E(z^*(x_1,m_1))\geq E(z^*(x_2,m_1))\geq E(z^*(m_1,m_2)).\]

Finally, we infer with the help of \eqref{eqn:minBerge} that
\begin{align*}
E(z^*(x_1,x_2))&>E(z^*(x_1,m_1))\geq E(x_1)
\end{align*}
and thus $z^*(x_1,x_2)\neq x_1$ as claimed in (f).
\end{bew}

\begin{bem}\label{bem:minimale Pfade}
It is useful to point out the following consequence of the previous proposition. If, for an arbitrary state $s$ and any two distinct metastable states $m,n\in M^{(i)}$, there exists a minimal path $\gamma$ from $s$ to $n$ that hits a state $r$ with $E(z^*(r,m))<E(z^*(r,n)$, then there is also a minimal path from $s$ to $n$ that passes through $m$. Namely, if we replace the segment from $r$ to $n$ of the former path by the concatenation of two minimal paths from $r$ to $m$ and from $m$ to $n$, then the maximal energy of this new path is 
\begin{align*}
E(z^{*}(s,n))\vee E(z^{*}(r,m))\vee E(z^{*}(m,n))&\leq E(z^{*}(s,n))\vee E(z^{*}(r,m))\vee E(z^*(r,n))\\
&=E(z^{*}(s,n))\vee E(z^{*}(r,n))\\
&=E(z^*(s,n)),
\end{align*}
by Proposition \ref{prop:essentiellerSattel}(b), whence the new path has to be minimal from $s$ to $n$ as well. This yields two facts:
\begin{itemize}
\item[(a)] A minimal path from $s$ to $n$, where $s\leadsto n$ at level $i$, hits $V^{(i)}_<(n)$ before it hits any $r$ with $E(z^*(r,m))<E(z^*(r,n))$ for some $m\in M^{(i)}$. Otherwise, since the subpath from $r$ to $m$ can be chosen to stay in $\{t|E(z^*(t,m))<E(z^*(t,n))\}$ and thus $E(z^*(s,m))=E(z^*(s,n))$, there would be a path from $s$ to $m$ not hitting $V^{(i)}_<(n)$.
\item[(b)] If $s\leadsto n$ at level $i$ and $m\in M^{(i)}\backslash \{n\}$ with $E(z^*(s,n))=E(z^*(s,m))$, then a minimal path from $s$ to $m$ does not only hit $V^{(i)}_<(n)$ at some time, but in fact earlier than any other valley $V^{(i)}_<(m'), m'\in M^{(i)}\backslash \{n\}$.
\end{itemize}
\end{bem}

\begin{lem}\label{lem:no inner points delayed}
Let $1\le i<j\le\mathfrak{n}$, $m=m^{(i)}$ and $s\in V^{(i)}(m)$. Then $s\in V^{(j)}_{<}(m')$ for some $m'\in M^{(j)}$ implies $l(i)\le j$, $m\in V^{(j)}_{<}(m')$ and thus $V^{(i)}(m)\subset V^{(j)}(m')$.
\end{lem}

In other words, whenever $V^{(i)}(m^{(i)})$ contains an element $s$ which at some higher level $j$ belongs to some $V^{(j)}_{<}(m')$, $m'\in M^{(j)}$, the same must hold true for $m^{(i)}$ itself implying $V^{(i)}(m^{(i)})\subset V^{(j)}(m')$. Conversely, this guarantees that $V^{(i)}(m^{(i)})$ will have no common elements with any $V^{(j)}_{<}(m')$ at levels $j<l(i)$ where it has not yet been merged into a valley of higher order.

\begin{bew}
Let us first note that, under the given assumptions,
\begin{equation*}
E(z^{*}(s,m))\le E(z^{*}(s,m'))<E(z^{*}(s,n))
\end{equation*}
for all $n\in M^{(j)}\backslash\{m'\}$, whence
\[E(z^{*}(m,n))\le E(z^{*}(s,m))\vee E(z^{*}(s,n))=E(z^{*}(s,n))\le E(z^{*}(s,m))\vee E(z^{*}(m,n))\]
entails $E(z^{*}(m,n))=E(z^{*}(s,n))$ for all such $n$. Using this fact, we find that
\begin{align*}
E(z^{*}(m,m'))&\le E(z^{*}(s,m))\vee E(z^{*}(s,m'))<E(z^{*}(s,n))=E(z^{*}(m,n))
\end{align*}
for all $n\in M^{(j)}\backslash\{m'\}$, which implies $m\leadsto m'$ at level $j$ and thus $l(i)\le j$ as well as the other assertions.
\end{bew}

\begin{prop}
For every $m\in M^{(i)}$ and $1\leq i \leq \mathfrak{n}$, $V_{<}^{(i)}(m)$ is connected.
\end{prop}

\begin{bew}
Pick any $s\in V^{(i)}_<(m)$, any minimal path from $s$ to $m$ and finally any intermediate state $r$ along this path for which $r\in V^{(i)}_<(m)$ must be verified. For every $m'\in M^{(i)}\backslash \{m\}$, we find
\begin{align*}
E(z^*(r,m))&\leq E(z^*(r,s))\vee E(z^*(s,m))\\
&=E(z^*(s,m))\\
&<E(z^*(s,m'))\\
&\leq \underbrace{E(z^*(s,r))}_{<E(z^*(s,m'))}\vee E(z^*(r,m'))\\
&=E(z^*(r,m')),
\end{align*}
which shows $r\in V^{(i)}_<(m)$ as required.
\end{bew}

Note that we have even shown that a minimal path from a state in $V^{(i)}_<(m)$ to $m$ will never leave this set. We may expect and will indeed show as Proposition \ref{prop:Zusammenhängend} below that $V^{(i)}(m)$ is connected as well. The following lemma is needed for its proof.

\begin{lem}\label{lem:connected direct path}
Given $1\le i\le\mathfrak{n},\,m\in M^{(i)}$ and $s\leadsto m$ at level $i$, let $\gamma=(\gamma_{1}{,}...,\gamma_{k})\in\Gamma^{*}(s,m)$ be a path such that $E(z^*(\gamma_i,m))\leq E(z^*(\gamma_i,n))$ for all $n\in M^{(i)}\backslash\{m\}$, and which stays in $V_{<}^{(i)}(m)$ once hitting this set (such a $\gamma$ exists by Remark \ref{bem:minimale Pfade} (a)). Then $\gamma_{j}\leadsto m$ at level $i$ for each $j=1{,}...,k$.
\end{lem}

\begin{bew}
There is nothing to prove for $\gamma_{1}=s$ and any $\gamma_{j}\in V_{<}^{(i)}(m)$. So let $r$ be any other state visited by $\gamma$, pick an arbitrary $n\in M^{(i)}\backslash \{m\}$ with $E(z^*(r,n))=E(z^*(r,m))$ and then any minimal path $\tau$ from $r$ to $n$. Let $\sigma$ be the subpath of $\gamma$ from $s$ to $r$. We must show that $\tau$ hits $V^{(i)}_<(m)$. First, we point out that the maximal energy $E(z^*(s,r))\vee E(z^*(r,n))$ of $\sigma \tau$, the concatenation of $\sigma$ and $\tau$, satisfies
\[E(z^*(s,n))\leq E(z^*(s,r))\vee E(z^*(r,n))\leq E(z^*(s,m))\vee E(z^*(r,m))=E(z^*(s,m))\leq E(z^*(s,n)),\]
implying $\sigma\tau\in\Gamma^{*}(s,n)$ and, furthermore,
\[E(z^*(s,r))\vee E(z^*(r,n))=E(z^*(s,m))\vee E(z^*(r,n))=E(z^*(s,m))\vee E(z^*(r,m))=E(z^*(s,m)).\] 
Thus $\sigma\tau$ must hit $V^{(i)}_<(m)$. But since $\sigma$ does not hit $V_{<}^{(i)}(m)$ by assumption, we conclude that $\tau$ must hit $V_{<}^{(i)}(m)$. Since $\tau\in\Gamma^{*}(r,n)$ was arbitrary, we infer $r\leadsto m$ at level $i$.
\end{bew}

The next two propositions provide information on the shape of the valleys and their nested structure.

\begin{prop}\label{prop:Zusammenhängend}
For every $m\in M^{(i)}$ and $1\leq i \leq \mathfrak{n}$, $V^{(i)}(m)$ is connected.
\end{prop}

\begin{bew}
We use an inductive argument. If $i=1$, the assertion follows directly from the definition of the level-one valleys because any $s\in V^{(1)}(m)$, $m\in M^{(1)}$, may be connected to $m$ by a minimal path that eventually enters $V_{<}^{(1)}(m)$ without hitting any other $V_{<}^{(1)}(n)$ and is therefore completely contained in $V^{(1)}(m)$ by the previous lemma.

Turning to the inductive step, suppose the assertion holds true up to level $i-1$. Fix any $m\in M^{(i)}$ and notice that, by the inductive hypothesis, $V^{(i-1)}(m)$ as well as all $V^{(j)}(m^{(j)})$ with $l(j)=i$ and $m^{(j)}\leadsto m$ at level $i$ are connected. Now, since these $m^{(j)}$ as well as all $s\in N^{(i-1)}$ attracted by $m$ at level $i$ may be connected to $m$ by minimal paths as assumed in Lemma \ref{lem:connected direct path}, we conclude that $V^{(i)}(m)$ is also connected.
\end{bew}

The second proposition shows the nested structure of our construction of valleys.

\begin{prop}\label{prop:TaelerRekursion}
The following inclusions hold true:
\begin{enumerate}
\item[(a)] $V^{(1)}(m)\subseteq...\subseteq V^{(i)}(m)$ for each $m\in M^{(i)},\,1\leq i\leq \mathfrak{n}$.
\item[(b)] $V^{(i)}(m)\subseteq V^{(j)}(n)$ for each $1\le i<j\le\mathfrak{n}$, $n\in M^{(j)}$ and $m\in M^{(i)}\cap V^{(j)}(n)$.
\end{enumerate}
\end{prop}

\begin{bew}
Since there is nothing to show for (a) we move directly to (b). But if $m\in M^{(i)}\cap V^{(j)}(n)$, then the definition of valleys ensures the existence of $1\le k\le j-i$ and of $n_{1}{,}...,n_{k-1}\in M^{(j)}\backslash M^{(i)}$ such that $n_{p-1}\leadsto n_{p}$ at level $l_{p}$ for each $p=1{,}...,k$ and levels $i<l_{1}<...<l_{k}=j$, where $n_{0}:=m$ and $n_{k}:=n$. As a consequence,
\[V^{(i)}(m)\subseteq V^{(l_{1})}(n_{1})\subseteq...\subseteq V^{(l_{k-1})}(n_{k-1})\subseteq V^{(j)}(n)\]
which proves the asserted inclusion.
\end{bew}

To finish the analysis of the shape of the valleys we show that they have the following important property: a special class of minimal paths from the inside of any $V^{(i)}(m)$ to the outside of it must hit its interior $V^{(i)}_<(m)$. But in order to show this we must first verify that all states attracted by $m$ at level $i$ belong to $V^{(i)}(m)$.

\begin{lem}\label{lem:valleys and attraction}
For each $1\le i\le\mathfrak{n}$ and $m\in M^{(i)}$, we have that
\[\left\{s\in \mathcal{S}\Big| s\leadsto m \textrm{ at level i}\right\}\subset V^{(i)}(m)\subset\left\{s\in\mathcal{S}\Big|E(z^{*}(s,m))\le E(z^{*}(s,m'))\text{ for all }m'\in M^{(i)}\right\}.\]
\end{lem}

\begin{bew}
For the second inclusion it suffices to refer to Proposition \ref{prop:essentiellerSattel}(c). The first inclusion being obviously true for $s\in N^{(i-1)}$, we turn directly to the case when
\[s\leadsto n_1 \textrm{ at level } l_1, \quad n_1\leadsto n_2 \textrm{ at level }l_2, \quad ... \quad n_{k-1}\leadsto n_k \textrm{ at level }l_k\]
with $k\geq 1$ and $1\leq l_1\leq ... \leq l_k\leq i-1$. Here, $n_1$ denotes the first minimum to which $s$ is attracted (thus $s\in V^{(l_{1})}(n_{1})$), while $n_k$ is the last minimum of this kind in the sequence. We may assume without loss of generality that $n_j\neq m$ for all $j$, for otherwise the assertion is clear. 

We show now that $n_1\leadsto m$ at level $i$ which in turn implies $n_j\leadsto m$ at level $i$ for all $1\leq j\leq k$. As a consequence, $n_k\in V^{(i)}(m)$ and thus $s\in V^{(i)}(m)$. If $E(z^*(n_1,m))<E(z^*(n_1,m'))$ for all $m'\in M^{(i)}\backslash \{m\}$, the assertion is proved. Hence suppose $E(z^*(n_1,m))\geq E(z^*(n_1,m'))$ for some $m'\in M^{(i)}\backslash \{m\}$. Then
\begin{align*}
E(z^*(s,m'))&\leq E(z^*(s,n_1))\vee E(z^*(n_1,m'))\\
&\leq E(z^*(s,n_1))\vee E(z^*(n_1,m))\\
&\leq E(z^*(s,n_1))\vee E(z^*(s,m))\\
&=E(z^*(s,m))\\
&\leq E(z^*(s,m')),
\end{align*}
implies $E(z^*(s,m))=E(z^*(s,m'))$ and also that the concatenation of any minimal path $\gamma$ from $s$ to $n_1$ and any minimal path $\tau$ from $n_1$ to $m'$ (with maximal energy $E(z^*(s,n_1))\vee E(z^*(n_1,m'))$) constitutes a minimal path from $s$ to $m'$ and must therefore hit $V^{(i)}_<(m)$. Note that we can choose $\gamma$ to stay in $V^{(l_{1})}(n_1)$ since $s\in V^{(l_{1})}(n_{1})$ and $V^{(l_{1})}(n_{1})$ is connected. Now, if $\tau$ hits $V^{(i)}_<(m)$, then $E(z^*(n_1,m))=E(z^*(n_1,m'))$ and we are done. Otherwise, $\gamma$ hits $V^{(i)}_<(m)$ implying $V^{(j)}(n_1)\cap V^{(i)}_<(m)\neq \emptyset$. Now use Lemma \ref{lem:no inner points delayed} to conclude $n_1\in V^{(i)}_<(m)$ and therefore $n_1\leadsto m$ at level $i$. This completes the argument for the first inclusion.
\end{bew}

We provide too further lemmata that will be needed later on.

\begin{lem}\label{lem:Anziehung}
Let $m\in M^{(i)}, x\leadsto m$ at level $i$ and $y\notin V^{(i)}(m)$. Then either every minimal path from $x$ to $y$ hits the set $V^{(i)}_<(m)$, or $E(z^*(x,y))>E(z^*(x,m))$. 
\end{lem}

\begin{bew}
Suppose there is a minimal path $\gamma$ from $x$ to $y$ avoiding $V^{(i)}_<(m)$. Since $y\notin V^{(i)}(m)$, it is not attracted by $m$ at level $i$ implying the existence of some $m'\in M^{(i)}$ with $E(z^*(y,m'))\leq E(z^*(y,m))$ and of some $\tau\in \Gamma^*(y,m')$ avoiding $V^{(i)}_<(m)$. Hence, the concatenation $\gamma\tau$ avoids $V^{(i)}_<(m)$ and must therefore have maximal energy larger than $E(z^*(x,m))$. Consequently,
\begin{align*}
E(z^*(x,m))
&<E(z^*(x,y))\vee E(z^*(y,m'))\\
&\leq E(z^*(x,y))\vee E(z^*(y,m))\\
&\leq E(z^*(x,y))\vee E(z^*(x,m)),
\end{align*}
and thus $E(z^*(x,y))>E(z^*(x,m))$.
\end{bew}

In order to state the second lemma, let us define the outer part $\partial^{+}V$ of a valley $V$ to be the set of those states outside $V$ which are adjacent to a state in $V$. With the help of the previous result, we can easily show that $\partial^{+}V$ contains only non-assigned states at any level where $V$ has not yet been merged into a larger valley.

\begin{lem}\label{lem:outer boundary nonassigned}
For any $1\leq i, j\leq \mathfrak{n}$ and $m=m^{(i)}$ with $l(i)> j$, the outer part $\partial^{+}V$ of the valley $V:=V^{(j\wedge i)}(m)$ is a subset of $N^{(j)}$ and $E(z^*(s,m))=E(s)$ for every $s\in \partial^+V$.
\end{lem}

\begin{bew}
First, let $s\in\partial^{+}V$ and suppose that $s\notin N^{(j)}$. Then $s\leadsto m'$ at level $k$, in particular $s\in V^{(k)}(m')$ for some $m'\in M^{(k)}$ and $k\le j$. Pick any $r\in V$ with $r\sim s$ and note that $r\in\partial^{+}V^{(k)}(m')$. W.l.o.g.\ we may assume that $r\leadsto m$ at level $j\wedge i$. Then Lemma \ref{lem:Anziehung} (with $x=r$ and $y=s$) ensures that either $E(z^{*}(r,s))>E(z^{*}(r,m)\ge E(r)$, thus $z^{*}(r,s)=s$ and $E(r)<E(s)$, or $r\in V_<(m)$ and, for some $n\in M^{(j\wedge i)}$,
\begin{equation}\label{eqn:outer boundary nonassigned}
\begin{split}
E(z^*(r,m))
&<E(z^*(r,n))\\
&\leq E(z^*(r,s))\vee E(z^*(s,n))\\
&\leq E(z^*(r,s))\vee E(z^*(s,m))\\
&\leq E(z^*(r,s))\vee E(z^*(r,m))\\
&=E(s)\vee E(z^*(r,m)),
\end{split}
\end{equation}
and thus again $E(r)<E(s)$. On the other hand, by the very the same lemma (now with $x=s$ and $y=r$), we infer $E(r)>E(s)$ which is clearly impossible. Consequently, $s$ must be non-assigned at level $j$ as claimed.

For the second assertion take again $s\in \partial^+V$ and a minimal path $\gamma=(s,r{,}...,m)$ from $s$ to $m$ with $r\in V$. Again, by use of Lemma \ref{lem:Anziehung}, we find either $E(z^*(r,m))<E(z^*(r,s))=E(s)$ or $r\in V_<(m)$, leading analogously to equation \eqref{eqn:outer boundary nonassigned} to $E(z^*(r,m))< E(s)\vee E(z^*(r,m))$ and thus $E(z^*(r,m))<E(s)$. In conclusion, both cases result in 
\[E(z^*(s,m))=E(s)\vee E(z^*(r,m))=E(s),\]
finishing the proof.
\end{bew}

The reader may wonder why valleys are defined here via essential saddles and not via the at first glance more natural overall energy barriers, viz.
\begin{equation}\label{eqn:I}
I(s,m):=\inf_{\gamma\in\Gamma(s,m)}I(\gamma_1{,}...,\gamma_{|\gamma|})
\end{equation}
with
\[I(s_1{,}...,s_n):=\sum_{i=1}^{n}(E(s_i)-E(s_{i-1}))^+\]
for a state $s$ in a valley and the pertinent minimum $m$. This latter quantity, also called \emph{activation energy}, is indeed an important parameter in \cite{Scop}. The reason for our definition becomes clear when regarding the last proposition which shows the nested structure of valleys of increasing order and which may fail to hold when choosing an alternative definition based on the activation energy. Valleys that are formed in one step may then be ripped apart in the next one. This happens, for instance, if there is just one large saddle along the path to a metastable state and several small ones, lower than the essential saddle, along the paths to another minimum such that their total sum is larger than the big saddle. In further support of our approach, it will be seen later that the essential saddles are the critical parameters for the behavior of the aggregated chain (see Theorem \ref{prop:AbschaetzungUeWkeiten}).

In a nutshell, by going from $(V^{(i)}(m))_{m\in M^{(i)}}$ to $(V^{(i+1)}(m))_{m\in M^{(i+1)}}$, some valleys are merged into one (with only the smaller minima retained as metastable states) and additionally those states from $N^{(i)}$ are added which at level $i$ were attracted by metastable states now all belonging to the same valley. This induces the following tree-structure on the state space:
\begin{itemize}
\item Fix $\varnothing=s_0$.
\item The first generation of the tree consists of all $m\in M^{(\mathfrak{n}-1)}\cup N^{(\mathfrak{n}-1)}$ and are thus connected to the root.
\item The second generation of the tree consists of all $m\in M^{(\mathfrak{n}-2)}\cup N^{(\mathfrak{n}-2)}$, and $m$ is connected to the node $k$ of the first generation for which $E(z^*(m,k))$ is minimal or to itself (in the obvious sense).
\item This continues until in the $\mathfrak{n}^{\textrm{th}}$ generation each state is listed and connected either with its unique point of attraction in the previous generation or with itself.
\end{itemize}

\begin{bsp}
For the energy function of Example \ref{bsp:Energielandschaft} and depicted in Figure \ref{fig:PEL}, the described tree is shown in Figure \ref{fig:Baum}. 
\begin{figure}[htb]
\begin{center}
\includegraphics[trim=5.5cm 22cm 8cm 3cm, clip, width = 6cm]{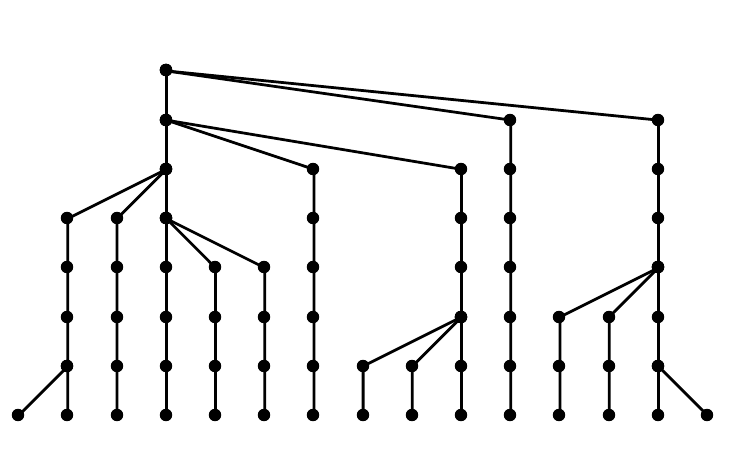}
\end{center}
\caption{The tree belonging to Figure \protect\ref{fig:PEL}}
\label{fig:Baum}
\end{figure}
The sets of non-assigned states at the different levels are
\begin{align*}
N^{(1)}&=\{3,5,7,9,11,13\}, \ &N^{(2)}&=\{3,5,7,11,13\}, &\\
N^{(3)}&=\{3,5,7,11\}, &N^{(4)}&=\{3,7,11\}, &\\ 
N^{(5)}&=\{7,11\}, \ &N^{(6)}&=\{11\}, \\
N^{(7)}&=\emptyset. &&
\end{align*} 
\end{bsp}

At each level $i$ of such a tree the subtree rooted at a node $m\in M^{(i)}$ consists of the states in the valley $V^{(i)}(m)$. A similar graph-theoretical modeling in order to visualize high dimensional energy landscapes has been used, for example, by {\sc Okushima} et~al. in \cite{OkNiIkSh}. These authors work with saddles of paths as well. In contrast to our approach, every possible path, that is every possible saddle is represented as a node in the tree. But as we will see, in the limit of low temperatures $(\beta\to \infty)$ the essential saddle is all we need.

\vspace{.2cm}
Now there are two fundamental directions for further investigations:

\begin{itemize}
\item[(1)] \textsc{Microscopic View}: 
What happens while the process visits a fixed valley $V$ (of arbitrary level)? \\ 
In Section \ref{sec:MikroskopischerProzess}, we will show that during each visit of a valley $V$ its minimum will be reached with probability tending to 1 as $\beta\to \infty$, and we also calculate the expected residence time in $V$, establish Property 3 stated in the Introduction and comment briefly on quasi-stationarity.
\item[(2)] \textsc{Macroscopic View}: 
How does the process jump between the valleys? \\ 
In Section \ref{sec:MakroskopischerProzess}, by drawing on the results of Section \ref{sec:MikroskopischerProzess}, we will show that an appropriate aggregated chain is Markovian in the limit as $\beta\to \infty$ and calculate its transition probabilities. With this we will finally be able to provide the definition of MB and establish Properties 1, 2 and 4 listed in the Introduction.
\end{itemize}

%
%

\section{Microscopic View: Fixing a Valley}\label{sec:MikroskopischerProzess}

Based on the provided definition of valleys of different orders, we are now going to study the process when moving in a fixed valley. 

\subsection{Trajectories for \texorpdfstring{$\beta\to\infty$}{ }}

The first goal in our study of the microscopic process and also the basic result for the subsequent analysis of the macroscopic process deals with the probabilities of reaching certain states earlier than others. From this we will conclude that in the limit $\beta\to \infty$ (which is the low temperature limit in the Metropolis Algorithm) the process, when starting somewhere in a valley, will visit its minimum before leaving it. 

\vspace{.1cm}
For $A\subset \mathcal{S}$ and $x\in \mathcal{S}$, we define 
\[\tau_A:=\inf\{n\geq 1|X_n\in A\},\quad\tau_x:=\tau_{\{x\}}\quad\text{and}\quad\mathcal{N}(x):=\{y\in \mathcal{S}|p(x,y)>0\}.\]
Two states $x,y$ with $p(x,y)>0$ are called neighbors ($x\sim y$) and $\mathcal{N}(x)$ the neighborhood of $x$. Hence, $\mathcal{S}$ may (and will) be viewed as a graph hereafter with edge set $\{(x,y)|x\sim y\}$. Given any subgraph $\Delta$, we will write $\widetilde{\mathbf{P}}$ for the transition matrix of the chain restricted to $\Delta$ ($\tilde{p}(r,s)=p(r,s)$ for all distinct $r,s\in \Delta$, $\tilde{p}(r,r)=1-\sum_{s\in \Delta, s\neq r}\tilde{p}(r,s)$) and $\widetilde{\mathbb{P}}_{x}$ for probabilities when regarding this restricted chain starting at $x\in\Delta$.

\subsubsection{Hitting probabilities}

\begin{thm}\label{thm:DriftZumMinimum}
Let $x, y, z\in \mathcal{S}$ be any pairwise distinct states satisfying $E(z^*(x,z))>E(z^*(x,y))$ and $z^*(x,z)\neq x$. Then there exist nonnegative constants $K(\beta)$ satisfying $\sup_{\beta>0}K(\beta)<\infty$ (and given more explicitly in Proposition \ref{prop:AbschaetzungUeWkeiten} below) such that
\[\mathbb{P}_{x}(\tau_{z}<\tau_{y})\leq K(\beta)\,e^{-\beta(E(z^*(x,z))-E(z^*(x,y))-7\gamma_{\beta})}=:\varepsilon(x,y,z,\beta) \overset{\beta\to\infty}{\longrightarrow} 0.\]
\end{thm}

Thus in the limit of low temperatures $(\beta\to\infty)$, only the smallest of all possible energy barriers affects the speed of a transition. In particular, we have the following result which is preliminary to the subsequent one.

\begin{thm}\label{thm:DriftZumMinimumTäler}
Given distinct $x,y\in\mathcal{S}$ and $m\in M^{(i)}$ such that $x\leadsto m$ at level $i$ and $y\notin V^{(i)}(m)$, let $B:=\{z\in \mathcal{S}|E(z^*(x,z))>E(z^*(x,m))\}$. Then it holds true that
\begin{align*}
\mathbb{P}_{x}&(\tau_{y}<\tau_{m})\\
&\le\ \varepsilon(x,m,y,\beta)\,\mathds{1}_{B}(y)+\left(\sum_{z:E(z)>E(z^*(x,y))}\varepsilon(x,y,z,\beta)+\sum_{z\in V^{(i)}_<(m)}\varepsilon(z,m,y,\beta)\right)\,\mathds{1}_{B^{c}}(y)\\
&=:\ \tilde{\varepsilon}(x,m,y,\beta)\ \overset{\beta\to\infty}{\longrightarrow}\ 0.
\end{align*}
\end{thm}

\begin{thm}\label{thm:DriftZumMinimumTäler2}
Given $m\in M^{(i)}$, $x\in V^{(i)}(m)$ and $y\notin V^{(i)}(m)$, let $k\le i$ be such that
\begin{align*}
&m_{0}:=x\leadsto m_{1}\text{ at level }l_{1},\quad m_{1}\leadsto m_{2}\text{ at level }l_{2},\quad...\quad
m_{k-1}\leadsto m_{k}=m\text{ at level }l_{k}
\end{align*}
for suitable $1\le l_{1}<...<l_{k}\le i$, $m_{j}\in M^{(l_{j})}$ for $j=1{,}...,k$ determined by the construction in Definition \ref{def:Tal}. Then
\begin{align*}
\mathbb{P}_{x}(\tau_{y}<\tau_{m})\ &\leq\ \sum_{j=1}^{k}\mathbb{P}_{m_{j-1}}(\tau_{y}<\tau_{m_{j}})\leq\ \sum_{j=1}^{k}\tilde{\varepsilon}(m_{j-1},m_{j},y,\beta)
\ \overset{\beta\to\infty}{\longrightarrow}\ 0.
\end{align*}
\end{thm}

For the valleys as defined here this confirms Property 1 stated in the Introduction: If $\beta$ is sufficiently large, then with high probability the minimum of a valley is visited before this valley is left.

All three theorems are proved at end of the next subsection after a number of auxiliary results.

\subsubsection{Auxiliary results and proofs}

The proof of Theorem \ref{thm:DriftZumMinimum} will be accomplished by a combination of two propositions due to \textsc{Bovier} et al.~\cite{Bov} for a more special situation. We proceed with a reformulation of the first one in a weaker form and under weaker assumptions.

\begin{prop}[compare Theorem 1.8 in \cite{Bov}]\label{prop:AbschaetzungUeWkeiten}
Let $x, y, z\in\mathcal{S}$ be pairwise distinct such that $z^*(x,z)\neq x$. Then
\begin{align*}
\mathbb{P}_x(\tau_z<\tau_x)&\leq K(\beta)\,|\mathcal{S}|^{-1}\,e^{-\beta(E(z^*(x,z))-E(x)-2\gamma_{\beta})},\\
\mathbb{P}_x(\tau_y<\tau_x)&\geq |\mathcal{S}|^{-1}\,e^{-\beta (E(z^*(x,y))-E(x)+5\gamma_{\beta})},
\end{align*}
where $K(\beta):=|\mathcal{S}|\max_{r\in\mathcal{S}}|\mathcal{N}(r)|\left(|\mathcal{S}|e^{-\beta(\min_{a\neq b:E(a)>E(b)}(E(a)-E(b))-2\gamma_{\beta})}+1\right)$.
\end{prop}

The proof requires several lemmata, the first of which may already be found in \cite[Theorem 6.1]{Lig} and is stated here in the notation of \cite{Bov}.

\begin{lem}[see Theorem 2.1 in \cite{Bov}]\label{lem:Dirichlet form connection}
Defining 
\[\mathcal{H}^x_z:=\{h:\mathcal{S}\rightarrow [0,1]\,|\, h(x)=0,h(z)=1\}\]
and the Dirichlet form
\[\mathcal{E}(h):=\sum_{r\sim s\in\mathcal{S}}\pi(r)p(r,s) (h(r)-h(s))^2,\]
we have
\[\mathbb{P}_x(\tau_z<\tau_x)=\frac{1}{2\pi(x)}\inf_{h\in \mathcal{H}^x_z}\mathcal{E}(h).\]
\end{lem}

\begin{lem}[see Lemma 2.2 in \cite{Bov}]\label{lem:EingschränkteÜWkeiten}
For any subgraph $\Delta\subset\mathcal{S}$ containing $x,y$ and corresponding transition matrix $\widetilde{\mathbf{P}}$, we have
\[\mathbb{P}_x(\tau_y<\tau_x)\geq \widetilde{\mathbb{P}}_x(\tau_y<\tau_x).\]
\end{lem}


\begin{lem}[see Lemma 2.5 in \cite{Bov}]\label{lem:Lemma 2.5 Bovier}
If $\Delta=(\omega_0{,}...,\omega_k)$ is any one-dimensional subgraph of $\mathcal{S}$, then
\[\widetilde{\mathbb{P}}_{\omega_0}(\tau_{\omega_k}<\tau_{\omega_0}) =\left(\sum_{i=1}^k\frac{\pi(\omega_0)}{\pi(\omega_i)}\frac{1}{p(\omega_i,\omega_{i-1})}\right)^{-1}.\]
\end{lem}

\begin{bew}\emph{(of Proposition \ref{prop:AbschaetzungUeWkeiten})} In view of Lemma \ref{lem:Dirichlet form connection}, we must find an appropriate function $h$ for the upper bound. Let us define
\[\mathcal{R}:=\{s\in \mathcal{S}\,|\,\exists ~ \gamma\in \Gamma(x,z):\argmax_{i=1{,}... , |\gamma|}E(\gamma_i)=s\},\]
the set of all peak states (with respect to the energy function) along self-avoiding paths from $x$ to $z$. Obviously $z^*(x,z)\in \mathcal{R}$ and, by non-degeneracy, $E(z^*(x,z))<E(s)$ for all $s\in \mathcal{R}\backslash\{z^*(x,z)\}$. The set $\mathcal{R}$ divides $\mathcal{S}$ into a set $\mathcal{R}_x$ containing $x$ (and consisting of those $y$ that can be reached from $x$ without hitting $\mathcal{R}$) and a set $\mathcal{R}_z$ containing $z$. Now choose 
\[h(r):=\mathds{1}_{\mathcal{R}_z\cup \mathcal{R}}(r), ~r\in\mathcal{S}.\]
Note that only neighboring states contribute to the Dirichlet form $\mathcal{E}(h)$ and that, for any two such $r\sim s$, the term $(h(r)-h(s))^2$ is positive iff one of these states is in $\mathcal{R}$ and the other one in $\mathcal{R}_x$, in which case the squared difference equals 1. By invoking Lemma \ref{lem:statVert}, we obtain
\begin{align*}
\mathbb{P}_x(\tau_z<\tau_x)\ 
&=\ \frac{1}{2\pi(x)}\inf_{h\in \mathcal{H}^x_z}\mathcal{E}(h)\\
&\leq\ \frac{\pi(z^*(x,z))}{\pi(x)}\left(\sum_{\mathcal{R}\ni r\neq z^*(x,z)}\,
\sum_{s\sim r}\frac{\pi(r)}{\pi(z^*(x,z))}\,p(r,s)+|\mathcal{N}(z^*(x,z))|\right)\\
&\leq\ e^{-\beta (E(z^*(x,z))-E(x)-2\gamma_{\beta})}\\
&\hspace*{1cm}\times\left(\sum_{\mathcal{R}\ni r\neq z^*(x,z)}\,
\sum_{s\sim r}e^{-\beta(\min_{a\neq b:E(a)>E(b)} (E(a)- E(b))-2\gamma_{\beta})}+|\mathcal{N}(z^*(x,z))|\right)\\
&\leq\ e^{-\beta (E(z^*(x,z))-E(x)-2\gamma_{\beta})}\\
&\hspace*{1cm}\times\,\max_{r\in\mathcal{S}}|\mathcal{N}(r)|\,
\left(|\mathcal{S}|e^{-\beta(\min_{a\neq b:E(a)>E(b)}(E(a)- E(b))-2\gamma_{\beta})}+1\right).
\end{align*}

Lemmata \ref{lem:EingschränkteÜWkeiten} and \ref{lem:Lemma 2.5 Bovier} will enter in the proof of the lower bound. Consider the chain restricted to the one-dimensional subgraph given by a minimal path $\rho=(s_{1}{,}...,s_{|\rho|})$ from $x$ to $y$. Then
\begin{align*}
\mathbb{P}_x(\tau_y<\tau_x)\ 
&\geq\ \left(\sum_{i=1}^{|\rho|}\frac{\pi(x)}{\pi(s_i)}\frac{1}{p(s_i,s_{i-1})}\right)^{-1}\\
&\geq\ \frac{\pi(z^*(x,y))}{\pi(x)}\left(\sum_{i=1}^{|\rho|}\frac{\pi(z^*(x,y))}{\pi(s_i)}e^{\beta((E(s_{i-1})-E(s_i))^{+} +\gamma_{\beta})}\right)^{-1}\\
&\geq\ e^{-\beta(E(z^*(x,y))-E(x) +2\gamma_{\beta})}\left(\sum_{i=1}^{|\rho|}e^{-\beta(E(z^*(x,y))-E(s_i)-(E(s_{i-1})-E(s_i))^{+}-3\gamma_{\beta})}\right)^{-1}\\
&\geq\ e^{-\beta(E(z^*(x,y))-E(x)+5\gamma_{\beta})}\,\frac{1}{|\mathcal{S}|}.
\end{align*}
This completes the proof of Proposition \ref{prop:AbschaetzungUeWkeiten}.
\end{bew}

We proceed to the second proposition needed to prove Theorem \ref{thm:DriftZumMinimum}.

\begin{prop}[see Corollary 1.6 in \cite{Bov}]\label{prop:Cor 1.6 Bovier}
Given $I\subset\mathcal{S}$ and distinct $x,z\in\mathcal{S}\backslash I$,
\[\mathbb{P}_x(\tau_z<\tau_I)=\frac{\mathbb{P}_x(\tau_z<\tau_{I\cup\{x\}})}{\mathbb{P}_x(\tau_{I\cup\{z\}}<\tau_x)}\]
holds true.
\end{prop}

With the help of Propositions \ref{prop:AbschaetzungUeWkeiten} and \ref{prop:Cor 1.6 Bovier}, the proof of Theorem \ref{thm:DriftZumMinimum} can now be given quite easily.

\begin{bew}\emph{(of Theorem \ref{thm:DriftZumMinimum})}
By first using the previous result and then Proposition \ref{prop:AbschaetzungUeWkeiten} (with $K(\beta)$ as defined there), we find
\begin{align*}
\mathbb{P}_{x}(\tau_{z}<\tau_{y})\ 
&=\ \frac{\mathbb{P}_{x}(\tau_{z}<\tau_{\{x,y\}})}{\mathbb{P}_{x}(\tau_{\{z,y\}}<\tau_{x})}\\
&\leq\ \frac{\mathbb{P}_{x}(\tau_{z}<\tau_{x})}{\mathbb{P}_{x}(\tau_{y}<\tau_{x})}\\
&\leq\ K(\beta)\,e^{-\beta(E(z^*(x,z))-E(z^*(x,y))-7\gamma_{\beta})}.
\end{align*}
The argument is completed by noting that $E(z^*(x,z))>E(z^*(x,y))$ and $K(\beta)\ge 0$ converges to
\[|\mathcal{S}|\max_{r\in\mathcal{S}}|\mathcal{N}(r)|\]
as $\beta\to\infty$.
\end{bew}

\begin{bew}\emph{(of Theorem \ref{thm:DriftZumMinimumTäler})}
If $B$ occurs, the asserted bound follows directly from Theorem \ref{thm:DriftZumMinimum}. Proceeding to the case when $B^{c}$ occurs, i.e.\ $E(z^*(x,y))\leq E(z^*(x,m))$, we first point out that
\begin{align*}
\mathbb{P}_{x}(\tau_{y}<\tau_{m})\ 
&=\ \mathbb{P}_{x}(\tau_{y}<\tau_{m}, E(X_n)>E(z^*(x,y)) \text{ for some }n\leq \tau_{y})\\
&\qquad+\ \mathbb{P}_{x}(\tau_{y}<\tau_{m},E(X_n)\leq E(z^*(x,y))\text{ for all }n\leq \tau_{y})\\ 
&=:\ P_{1}+P_{2}.
\end{align*}
For all $z\in \mathcal{S}$ with $E(z)>E(z^*(x,y))$, we have $z^*(x,z)\neq x$ and $E(z^*(x,z))>E(z^*(x,y))$, for
\[E(z^*(x,z))\geq E(z)>E(z^*(x,y))\geq E(x).\]
Therefore, by an appeal to Theorem \ref{thm:DriftZumMinimum},
\begin{align*}
P_{1}\ 
&\leq\ \mathbb{P}_{x}(\tau_{z}<\tau_{y}\text{ for some $z$ with $E(z)>E(z^*(x,y))$})\\
&\leq\ \sum_{z:E(z)>E(z^*(x,y))}\mathbb{P}_{x}(\tau_z<\tau_{y})\\
&\leq\ \sum_{z:E(z)>E(z^*(x,y))}\varepsilon(x,y,z,\beta).
\end{align*}
To get an estimate for $P_{2}$, note that every minimal path from $x$ to $y$ must pass through $V^{(i)}_<(m)$ (Lemma \ref{lem:Anziehung}). With this observation and by another appeal to Theorem \ref{thm:DriftZumMinimum}, we infer 
\begin{align*}
\mathbb{P}_{x}(\tau_{y}<\tau_{m},E(X_n)\leq E(z^*(x,y))\text{ for all }n\leq\tau_{y})\ 
&\leq \sum_{z\in V^{(i)}_<(m)}\mathbb{P}_z(\tau_{y}<\tau_ {m})\\
&\leq\ \sum_{z\in V^{(i)}_<(m)}\varepsilon(z,m,y,\beta),
\end{align*}
having further utilized that (by Proposition \ref{prop:essentiellerSattel}(f) and (e))
$z^*(z,y)\neq z$ and
\begin{align*}
E(z^*(z,m))&<E(z^*(z,m'))\\
&\leq E(z^*(z,y))\vee E(z^*(y,m'))\\
&\leq E(z^*(z,y))\vee E(z^*(y,m))\\
&\leq E(z^*(z,y))\vee E(z^*(z,m))\\
&=E(z^*(z,y))
\end{align*}
for some $m'\in M^{(i)}\backslash \{m\}$ with $E(z^*(y,m))\geq E(z^*(y,m'))$, which must exist since $y\notin V^{(i)}(m)$.
\end{bew}

\begin{bew}\emph{(of Theorem \ref{thm:DriftZumMinimumTäler2})}
We first note that $y\notin V^{(l_j)}(m_j)$ for all $1\leq j\leq k$. With $m_{0}{,}...,m_{k}$ as stated in the theorem (recall $m_{0}=x$ and $m_{k}=m$), we obtain
\begin{align*}
\Prob_{x}(\tau_{y}<\tau_{m})
&=\Prob_{x}(\tau_{m_{1}}<\tau_{y}<\tau_{m})+\Prob_{x}(\tau_{y}<\tau_{m_1}\wedge\tau_{m})\\
&\leq\Prob_{m_{1}}(\tau_{y}<\tau_{m})+\Prob_{x}(\tau_{y}<\tau_{m_{1}})\\
&\leq\Prob_{m_{2}}(\tau_{y}<\tau_{m})+\Prob_{m_{1}}(\tau_{y}<\tau_{m_{2}})+\Prob_{x}(\tau_{y}<\tau_{m_{1}})\\
&\hspace{1cm}\vdots\\
&\leq\sum_{j=1}^{k}\Prob_{m_{j-1}}(\tau_{y}<\tau_{m_{j}}).
\end{align*}
Finally use Theorem \ref{thm:DriftZumMinimumTäler} to infer
\[\Prob_{m_{j-1}}(\tau_{y}<\tau_{m_{j}})\le\tilde{\varepsilon}(m_{j-1},m_{j},y,\beta)\]
for each $j=1{,}...,k$.
\end{bew}

\subsection{Quasi-stationarity and exit time}

Naturally, several other questions concerning the behavior of the process when moving in a fixed valley are of interest, and \emph{quasi-stationarity} may appear as one to come up with first. For a given valley $V$ (of any level), a quasi-stationary distribution $\nu=(\nu(j))_{j\in V}$ is characterized by the \emph{quasi-invariance}, viz.
\begin{equation}\label{eq:quasi-invariance}
\mathbb{P}_{\nu}(X_{n}=j|\tau_{\mathcal{S}\backslash V}>n)=\nu(j)\quad\text{for all }j\in V,
\end{equation}
but also satisfies
\begin{equation}\label{eq:asympt quasi-invariance}
\lim_{n\to\infty}\,\mathbb{P}_{\mu}(X_{n}=j|\tau_{\mathcal{S}\backslash V}>n)=\nu(j)\quad\text{for all }j\in V
\end{equation}
if $\mu$ is an arbitrary distribution with $\mu(V)=1$. The latter property renders uniqueness of $\nu$. Since $\mathcal{S}$ is finite, existence of $\nu$ follows by an old result due to \textsc{Darroch \& Seneta} \cite{Dar}. It is obtained as the normalized eigenvector of the Perron-Frobenius eigenvalue $\lambda=\lambda(V)$ of a modification of $\mathbf{P}$. This eigenvalue $\lambda$ is also the probability for the chain to stay in $V$ at least one step when started with $\nu$, thus $\mathbb{P}_{\nu}(\tau_{V^{c}}>1)=\lambda$. As an immediate consequence, one finds that the exit time $\tau_{V^{c}}$ has a geometric distribution with parameter $1-\lambda$ under $\mathbb{P}_{\nu}$. In the present context, this naturally raises the question how the parameter $\lambda$ relates to the transition probabilities or the energies of the valley $V$. 
A simple probabilistic argument shows the following basic and intuitively obvious result concerning the eigenvalues associated with the nesting $V^{(1)}(m)\subset...\subset V^{(i)}(m)$ (Proposition \ref{prop:TaelerRekursion}) for any $1\leq i\leq \mathfrak{n}$ and $m\in M^{(i)}$.

\begin{prop}\label{prop:Eigenwerte}
Fixing any $1\leq i\leq \mathfrak{n}$ and $m\in M^{(i)}$, let $\lambda^{(j)}:=\lambda(V^{(j)}(m))$ for $j=1{,}...,i$. Then $\lambda^{(1)}\leq...\leq \lambda^{(i)}$.
\end{prop}

\begin{bew}
Write $\nu_{j}$ as shorthand for the quasi-stationary distribution on $V^{(j)}(m)$ and $T_{j}$ for $\tau_{\mathcal{S}\backslash V^{(j)}(m)}$. Plainly, $T_{j}\le T_{j+1}$
\begin{align}\label{eqn:NestingEigenvalues}
(\lambda^{(j)})^{n}\ 
=\ \mathbb{P}_{\nu_{j}}(T_{j}>n)\ 
&\le\ \mathbb{P}_{\nu_{j}}(T_{j+1}>n)\nonumber\\
&=\ \int_{\{T_{j+1}>k\}}\mathbb{P}_{X_{k}}(T_{j+1}>n-k)\ d\mathbb{P}_{\nu_{j}}\nonumber\\
&=\ \mathbb{P}_{\nu_{j}}(T_{j+1}>k)\,\mathbb{P}_{\mu_{k}}(T_{j+1}>n-k),
\end{align}
where $\mu_{k}(x):=\mathbb{P}_{\nu_{j}}(X_{k}=x|T_{j+1}>k)$ for $x\in V^{(j+1)}$. Since $\mathcal{S}$ is finite and by virtue of \eqref{eq:asympt quasi-invariance}, we have that $\mu_{k}\le 2\nu_{j+1}$ when choosing $k$ sufficiently large. For any such $k$, we find that \eqref{eqn:NestingEigenvalues} has upper bound
\[2\,\mathbb{P}_{\nu_{j}}(T_{j+1}>k)\,\mathbb{P}_{\nu_{j+1}}(T_{j+1}>n-k)\ 
=\ 2\,\mathbb{P}_{\nu_{j}}(T_{j+1}>k)\,(\lambda^{(j+1)})^{n-k}.\]
Hence, we finally conclude
\[\lambda^{(j)}\ \le\ \Big(2\,\mathbb{P}_{\nu_{j}}(T_{j+1}>k)\,(\lambda^{(j+1)})^{-k}\Big)^{1/n}\,\lambda^{(j+1)}\]
and thereby the assertion upon letting $n\to\infty$.
\end{bew}


An alternative matrix-analytic proof draws on an old result by \textsc{Frobenius} \cite{Frob}, here cited from \cite[Chapter III, \S 2, Lemma 2]{Gan}.

\begin{lem}\label{lem:Gantmacher lemma}
If $A=(a_{ij})$ and $C=(c_{ij})$ denote two real $k\times k$-matrices such that $A$ is nonnegative and irreducible with maximal eigenvalue $\lambda_A^*$ and $|c_{ij}|\leq a_{ij}$ for all $1\leq i,j\leq k$, then $|\lambda|\le\lambda_A^*$ for all eigenvalues $\lambda$ of $C$.
\end{lem}

\noindent
\textit{Second proof of Proposition \ref{prop:Eigenwerte}:}
For any fixed valley $V$, collaps all states $s\notin V$ into an absorbing state (grave) $\Delta$ which leaves transition probabilities between states in $V$ unchanged. A proper rearrangement of states allows us to assume that the new transition matrix has the form
\begin{equation*}
\mathbf{P}=\begin{pmatrix} 1 & {\bf 0} \\ \mathbf{p} & \mathbf{Q}\end{pmatrix}
\end{equation*}
for a $|V|\times 1$-column vector $\mathbf{p}\neq 0$ and a nonnegative, substochastic and irreducible $|V|\times|V|$-matrix $\mathbf{Q}$. Now, for any $2\le j\le i$, let $A$ be this matrix $\mathbf{Q}$ when $V=V^{(j)}(m)$, and $D$ be this matrix when $V=V^{(j-1)}(m)$. Then, obviously,
\[A:=\begin{pmatrix} A_1 & A_2 \\ A_3 & D  \end{pmatrix}\]
and $A$ is irreducible and nonnegative with maximal eigenvalue $\lambda^{(j)}$. Defining further
\[C:=\begin{pmatrix} 0 & 0 \\ 0 & D \end{pmatrix}.\] 
the largest eigenvalue of $C$ equals the largest eigenvalue of $D$, thus $\lambda^{(j-1)}$. Finally, the desired conclusion follows from the previous lemma, since $|c_{ij}|=c_{ij}\leq a_{ij}$ for all $i\leq i,j\leq k$.\hfill $\square$\break
\vspace{.2cm}

Another question is how long a given valley is visited and thus about its exit time. There is an extensive literature on exit problems for different kinds of stochastic processes. We mention \cite[Ch.\ XI.2]{vKam} and \cite[Ch.\ 4, \S 4, Theorem 4.1]{FrWe} as two related to our work. The latter one studies perturbed systems on a continuous space. We can discretize their argument to get, with use of the main theorem in \cite{Scop}, a nice result on the time needed to leave a valley $V^{(i)}(m)$ for any fixed $1\leq i\leq \mathfrak{n}$ and $m\in M^{(i)}$. This result is more explicit than the one in \cite[Proposition 4.6]{OlSc96}.

\begin{defn}\label{defn:Stoppingtimes}
For $1\leq i\leq \mathfrak{n},\, N:=N^{(i)}$, we define the following stopping (entrance/exit) times:
\begin{align*}
\xi_0^{(i)}&:=\tau_{N^c}\\
\zeta_{n}^{(i)}&:=\inf\left\{k\geq \xi_n^{(i)}|X_k\in N\right\}\\
\xi_{n+1}^{(i)}&:=\inf\left\{k\geq \zeta_{n}^{(i)}|X_k\in N^c\right\}, ~n\geq 0.
\end{align*}
\end{defn}

The entrance times $\xi^{(i)}_n$ mark the epochs when a new valley is visited, while the exit times $\zeta^{(i)}_n$ are the epochs at which a valley is left. The reader should notice that we do not restrict ourselves to valleys of order $i$ but include those valleys which up to order $i$ have not yet been absorbed by some larger valley. Exit and entrance times never coincide since there is no way to go from one valley to another without hitting a non-assigned state - crests are always non-assigned (see Lemma \ref{lem:outer boundary nonassigned}).

In this section, we will focus on $\zeta_0^{(i)}$ for any fixed $i$, thus writing $\zeta_0:=\zeta_0^{(i)}$ hereafter, but later for the macroscopic process the other times will be needed as well.

\vspace{.2cm}
For each valley $V^{(i)}(m),\,m\in M^{(i)}$, let us define
\begin{equation}\label{eqn:def min state outer boundary}
s_m\ =\ s_m^{(i)}\ :=\ \argmin_{s\in \partial^+V^{(i)}(m)}E(s)\ =\ \argmin_{s\in \partial^+V^{(i)}(m)}E(z^*(m,s)),
\end{equation}
where the second equality follows from Lemma \ref{lem:outer boundary nonassigned}. 

\begin{thm}\label{thm:Verlassenszeit}
Let $m\in M^{(i)}$. Then
\[\lim_{\beta\rightarrow \infty}\frac{1}{\beta}\ln\mathbb{E}_{r}\zeta_0\ =\ E(s_m)-E(m)\]
for any $r\in V^{(i)}(m)$.
\end{thm}

For the upper bound, we need a result from \cite{Scop}, which in our notation is:

\begin{prop}[see Main Theorem (iv) in \cite{Scop}]\label{prop:Scop}
For any $1\leq i\leq \mathfrak{n},\, \beta$ sufficiently large and $t> 2^{i-1}\exp(\beta(E(s_{m^{(i-1)}})-E(m^{(i-1)})+2i|\mathcal{S}|\gamma_{\beta}))$,
\[\sup_{x\in M^{(i-1)}}\Prob_x(\tau_{M^{(i)}}>t)\leq \exp(-\Delta\beta)\]
holds true with a positive constant $\Delta$, where $M^{(0)}=\mathcal{S}$ should be recalled.
\end{prop}

This will now be used to show the following result.

\begin{lem}\label{lem:Tails2}
Fix $1\leq i\leq \mathfrak{n},\, m\in M^{(i)}$ and $r\in V^{(i)}(m)$. Then, for any $\beta$ sufficiently large and $t>2^{i}\exp(\beta(E(s_{m})-E(m)+2(i+1)|\mathcal{S}|\gamma_{\beta}))$, it holds true that
\[\Prob_r(\zeta_0<(i+1)t)\ \geq\ \frac{1}{4}.\]
\end{lem}

\begin{bew}
Let us first note that we can always arrange for $m$ being equal to $m^{(i)}$ by sufficiently decreasing the energy function at any $m'\in M^{(i)}\backslash\{m\}$ so as to make $E(s_{m})-E(m)$ minimal among all states in $M^{(i)}$. This affects neither the valley $V^{(i)}(m)$ and its outer boundary nor the distribution of $\zeta_0$ when starting in $m$, for this distribution does not depend on the energy landscape outside of $V^{(i)}(m)\cup \partial^+V^{(i)}(m)$. When applying the previous proposition, the constant $\Delta$ may have changed but is still positive which suffices for our purposes. So let $m=m^{(i)}$ hereafter.

\vspace{.2cm}
Fix $t>2^i\exp(\beta(E(s_{m})-E(m)+2(i+1)|\mathcal{S}|\gamma_{\beta}))$ and $T:=it$. Since 
\[E(s_{m})-E(m)\geq E(s_{m^{(j)}})-E(m^{(j)})\]
for every $1\leq j\le i$, we infer
\begin{align*}
\Prob_r(\tau_{M^{(i)}}\leq T)\ 
&\geq\ \Prob_r(\tau_{M^{(i)}}\leq T, \tau_{M^{(1)}}\leq t)\\
&\geq\ \Prob_r(\tau_{M^{(1)}}\leq t)\,\inf_{x\in M^{(1)}}\Prob_{x}(\tau_{M^{(i)}}\leq (i-1)t)\\
&\geq\ \Prob_r(\tau_{M^{(1)}}\leq t)\,\inf_{x\in M^{(1)}}\Prob_{x}(\tau_{M^{(2)}}\leq t)\, 
\inf_{x\in M^{(2)}}\Prob_{x}(\tau_{M^{(i)}}\leq (i-2)t)\\
&\hspace*{1cm}\vdots \\
&\geq\ \prod_{j=1}^{i}\inf_{x\in M^{(j-1)}}\Prob_{x}(\tau_{M^{(j)}}\leq t)\\
&\geq\ \big(1-\exp(-\Delta\beta)\big)^i\\
&\geq\ \frac{3}{4}
\end{align*}
for $\beta$ sufficiently large. Furthermore, for $\beta$ so large that $\Prob_r(\tau_{M^{(i)}}<\tau_m)\leq 1/4$, we find that
\begin{align*}
\Prob_r(\tau_{M^{(i)}}\leq T)\ 
&=\ \Prob_r(\tau_{M^{(i)}}=\tau_m\leq T)+\Prob_r(\tau_{M^{(i)}}\leq T, \tau_{M^{(i)}}<\tau_m)\\
&\leq\ \Prob_r(\tau_m\leq T)+\Prob_r(\tau_{M^{(i)}}<\tau_m)\\
&\leq\ \Prob_r(\tau_m\leq T)+\frac{1}{4}.
\end{align*}
By combining both estimates, we obtain
\begin{equation*}\label{eqn:lemTails1}
\Prob_r(\tau_m\leq T)\ \geq\ \Prob_r(\tau_{M^{(i)}}\leq T)-\frac{1}{4}\ \geq\ \frac{1}{2}.
\end{equation*}
Hence, state $m$ is hit in time $T$ with at least probability $1/2$ when starting in $r$. Since $m=m^{(i)}$, we further have
\begin{align*}\label{eqn:lemTails2}
\Prob_m(\zeta_0\leq t)\ \geq\ \Prob_m(\tau_{M^{(i+1)}}\leq t)\ \geq\ 1-\exp(-\Delta\beta)\ \geq\ \frac{1}{2}
\end{align*}
for $\beta$ sufficiently large. Hence, state $s_m$ is hit in time $t$ with at least probability $1/2$ when starting in $m$. By combining the estimates, we finally obtain
\begin{align*}
\mathbb{P}_{r}(\zeta_0\leq(i+1)t)\ 
&\geq\ \mathbb{P}_{r}(\zeta_0\leq T+t|\tau_{m}\leq T)\,\mathbb{P}_{r}(\tau_{m}\leq T)\\
&\geq\ \mathbb{P}_{r}(\tau_{m}\leq T)\,\mathbb{P}_{m}(\zeta_0\leq t)\\
&\geq\ \frac{1}{4},
\end{align*}
which proves our claim.
\end{bew}

\begin{bew}\textit{(of Theorem \ref{thm:Verlassenszeit})}
Using the lemma just shown, we infer
\begin{align*}
\mathbb{E}_{r}(\zeta_0)\ 
&\leq\ (i+1)t\sum_{n\geq0}(n+1)\,\mathbb{P}_{r}\left(n(i+1)t\leq\zeta_0<(n+1)(i+1)t\right)\\
&=\ (i+1)t\sum_{n\geq 0}(n+1)
\Big(\mathbb{P}_{r}\left(\zeta_0\geq n(i+1)t\right)-\mathbb{P}_{r}\left(\zeta_0\geq(n+1)(i+1)t\right)\Big)\\
&=\ (i+1)t\sum_{n\geq0}\mathbb{P}_{r}\left(\zeta_0\geq n(i+1)t\right)\\
&\leq\ (i+1)t\sum_{n\geq0}\left(\max_{x\in V}\mathbb{P}_x\left(\zeta_0\geq (i+1)t\right)\right)^n\\
&\leq\ (i+1)t\sum_{n\geq0}\left(\frac{3}{4}\right)^n\\
&=\ 4(i+1)t,
\end{align*}
where $t:=2^i\exp(\beta(E(s_{m})-E(m)+2(i+1)|\mathcal{S}|\gamma_{\beta}))+1$. Since $\gamma_{\beta}\rightarrow 0$, we get in the limit
\[\lim_{\beta\rightarrow \infty}\frac{1}{\beta}\ln\mathbb{E}_{r}\zeta_0\ \leq\ E(s_m)-E(m)\]
for all $r\in V^{(i)}(m)$.

\vspace{.2cm}
Turning to the lower bound, define a sequence of stopping times, viz.\ $\rho_0:=0$ and
\[\rho_n:=\inf\{k>\rho_{n-1}|X_k=m\text{ or }X_{k}\in\partial^+ V\}\]
for $n\ge 1$. Then $Z_n:=X_{\rho_n}$, $n\ge 0$, forms a Markov chain the transition probabilities of which when starting in $m$ can be estimated with the help of Proposition \ref{prop:AbschaetzungUeWkeiten}, namely
\begin{align*}
\mathbb{P}(Z_1\in\partial^+ V^{(i)}(m)|Z_0=m)\ 
&=\ \mathbb{P}_{m}(\rho_1=\zeta_0)\\
&=\ \mathbb{P}_{m}(\zeta_0<\tau_m)\\
&\leq\ \sum_{s\in \partial^+V}\Prob_m(\tau_s<\tau_m)\\
&\leq\ K(\beta)\,e^{-\beta(\min_{s\in \partial^+V}E(z^*(m,s))-E(m)-2\gamma_{\beta})}\\
&=\ K(\beta)\,e^{-\beta(E(s_m)-E(m)-2\gamma_{\beta})}
\end{align*}
where $K(\beta)\to K\in (0,\infty)$ as $\beta\to\infty$. Further defining $\nu:=\inf\{k\geq0| Z_k\in \partial^+ V\}$, this implies in combination with a geometric trials argument that
\[\mathbb{P}_{m}(\nu>n)\ \geq\ \left(1-K(\beta)\,e^{-\beta(E(s_m)-E(m)-2\gamma_{\beta})}\right)^{n-1}.\] 
As a consequence,
\begin{align*}
\mathbb{E}_{m}\zeta_0\ 
&=\ \sum_{n\geq 1} \mathbb{E}_{m}(\underbrace{\rho_n-\rho_{n-1}}_{\geq 1})\mathds{1}_{\{\nu\geq n\}}\ 
\geq\ \sum_{n\geq 1} \mathbb{P}_{m}(\nu\geq n)\ 
\geq\ K(\beta)^{-1}\,e^{\beta(E(s_m)-E(m)-2\gamma_{\beta})}.
\end{align*}
For arbitrary $r\in V$, we now infer
\begin{align*}
\mathbb{E}_{r}\zeta_0\ 
&=\ \mathbb{E}_{r}\zeta_0\mathds{1}_{\{\zeta_0\leq \rho_1\}}+\mathbb{E}_{r}\zeta_0\mathds{1}_{\{\zeta_0>\rho_1\}}\\
&\geq\ \mathbb{E}_{r}\big(\mathbb{E}_r(\zeta_0\mathds{1}_{\{\zeta_0>\rho_1\}}|X_{\rho_{1}}=m)\big)\\
&\geq\ \mathbb{E}_{r}\mathds{1}_{\{\zeta_0>\rho_1\}}\mathbb{E}_{m}\zeta_0\\
&\geq\ \mathbb{P}_{r}(\zeta_0>\rho_1)K(\beta)^{-1}\,e^{\beta(E(s_m)-E(m)-2\gamma_{\beta})}\\
&\geq\ \frac{1}{2}\,K(\beta)^{-1}\,e^{\beta(E(s_m)-E(m)-2\gamma_{\beta})}
\end{align*}
for all sufficiently large $\beta$, because $\lim_{\beta\to\infty}\mathbb{P}_{r}(\zeta_0>\rho_1)=1$ (Theorem \ref{thm:DriftZumMinimumTäler2}). Finally, by taking logarithms and letting $\beta$ tend to $\infty$, we arrive at the inequality
\[\lim_{\beta\rightarrow \infty}\frac{1}{\beta}\ln\mathbb{E}_{r}\zeta_0\ \geq\ E(s_m)-E(m)\]
which completes the proof.
\end{bew} 

In \cite{He08}, $E(s_m)-E(m),\,m\in M^{(i)},$ is referred to as the depth of the valley $V(m)$. Therefore, Property 3 from the Introduction holds true and we can relate thermodynamics of the system (energies) to dynamics of the chain (holding times) in a very precise way. Especially, there is no universal scale for the times spent in different valleys because in general they differ exponentially.

%
%

\section{Macroscopic View: Transitions between valleys}\label{sec:MakroskopischerProzess}

With the help of the nested state space decompositions into valleys of different orders and around bottom states of different stability, we will now be able to provide an appropriate definition of the metabasins (MB) that has been announced and to some extent discussed in the Introduction. We will further define and study macroscopic versions of the original process $X=(X_{n})_{n\ge 0}$. These are obtained by choosing different levels of aggregation in the sense that they keep track only of the valleys of a chosen level that are visited by $X$. The motivation behind this approach is, on the one hand, to exhibit strong relations between properties of the energy landscape and the behavior of $X$ (as in Theorem \ref{thm:Verlassenszeit}) and, on the other hand, to describe essential features of this process by looking at suitable macroscopic scales. 

\vspace{.1cm}
In the subsequent definition of aggregated versions of $X$, we will distinguish between two variants:
\begin{enumerate}
\item[$\bullet$] A \emph{time-scale preserving} aggregation that, for a fixed level and each $n$, keeps track of the valley the original chain visits at time $n$ and thus only blinds its exact location within a valley.
\item[$\bullet$] An \emph{accelerated} version that, while also keeping track of the visited valleys, further blinds the sojourn times within a valley by counting a visit just once.
\end{enumerate}

Actually, the definition of these aggregations at a chosen level $i$ is a little more complicated because their state space, denoted $\mathcal{S}^{(i)}$ below and the elements of which we call \emph{level $i$ metastates}, also comprises the non-assigned states at level $i$ as well as the minima of those valleys that were formed at an earlier level and whose merger is pending at level $i$ because their minima are not attracted at this level.

\begin{defn}\label{def:macroscopic processes}
Fix $1\leq i\leq \mathfrak{n}$, let $\mathcal{S}^{(i)}:=\{m^{(j)}\in M^{(1)}|\,l(j)>i\}\cup N^{(i)}$ and 
\begin{align*}
V^{(i)}(s):=
\begin{cases} 
V^{(i)}(m^{(j)}),&\text{if }s=m^{(j)} \textrm{ for some $j\geq i$}\\
V^{(j)}(m^{(j)}),&\text{if }s=m^{(j)} \textrm{ for some $j<i$}\\ 
\hfill\{s\},&\textrm{if $s\in N^{(i)}$}
\end{cases}
\end{align*}
for $s\in \mathcal{S}^{(i)}$. Then define
\begin{align*}
\overline{Y}_n^{(i)}\ &:=\ \sum_{s\in \mathcal{S}^{(i)}}s\,\mathds{1}_{\{X_{n}\in V^{(i)}(s)\}},\quad n\ge 0,\\
\text{and}\quad 
Y^{(i)}_n\ &:=\ \overline{Y}^{(i)}_{\sigma_n},\quad n\ge 0,
\end{align*}
where $\sigma_{0}=\sigma_{0}^{(i)}:\equiv 0$ and
\begin{equation*}
\sigma_{n} = \sigma_{n}^{(i)} := \inf\left\{k>\sigma_{n-1}\Big|\overline{Y}_k^{(i)}\neq \overline{Y}_{k-1}^{(i)}\right\}
\end{equation*}
for $n\ge 1$. We call $\overline{Y}^{(i)}=(\overline{Y}_n^{(i)})_{n\ge 0}$ and $Y^{(i)}=(Y_{n}^{(i)})_{n\ge 0}$ the \emph{aggregated chain (AC)} and the \emph{accelerated aggregated chain (AAC) (at level $i$)} associated with $X=(X_n)_{n\geq 0}$.
\end{defn}

So, starting in an arbitrary valley, the original chain stays there for a time $\zeta_0=\zeta_{0}^{(i)}$ (as defined in Definition \ref{defn:Stoppingtimes}) before it jumps via some non-assigned states $k_1{,}..., k_l$ (staying a geometric time in each of these states) to another valley at time $\xi_1=\xi_{1}^{(i)}$. There it stays for $\zeta_1-\xi_1$ time units before it moves on in a similar manner. By going from $X$ to its aggregation $\overline{Y}^{(i)}$ at level $i$, we regard the whole valley $V^{(i)}(s)$ for $s\in \mathcal{S}^{(i)}$ as one single metastate and therefore give up information about the exact location of $X$ within a valley. $\overline{Y}^{(i)}$ is a jump process on $\mathcal{S}^{(i)}$ with successive sojourn times $\sigma_{n+1}-\sigma_n,\,n\geq 0$, which do not only depend on the valley but also on the states of entrance and exit. The AAC then is the embedded chain, viz.
\[\overline{Y}^{(i)}_n=\sum_{j\geq 0}Y^{(i)}_j\mathds{1}_{\{\sigma_j\leq n<\sigma_{j+1}\}},\]
giving the states only at jumps epochs: starting from the minimum of a first valley it moves to states $k_1{,}...,k_l\in N^{(i)}$ and then proceeds to the minimum of a second valley, and so on.

Of course, at small temperatures the time spent in a non-assigned state or in a valley around a low order metastable state is very small compared to the time spent in a valley around a metastable state of higher order. Thus, such states can be seen as instantaneous and of little importance for the evolution of the process. We account for them nonetheless for two reasons. First, in the path-dependent definition mentioned in the Introduction and used in Physics, they build small MB of great transitional activity of the process and are thus relevant in view of our goal to provide a definition of MB that conforms as much as possible to a path-dependent one. Second, a complete partitioning of the state space that is an assignment of every $s\in \mathcal{S}$ to a metastate via a global algorithm fails when merely focusing on $\left\{V^{(i)}(m), m\in M^{(i)}\right\}$ because there is neither an obvious nor natural way how to assign non-assigned states to them.

\vspace{.2cm}
The \emph{incoherent scattering function} and its associated relaxation time, for $X$ defined by
\[S(q,n):=\mathbb{E}_{\pi}\cos\big(q|X_n-X_0|\big)\]
(with $|\cdot|$ being Euclidean distance in phase space) and
\[\tau_q(\varepsilon):=\inf\{n|S(q,n)\leq \varepsilon\},\quad\varepsilon>0,\]
respectively, may serve as an example which shows the strong relation between the behavior of the original process and its macroscopic versions. For more detailed information on the meaning and relevance of $S(q,n)$ as a measure of incoherent scattering between the initial state of a glass-forming system and its state $n$ time steps onward, we refer to the survey by \textsc{Heuer} \cite{He08}.

\begin{prop}
For each $1\leq i\leq \mathfrak{n}$ there is a constant $\Delta(i)$ such that
\begin{equation}\label{eqn:ScatteringFunction}
\sup_{n\ge 0}\,\Prob_{\pi}\left(X_n\neq \overline{Y}^{(i)}_n\right)\leq e^{-\Delta(i)\beta}.
\end{equation}
As a consequence, for any given $\varepsilon>0$, the incoherent scattering functions of $X$ and $\overline{Y}^{(i)}$ differ by at most $4\varepsilon$ for $\beta$ sufficiently large.
\end{prop}

\begin{bew}
Use Lemma \ref{lem:statVert} to infer
\begin{align*}
\Prob_{\pi}(X_n\neq \overline{Y}^{(i)}_n)
&=\sum_{s\in \mathcal{S}^{(i)}\backslash N^{(i)}}\Prob_{\pi}(\overline{Y}^{(i)}_n=s)
\sum_{x\in V^{(i)}(s)\backslash \{s\}}\frac{\pi(x)}{\pi(V^{(i)}(s))}\\
&\leq \sum_{s\in \mathcal{S}^{(i)}\backslash N^{(i)}}\Prob_{\pi}(\overline{Y}^{(i)}_n=s)
\sum_{x\in V^{(i)}(s)\backslash \{s\}}\frac{\pi(x)}{\pi(s)}\\
&\leq \sum_{s\in \mathcal{S}^{(i)}\backslash N^{(i)}}\Prob_{\pi}(\overline{Y}^{(i)}_n=s)
\sum_{x\in V^{(i)}(s)\backslash \{s\}}e^{-\beta(E(x)-E(s)-2\gamma_{\beta})}\\
&\leq \max_{s\in \mathcal{S}^{(i)}\backslash N^{(i)}}|V^{(i)}(s)|
\max_{x\in V^{(i)}(s)\backslash \{s\}}e^{-\beta(E(x)-E(s)-2\gamma_{\beta})}.
\end{align*}
This proves equation \eqref{eqn:ScatteringFunction} because $E(x)>E(s)$ for each $x\in V^{(i)}(s)\backslash\{s\}, s\in \mathcal{S}^{(i)}\backslash N^{(i)}$. Now let $\beta$ be so large that $e^{-\Delta(i)\beta}\leq \varepsilon$ for a given $\varepsilon>0$ and observe that
\begin{align*}
\mathbb{E}_{\pi}&\cos\big(q|X_n-X_0|\big)\\
&=\ \int_{\{X_n=\overline{Y}_{n}^{(i)},X_0=\overline{Y}^{(i)}_0\}}\cos\big(q|X_n-X_0|\big)\ d\Prob_{\pi}\ 
+\ \int_{\{X_n\neq\overline{Y}_n^{(i)}\}\cup \{X_0\neq \overline{Y}^{(i)}_0\}}\cos\big(q|X_n-X_0|\big)\ d\Prob_{\pi}\\
&\leq\ \int_{\{X_n=\overline{Y}_n^{(i)},X_0=\overline{Y}^{(i)}_0\}}\cos\big(q|X_n-X_0|\big)\ d\Prob_{\pi}\ 
+\ \Prob_{\pi}(X_n\neq \overline{Y}_n^{(i)})\ +\ \Prob_{\pi}(X_0\neq \overline{Y}^{(i)}_0)\\
&=\ \mathbb{E}_{\pi}\cos\left(q|\overline{Y}_n^{(i)}-\overline{Y}_0^{(i)}|\right)\ 
-\ \int_{\{X_n\neq \overline{Y}_n^{(i)}\}\cup \{X_0\neq\overline{Y}^{(i)}_0\}}
\cos\left(q|\overline{Y}_n^{(i)}-\overline{Y}_0^{(i)}|\right)\ d\Prob_{\pi} +\ 2\varepsilon\\
&\leq\ \mathbb{E}_{\pi}\cos\left(q|\overline{Y}_n^{(i)}-\overline{Y}_0^{(i)}|\right)\ +\ 4\varepsilon
\end{align*}
and, by a similar argument,
\begin{align*}
\mathbb{E}_{\pi}\cos\left(q|\overline{Y}_n^{(i)}-\overline{Y}_0^{(i)}|\right)
\ \leq\ \mathbb{E}_{\pi}&\cos\big(q|X_n-X_0|\big)\ +\ 4\varepsilon.
\end{align*}
This completes the proof.
\end{bew}

\subsection{(Semi-)Markov Property}\label{subsec:MarkovProperty}

In general, both aggregated chains are no longer Markovian. Transition probabilities of the AAC not only depend on the current state, i.e.\ the current valley, but also on the entrance state into that valley, whereas transition probabilities of the AC depend on the current sojourn times which in turn depend on the previous, the present and the next state. On the other hand, since valleys are defined in such a way that asymptotically almost surely (a.a.s.), i.e., with probability tending to one as $\beta \to\infty$, the minimum will be reached from anywhere inside the valley before the valley is left, and since, furthermore, the exit state on the outer boundary a.a.s.\ equals the one with the smallest energy, the AAC will be shown below to converge to a certain Markov chain on $\mathcal{S}^{(i)}$. Also, the sojourn times depend on the past only via the last and the current state. This means that the AC converges to a semi-Markov chain (for semi-Markov chains see for example \cite{BaLi08}):

\begin{defn}\label{def:SemiMarkov}
Given any nonempty countable set $\mathcal{S}$, let $(M_{n},T_{n})_{n\geq 0}$ be a bivariate temporally homogeneous Markov chain on $\mathcal{S}\times\mathbb{N}$, with transition kernel $Q(s,\cdot)$ only depending on the first component, viz., for all $n\geq 0, s\in S$ and $t\ge 0$,
\begin{equation}\label{eq:SMC def}
\Prob(M_{n+1}=s, T_{n+1}\leq t|M_{n},T_{n})=Q(M_{n},\{s\}\times [0,t])
\end{equation}
holds. Put $S_n:=\sum_{i=0}^{n}T_i$ for $n\ge 0$ and $\nu(t):=\max\{n\geq 0|S_n\leq t\} \ (\max\emptyset :=0)$ for $t\ge 0$. Then $Z_n:=M_{\nu(n)}, n\geq 0,$ is called \textit{semi-Markov chain} with \emph{embedded Markov chain} $(M_{n})_{n\ge 0}$ and \emph{sojourn} or \emph{holding times} $T_{0},T_{1}{,}...$.
\end{defn}

Note that equation \eqref{eq:SMC def} holds iff $M=(M_{n})_{n\ge 0}$ forms a temporally homogeneous Markov chain and the $(T_{n})_{n\geq0}$ are conditionally independent given $M$ such that the distribution of $T_{n}$ only depends on $M_{n-1},M_{n}$ for $n\ge 1$ (in a temporally homogeneous manner), and on $M_{0}$ for $n=0$. Note further that we have specialized to the case where holding times take values in $\mathbb{N}$ only (instead of $(0,\infty)$).

Recall from \eqref{eqn:def min state outer boundary} the definition of $s_m$ for $m\in \mathcal{S}^{(i)}\backslash N^{(i)}$ and notice that the second equality there entails $E(z^{*}(m,s_{m}))<E(z^*(m,s)$ for any $s\in\partial^+V^{(i)}(m) \backslash\{s_{m}\}$. Further recall from our basic assumptions that $p^{*}(r,s)=\lim_{\beta\to\infty}p(r,s)$ exists for all $r,s\in\mathcal{S}$ and is positive if $E(r)\geq E(s)$. The following result, revealing the announced convergence for AAC, confirms in particular that a valley $V(m)$, $m\in \mathcal{S}^{(i)}\backslash N^{(i)}$, is a.a.s.\ to be left via $s_{m}$.

\begin{prop} \label{prop:Übergänge}
For each $1\leq i\leq \mathfrak{n}$ and as $\beta\to\infty$, the level $i$ AAC $Y^{(i)}$ converges to a Markov chain $\wh{Y}^{(i)}=(\wh{Y}_{n}^{(i)})_{n\ge 0}$ on $\mathcal{S}^{(i)}$ with transition probabilities $\wh{p}(r,s)=\wh{p}_{i}(r,s)$ stated below, that is
\begin{equation*}
\lim_{\beta\to\infty}\Prob(Y^{(i)}_{n+1}=s|Y_{n}^{(i)}=r,Y^{(i)}_{n-1}=m_{n-1}{,}...,Y^{(i)}_0=m_0)\ =\ \wh{p}(r,s)
\end{equation*}
for all $m_{0}{,}...,m_{n-1},r,s\in\mathcal{S}^{(i)}$ and $n\ge 0$. We have $\wh{p}(r,\cdot):=\delta_{s_{r}}$ if $r\in \mathcal{S}^{(i)}\backslash N^{(i)}$, and
\begin{align*}
\wh{p}(r,\cdot)\ 
:=\ \frac{1}{1-p^*(r,r)}\left(\sum_{s\in \mathcal{N}(r)\cap N^{(i)}}p^*(r,s)\,\delta_s
+\sum_{s\in \mathcal{S}^{(i)}\backslash N^{(i)}}\!\!
\left(\sum_{r'\in \mathcal{N}(r)\cap V^{(i)}(s)}p^*(r,r')\right)\delta_{s}\right),
\end{align*}
if $r\in N^{(i)}$.
\end{prop}
$\wh{Y}^{(i)}=(\wh{Y}_{n}^{(i)})_{n\ge 0}$ is called the \emph{asymptotic jump chain at level $i$} hereafter. Note that, typically, it is not irreducible. It may have transient states, not necessarily non-assigned, and its irreducibility classes are of the form $\{m_{1}{,}...,m_{k},s\}$ for a collection $m_{1}{,}...,m_{k}\in\mathcal{S}^{(i)}\backslash N^{(i)}$ and some $s\in N^{(i)}$ satisfying $s=s_{m_{1}}=...=s_{m_{k}}$.

\begin{bew}
Fix $1\leq i\leq \mathfrak{n}$ and write $Y_{n}$ for $Y_{n}^{(i)}$. The first step is to verify that, as $\beta\to\infty$,
\begin{equation*}
\Prob(Y_{n+1}=s|Y_{n}=r,Y_{n-1}=m_{n-1}{,}...,Y_0=m_0)=\Prob_{r}(Y_{1}=s)+o(1)
\end{equation*}
for all $m_{0}{,}...,m_{n-1},r,s\in\mathcal{S}^{(i)}$ and $n\ge 0$. If $r\in N^{(i)}$, then $Y_n=X_{\sigma_n}$ and the Markov property of $X$ provide us with the even stronger result
\begin{equation*}
\Prob(Y_{n+1}=s|Y_{n}=r,Y_{n-1}=m_{n-1}{,}...,Y_0=m_0)=\Prob_{r}(Y_{1}=s).
\end{equation*}
A little more care is needed if $r\in\mathcal{S}^{(i)}\backslash N^{(i)}$. For any $s\in \mathcal{S}^{(i)},\,x\in V^{(i)}(r)$ and $n\geq 0$, we have
\begin{align*}
\mathbb{P}(Y_{n+1}=s&|Y_n=r,X_{\sigma_n}=x)\\
&=\ \mathbb{P}_x(Y_1=s,\tau_{r}<\sigma_1)+\mathbb{P}_x(Y_1=s,\tau_{r}>\sigma_1)\\
&=\ \mathbb{P}_{r}(Y_1=s)\,\mathbb{P}_x(\tau_{r}<\sigma_1)+\mathbb{P}_x(Y_1=s,\tau_{r}>\sigma_1).
\end{align*}
The last two summands can further be bounded by
\begin{align*}
\mathbb{P}_{r}(Y_1=s)\,\mathbb{P}_x(\tau_{r}<\sigma_1)\ \leq\ \mathbb{P}_{r}(Y_1=s)\quad\text{and}\quad 
\mathbb{P}_x(Y_1=s,\tau_{r}>\sigma_1)\ \leq\ \mathbb{P}_x(\sigma_1<\tau_{r}).
\end{align*}
For the last probability, Theorem \ref{thm:DriftZumMinimumTäler} ensures
\begin{align*}
\mathbb{P}_x(\sigma_1<\tau_{r})\ 
&\leq\ \sum_{z\in \partial^+ V^{(i)}(r)}\mathbb{P}_x(\tau_z<\tau_{r})\ 
\leq\ \sum_{z\in \partial^+ V^{(i)}(r)} \tilde{\varepsilon}(x,r,z,\beta)\ 
\stackrel{\beta\to\infty}{\longrightarrow}\ 0.
\end{align*}
Consequently, as $\beta\to\infty$,
\begin{align*}
\mathbb{P}_{r}(Y_1=s)\ 
&=\ \left(1-\sum_{z\in \partial^+ V^{(i)}(r)} \tilde{\varepsilon}(x,r,z,\beta)\right)\mathbb{P}_{r}(Y_1=s)+o(1)\\
&\leq\ \big(1-\mathbb{P}_x(\sigma_1<\tau_{r})\big)\,\mathbb{P}_{r}(Y_1=s)+o(1)\\
&\leq\ \mathbb{P}(Y_{n+1}=s|X_{\sigma_n}=x,Y_n=r)+o(1)\\
&\leq\ \mathbb{P}_{r}(Y_1=s)+\sum_{z\in \partial^+ V^{(i)}(r)} \tilde{\varepsilon}(x,r,z,\beta)+o(1)\\
&=\ \mathbb{P}_{r}(Y_1=s)+o(1),
\end{align*}
and therefore
\begin{align*}
\Prob&(Y_{n+1}=s|Y_{n}=r,Y_{n-1}=m_{n-1}{,}...,Y_0=m_0)\\
&=\ \sum_{x\in V^{(i)}(r)}\Prob(Y_{n+1}=s|X_{\sigma_n}=x, Y_n=r)\,\Prob(X_{\sigma_n}=x|Y_{n}=r,Y_{n-1}=m_{n-1}{,}...,Y_0=m_0)\\
&=\ \mathbb{P}_{r}(Y_1=s)+o(1).
\end{align*}
It remains to verify that $\mathbb{P}_{r}(Y_1=s)=\wh{p}(r,s)+o(1)$ for any $r,s\in\mathcal{S}^{(i)}$. If $r\in\mathcal{S}^{(i)}\backslash N^{(i)}$, then $\sigma_{1}=\tau_{N^{(i)}}$ and $Y_{1}=X_{\tau_{N^{(i)}}}$. Since $E(z^*(r,s_r))<E(z^*(r,s))$ for each $s_r\neq s\in N^{(i)}\cap \partial^+V^{(i)}(r)$, we now infer with the help of Theorem \ref{thm:DriftZumMinimum}
\begin{align*}
\Prob_{r}(Y_{1}\ne s_{r})\ 
&=\ \Prob_{r}\big(\tau_{s}<\tau_{s_{r}}\text{ for some }s\in N^{(i)}\backslash\{s_{r}\}\big)\\
&\ \le\ \sum_{s_{r}\ne s\in N^{(i)}}\Prob_{r}\big(\tau_{s}<\tau_{s_{r}})\\
&\le\ \sum_{s_{r}\ne s\in N^{(i)}}\tilde\varepsilon(r,s_{r},s,\beta)\ \\
&=\ o(1),
\end{align*}
as $\beta\to\infty$ and thus $\Prob_{r}(Y_{1}\in\cdot)\to\delta_{s_{r}}=\wh{p}(r,\cdot)$ as claimed. If $r\in N^{(i)}$, then either $Y_{1}=s\in\mathcal{N}(r)\cap N^{(i)}$, or $Y_{1}=s\in\mathcal{S}^{(i)}\backslash N^{(i)}$ and $X_{\sigma_{1}}=r'$ for some $r'\in\mathcal{N}(r)\cap V^{(i)}(s)$. It thus follows that
\begin{align*}
\Prob_{r}(Y_{1}=s)\ 
=\ \Prob_{r}(X_{\sigma_{1}}=s)\ 
=\ \frac{p(r,s)}{1-p(r,r)}\ 
=\ \frac{p^{*}(r,s)}{1-p^{*}(r,r)}+o(1)
\end{align*}
if $s\in\mathcal{N}(r)\cap N^{(i)}$, while
\begin{align*}
\Prob_{r}(Y_{1}=s)\ 
=\ \sum_{r'\in\mathcal{N}(r)\cap V^{(i)}(s)}\Prob_{r}(X_{\sigma_{1}}=r')\ 
=\ \sum_{r'\in\mathcal{N}(r)\cap V^{(i)}(s)}\frac{p^{*}(r,r')}{1-p^{*}(r,r)}+o(1)
\end{align*}
in the second case.
\end{bew}

Having shown that $Y^{(i)}$ behaves asymptotically as a Markov chain, viz.\ the jump chain $\wh{Y}^{(i)}$, it is fairly easy to verify with the help of the next simple lemma that the augmented bivariate AC  $\big(\overline{Y}_{n}^{(i)},\overline{Y}_{n+1}^{(i)}\big)_{n\ge 0}$ is asymptotically semi-Markovian.

\begin{lem}\label{lem:ConditionalWaitingTimes}
For each $\beta>0$, the sojourn times $\sigma_{n+1}-\sigma_{n}$, $n\ge 0$, of the AC $\overline{Y}^{(i)}$ are conditionally independent given $Y^{(i)}$. The conditional law of $\sigma_{n+1}-\sigma_{n}$ depends only on $(Y_{n-1}^{(i)},Y_{n}^{(i)},Y_{n+1}^{(i)})$ and satisfies
\begin{align}\label{eqn:ConditionalWaitingTimes}
&\Prob\big(\sigma_{n+1}-\sigma_{n}\in\cdot~|Y^{(i)}_{n-1}=x,Y^{(i)}_{n}=y,Y^{(i)}_{n+1}=z\big)\nonumber \\
&\hspace{.5cm}=\ Q((x,y,z),\cdot\, )\ 
:=\ \begin{cases}
\hfill\textit{Geom}(1-p(y,y)),
&\text{if }y\in N^{(i)}\\ 
\sum_{s\in V^{(i)}(y),s\sim x}\Prob_s(\sigma_1\in\cdot~|Y^{(i)}_1=z)\,\Prob_{x}(X_{\sigma_1}=s),
&\text{if }y\notin N^{(i)}
\end{cases}
\end{align}
for all $x,y,z\in\mathcal{S}^{(i)}$ with $\Prob(Y^{(i)}_{n-1}=x,Y^{(i)}_{n}=y,Y^{(i)}_{n+1}=z)>0$ and $n\ge 1$.
\end{lem}

\begin{bew}
The assertions follow easily when observing that, on the one hand, at least one state $y\in N^{(i)}$ must be visited between two states $x,z\in \mathcal{S}^{(i)}\backslash N^{(i)}$ (Lemma \ref{lem:outer boundary nonassigned}) and that, on the other hand, the original chain $X$ and its aggregation $\overline{Y}^{(i)}$ coincide at any epoch where a non-assigned state is hit, which renders the Markov property of $\overline{Y}^{(i)}$ at these epochs. Further details are omitted.
\end{bew}

In order to formulate the next result, let $0=\wh\sigma_{0}<\wh\sigma_{1}<...$ be an increasing sequence of random variables such that its increments $\wh\sigma_{n+1}-\wh\sigma_{n}, n\ge 0$, are conditionally independent given the asymptotic jump chain $\wh{Y}^{(i)}$. Moreover, let the conditional law of $\wh\sigma_{n+1}-\wh\sigma_{n}$ depend only on $\big(\wh{Y}^{(i)}_{n-1},\wh{Y}^{(i)}_{n},\wh{Y}^{(i)}_{n+1}\big)$ and be equal to $Q\big(\big(\wh{Y}^{(i)}_{n-1},\wh{Y}^{(i)}_{n},\wh{Y}^{(i)}_{n+1}\big),\cdot\,\big)$, with $Q$ as defined in \eqref{eqn:ConditionalWaitingTimes}. Then $((\wh{Y}^{(i)}_{n},\wh{Y}^{(i)}_{n+1}),\wh\sigma_{n+1})_{n\geq 0}$ forms a Markov renewal process and $\big(\wh{Y}^{(i)}_{\wh\nu(n)},\wh{Y}^{(i)}_{\wh\nu(n+1)}\big)_{n\ge 0}$ a semi-Markov chain, where $\wh\nu(n):=\sup\{k\ge 0|\wh\sigma_{k}\le n\}$.

\begin{prop}\label{prop:semi-Markov AC}
For each $1\le i\le\mathfrak{n}$, $((Y^{(i)}_{n},Y^{(i)}_{n+1}),\sigma_{n+1})_{n\geq 0}$ converges to the Markov renewal process $((\wh{Y}^{(i)}_{n},\wh{Y}^{(i)}_{n+1}),\wh\sigma_{n+1})_{n\geq 0}$ in the sense that
\begin{align*}
\lim_{\beta\to\infty}
\frac{\Prob_{y_0}\big(\big(Y^{(i)}_{k},Y^{(i)}_{k+1}\big)=(y_{k},y_{k+1}),\,\sigma_{k+1}=i_{k+1},\,0\le k\le n\big)}
{\Prob_{y_0}\big(\big(\wh{Y}^{(i)}_{k},\wh{Y}^{(i)}_{k+1}\big)=(y_{k},y_{k+1}),\,\wh\sigma_{k+1}=i_{k+1},\,0\le k\le n\big)}\ =\ 1
\end{align*}
for all $y_{0}{,}...,y_{n+1}\in\mathcal{S}^{(i)}$, $0<i_{1}<...<i_{n+1}$ and $n\ge 0$ such that the denominator is positive. Furthermore, $(\overline{Y}^{(i)}_{n},\overline{Y}^{(i)}_{n+1})_{n\geq 0}$ is asymptotically semi-Markovian in the sense that
\begin{align*}
\lim_{\beta\to\infty}
\frac{\Prob_{y_0}\big(\big(\overline{Y}_{k}^{(i)},\overline{Y}_{k+1}^{(i)}\big)=(y_{k},y_{k+1}),0\le k\le n\big)}
{\Prob_{y_0}\big(\big(\wh{Y}^{(i)}_{\wh\nu(k)},\wh{Y}^{(i)}_{\wh\nu(k+1)}\big)=(y_{k},y_{k+1}),0\le k\le n\big)}\ =\ 1
\end{align*}
for all $y_{0}{,}...,y_{n+1}\in\mathcal{S}^{(i)}$ and $n\ge 0$ such that the denominator is positive.
\end{prop}

\begin{bew}
The first assertion being obvious by Proposition \ref{prop:Übergänge}, note that it implies, with $\nu(n):=\sup\{k\ge 0|\sigma_{k}\le n\}$,
\begin{align*}
\lim_{\beta\to\infty}
\frac{\Prob_{y_0}\big(\big(Y_{\nu(k)}^{(i)},Y_{\nu(k+1)}^{(i)}\big)=(y_{k},y_{k+1}),0\le k\le n\big)}
{\Prob_{y_0}\big(\big(\wh{Y}^{(i)}_{\wh\nu(k)},\wh{Y}^{(i)}_{\wh\nu(k+1)}\big)=(y_{k},y_{k+1}),0\le k\le n\big)}\ =\ 1
\end{align*}
for all $y_{0}{,}...,y_{n+1}\in\mathcal{S}^{(i)}$ and $n\ge 0$ such that the denominator is positive. Therefore the second assertion follows when finally noting that
\[Y^{(i)}_{\nu(n)}=\sum_{j\geq 0}Y^{(i)}_j\mathds{1}_{\{\sigma_j\leq n<\sigma_{j+1}\}}=\overline{Y}^{(i)}_n\]
for each $n\ge 0$.
\end{bew}

So we have shown that, although aggregation generally entails the loss of the Markov property, here it leads back to processes of this kind (Markov or semi-Markov chains) in an asymptotic sense at low temperature regimes.

\subsection{Reciprocating Jumps}\label{subsec:ÜWkeiten}

As discussed to some extent in the Introduction, we want to find an aggregation level at which reciprocating jumps appear to be very unlikely so as to obtain a better picture of essential features of the observed process. To render precision to this informal statement requires to further specify the term ``reciprocating jump'' and to provide a measure of likelihood for its occurrence. It is useful to point out first that the original chain $X$ exhibits two types of reciprocating jumps:
\begin{description}
\item[Intra-valley jumps] which occur between states inside a valley (starting in a minimum the process falls back to it many times before leaving the valley).
\item[Inter-valley jumps] which occur between two valleys (typically, when the energy barrier between these valleys is much lower then the barrier to any other valley).
\end{description}
\begin{figure}[ht]
\begin{center}
\includegraphics[trim=3cm 22.5cm 5cm 3.5cm, clip, width = 10cm]{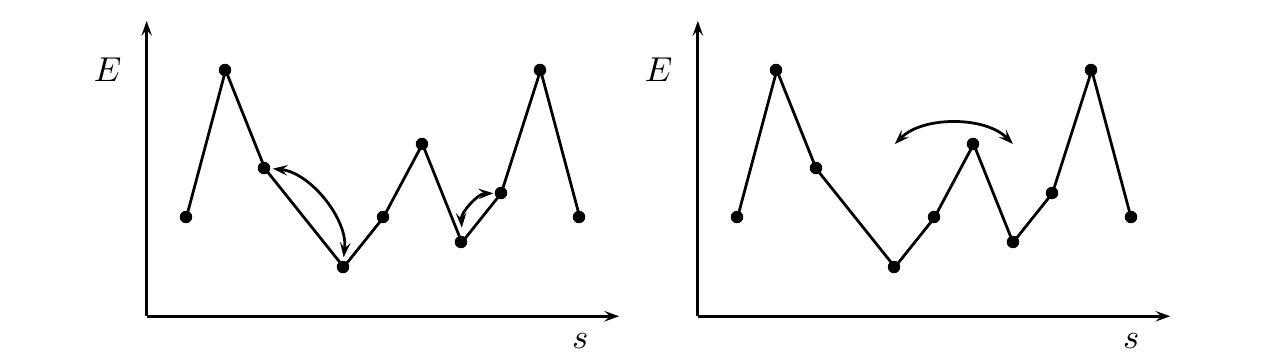}
\end{center}
\caption{Illustration of intra-valley jumps (left panel) versus inter-valley jumps (right panel).}
\label{fig:ReciprocatingJumps}
\end{figure}
Clearly, intra-valley jumps disappear by aggregating valleys into metastates, while inter-valley jumps may also be viewed as intra-valley jumps for higher order valleys and do occur when transitions between any two of them are much more likely than those to other valleys in which case they should be aggregated into one valley. This motivates the following definition.

\begin{defn}\label{def:reciprocating jumps of given order}
We say the process $(Y^{(i)}_n)_{n\in \mathbb{N}}$ exhibits \textit{reciprocating jumps of order} $\varepsilon>0$ if there exists a nonempty subset $A\varsubsetneq\mathcal{S}^{(i)}\backslash N^{(i)}$ with the following property: For each $m_{1}\in A$, there exists $m_{2}\in A$ such that
\begin{equation*}
\lim_{\beta\to \infty}\frac{1}{\beta}\left(\ln\left(
\Prob_{m_1}\left(X_{\xi_1}\in V^{(i)}(m_2)\right)\right)-\ln\left(\Prob_{m_1}\left(X_{\xi_1}\in V^{(i)}(m)\right)\right)\right)
\geq \varepsilon
\end{equation*}
for all $m\in\mathcal{S}^{(i)}\backslash (N^{(i)}\cup A)$. In other words, it is exponentially more likely to stay in $A$ than to leave it (ignoring intermediate visits to non-assigned states).
\end{defn}

In view of our principal goal to give a path-independent definition of MBs, we must point out that, by irreducibility, reciprocating jumps \emph{always} occur with positive probability at any nontrivial level of aggregation and can therefore never be ruled out completely. This is in contrast to the path-dependent version by \textsc{Heuer} \cite{He08} in which the non-occurrence of reciprocating jumps appears to be the crucial requirement. As a consequence,  Definition \ref{def:reciprocating jumps of given order} provides an alternative, probabilistic and verifiable criterion for reciprocating jumps to be sufficiently unlikely in a chosen aggregation.

\vspace{.2cm}
The following proposition contains further information on which valleys are visited successively by providing the probabilities of making a transition from $V^{(i)}(m)$ to $V^{(i)}(m')$ for any $m,m'\in\mathcal{S}^{(i)}\backslash N^{(i)}$. It is a direct consequence of the asymptotic results in the previous subsection, notably Proposition \ref{prop:Übergänge}.

\begin{prop} \label{prop:Übergänge2}
Let $m\in\mathcal{S}^{(i)}\backslash N^{(i)}, s_m\in \partial^+V^{(i)}(m)$ be as defined in \eqref{eqn:def min state outer boundary} (i.e., the state on the outer boundary of $V^{(i)}(m)$ with minimal energy). Then
\begin{equation*}
\lim_{\beta\to\infty}\mathbb{P}_m(X_{\xi_1}\in V^{(i)}(m))\ =\ \wh{p}(s_m,m)\ 
=\ \frac{\sum_{r\in\mathcal{N}(s_m)\cap V^{(i)}(m)}p^{*}(s_m,r)}{1-p^{*}(s_m,s_m)},
\end{equation*}
while
\begin{equation*}
\lim_{\beta\to\infty}\mathbb{P}_m(X_{\xi_1}\in V^{(i)}(m'))\ 
=\ \wh{p}(s_{m},m')
+\sum_{n\ge 1}\sum_{r_{1}{,}...,r_{n}\in N^{(i)}}\wh{p}(s,r_{1})\cdot...\cdot\wh{p}(r_{n-1},r_{n})\,\wh{p}(r_{n},m')
\end{equation*}
for any other $m'\in\mathcal{S}^{(i)}\backslash N^{(i)}$.
\end{prop}

The reader should notice that, as $\wh{p}(s,r)=0$ whenever $E(s)<E(r)$, the last sum actually ranges only over those non-assigned $r_{1}{,}...,r_{n}$ with $E(s_{m})>E(r_{1})>...>E(r_{n})>E(m')$.

\begin{bew}
Let us first point out that $\Prob_{r}(X_{\xi_{0}}\in V^{(i)}(m))=o(1)$ as $\beta\to\infty$ for any $r\in N^{(i)}$ such that $E(r)<E(s_{m})$. Namely, since the last property implies $r\not\in\partial^{+}V^{(i)}(m)$, any path from $r$ into $V^{(i)}(m)$ must traverse a state $s\in\partial^{+}V^{(i)}(m)$ with $E(s)\ge E(s_{m})>E(r)$, whence the probability for such a path goes to zero as $\beta\to\infty$. Noting further that $\mathbb{P}_m(Y_1^{(i)}\ne s_m)=o(1)$ as $\beta\to \infty$ by Proposition \ref{prop:Übergänge}, we now infer (with $\xi_{n}=\xi_{n}^{(i)}$)
\begin{align*}
\Prob_m(X_{\xi_1}\in V^{(i)}(m))\ 
&=\ \Prob_{s_m}(X_{\xi_0}\in V^{(i)}(m))+o(1)\\
&=\ \Prob_{s_m}(Y^{(i)}_1=m)
+\sum_{r\in\mathcal{N}(s_{m})\cap N^{(i)}}\frac{p(s_{m},r)}{1-p(s_m,s_m)}\,\Prob_{r}(X_{\xi_0}\in V^{(i)}(m))+o(1)\\
&=\ \wh{p}(s_{m},m)+o(1).
\end{align*}
The expression for $\wh{p}(s_{m},m)$ in terms of the $p^{*}(s_{m},r)$ may be read off directly from the formula given in Proposition \ref{prop:Übergänge}. For $m'\ne m$, $m'\in\mathcal{S}^{(i)}\backslash N^{(i)}$, we obtain in a similar manner
\begin{align*}
\Prob_m&(X_{\xi_1}\in V^{(i)}(m'))\\ 
&=\ \Prob_{s_{m}}(Y_{1}^{(i)}=m')\ 
+\ \sum_{n\ge 1}\sum_{r_{1}{,}...,r_{n}\in N^{(i)}} \Prob_{s_{m}}(Y_{1}^{(i)}=r_{1}{,}...,Y_{n}^{(i)}=r_{n},Y_{n+1}^{(i)}=m')+o(1)\\
&=\ \wh{p}(s_{m},m')+\sum_{n\ge 1}\sum_{r_{1}{,}...,r_{n}\in N^{(i)}}
\wh{p}(s_m,r_{1})\cdot...\cdot\wh{p}(r_{n-1},r_{n})\,\wh{p}(r_{n},m')+o(1),
\end{align*}
the last line by another appeal to the afore-mentioned proposition. 
\end{bew}

In essence, the previous result tells us that a valley $V^{(i)}(m')$ is neighbored to $V^{(i)}(m)$, that is, reachable with positive probability by the asymptotic jump chain $\wh{Y}^{(i)}$ (and thus by $Y^{(i)}$ at any temperature level $\beta$) without intermediately hitting any other valley, iff there exists at least one (in terms of energies) decreasing path in $N^{(i)}$ from $s_{m}$ to $m'$. For any other such pair of valleys, connected by a path through states in $N^{(i)}$, the transition probability decreases to zero exponentially in $\beta$. If this path can be chosen to be unimodal, here called \textit{uphill-downhill-path}, this can be stated in a very precise way as the next result shows.

\begin{lem}\label{lem:TransitionProbability}
Let $m_0,m_1\in \mathcal{S}^{(i)}\backslash N^{(i)}$ be two distinct local minima for some $0\leq i\leq \mathfrak{n}$. Suppose there exists a minimal path $\gamma=(\gamma_0{,}...,\gamma_k)$ from $s_{m_0}$ to $m_1$ not hitting any other valley but $V^{(i)}(m_{1})$ and such that $I(\gamma_0{,}...,\gamma_{k})=E(z^*(s,m_1))-E(s_{m_0})$. Then
\[\lim_{\beta\to\infty}\frac{1}{\beta}\ln\mathbb{P}_{m_0}(X_{\xi_1}\in V^{(i)}(m_1))\ 
=\ -\left(E(z^*(m_0,m_1))-E(s_{m_0})\right).\]
Without assuming the existence of $\gamma$ as stated, the result remains valid when replacing $=$ with $\le$.
\end{lem}

Note that $I(\gamma_0{,}... ,\gamma_{k})=E(z^*(s_{m_0},m_1))-E(s_{m_0})$ does indeed imply the already mentioned property that
\begin{align*}
E(\gamma_i)&>E(\gamma_{i-1})\quad \textrm{for $1\leq i\leq j$}\\
\text{and}\quad 
E(\gamma_i)&<E(\gamma_{i-1})\quad\textrm{for $j+1\leq i\leq k$}
\end{align*}
if $\gamma_j=z^*(s_{m_0},m_1)$. We call such a path an \textit{uphill-downhill-path} because it first straddles the energy barrier $E(z^{*}(s_{m_0},m_{1}))$ and then falls down to the local minimum $m_{1}$. The existence of such a path can be found in most 2- or higher dimensional energy landscapes.

\begin{bew}
With $\gamma$ as stated, a lower bound for $\mathbb{P}_{m_0}(X_{\xi_1}\in V^{(i)}(m_1))$ is easily obtained as follows:
\begin{align*}
\mathbb{P}_{m_0}(X_{\xi_1}\in V^{(i)}(m_1))\ 
&\geq\ \mathbb{P}_{m_0}(X_{\zeta_{0}+i}=\gamma_i,\, 0\leq i\leq k)\\
&\geq\ \mathbb{P}_{m_0}(X_{\zeta_{0}}=s_{m_0})\,e^{-\beta I(\gamma_0{,}..., \gamma_{k})-\gamma_{\beta}\beta|\mathcal{S}|}\\
&=\ (1+o(1))\,e^{-\beta\left(E\left(z^*\left(s_{m_0},m_1\right)\right)-E(s_{m_0})\right)-\gamma_{\beta}\beta|\mathcal{S}|}\\
&=\ (1+o(1))\,e^{-\beta\left(E\left(z^*\left(m_0,m_1\right)\right)-E(s_{m_0})\right)-\gamma_{\beta}\beta|\mathcal{S}|}.
\end{align*}

For an upper bound, which does not require the existence of a $\gamma$ as claimed, we decompose the event into disjoint sets depending on the number of visits $N$, say, to $m_0$ between $1$ and $\zeta_{0}=\zeta_{0}^{(i)}$. This leads to
\[\mathbb{P}_{m_0}(X_{\xi_1}\in V^{(i)}(m_1),N=0)=\mathbb{P}_{m_0}(\xi_1=\tau_{V^{(i)}(m_1)}<\tau_{m_0})\]
and, for $k\geq 1$,
\begin{align*}
\mathbb{P}_{m_0}&(X_{\xi_1}\in V^{(i)}(m_1),N=k)\\
&=\mathbb{P}_{m_0}(X_{\xi_1}\in V^{(i)}(m_1),\,|\{\tau_{m_0}< n\leq \zeta_0|X_n=m_0\}|=k-1,\,\tau_{m_0}<\zeta_0)\\
&=\mathbb{P}_{m_0}(X_{\xi_1}\in V^{(i)}(m_1),N=k-1)\,\mathbb{P}_{m_0}(\tau_{m_0}<\zeta_0)\\
&~~\vdots\\
&=\mathbb{P}_{m_0}(X_{\xi_1}\in V^{(i)}(m_1),N=0)\,\mathbb{P}_{m_0}(\tau_{m_0}<\zeta_0)^k\\
&=\mathbb{P}_{m_0}(\xi_1=\tau_{V^{(i)}(m_1)}<\tau_{m_0})\,\mathbb{P}_{m_0}(\tau_{m_0}<\zeta_0)^k.
\end{align*}
Consequently,
\begin{align*}
\mathbb{P}_{m_0}(X_{\xi_1}\in V^{(i)}(m_1))
&=\sum_{k\geq 0}\mathbb{P}_{m_0}(X_{\xi_1}\in V^{(i)}(m_1),N=k)\\
&=\sum_{k\geq 0}\mathbb{P}_{m_0}(\xi_1=\tau_{V^{(i)}(m_1)}<\tau_{m_0})\,\mathbb{P}_{m_0}(\tau_{m_0}<\zeta_0)^k\\
&=\frac{\mathbb{P}_{m_0}(\xi_1=\tau_{V^{(i)}(m_1)}<\tau_{m_0})}{\mathbb{P}_{m_0}(\zeta_0<\tau_{m_0})}.
\end{align*}
By invoking Proposition \ref{prop:AbschaetzungUeWkeiten}, we infer
\begin{align}\label{eqn:UpperBound}
\frac{\mathbb{P}_{m_0}(\xi_1=\tau_{V^{(i)}(m_1)}<\tau_{m_0})}{\mathbb{P}_{m_0}(\zeta_0<\tau_{m_0})}\ 
&\leq\ \frac{\sum_{r\in V^{(i)}(m_1)}\mathbb{P}_{m_0}(\tau_{r}<\tau_{m_0})}{\mathbb{P}_{m_0}(\tau_x<\tau_{m_0})}\nonumber \\
&\leq\ K(\beta)\sum_{r\in V^{(i)}(m_1)}e^{-\beta(E(z^*(m_0,r))-E(z^*(m_0,x))-7\gamma_{\beta})}
\end{align}
for all $x\in V^{(i)}(m_0)^c$, where
\begin{equation*}
K(\beta)\ 
=\ |\mathcal{S}|\left(|\mathcal{S}|\exp\left(-\beta\min_{a\neq b:E(a)>E(b)}(E(a)-E(b))+2\gamma_{\beta}\beta\right)+1\right)
\max_{r\in \mathcal{S}}|\mathcal{N}(r)|.
\end{equation*}
For any $r\in V^{(i)}(m_{1})$, we have $E(z^{*}(m_{0},r))\ge E(z^{*}(m_{1},r))$ and therefore
\[E(z^{*}(m_{0},m_{1}))\ \le\ E(z^{*}(m_{0},r))\vee E(z^{*}(r,m_{1}))\ =\ E(z^{*}(m_{0},r)).\]
Using this in  \eqref{eqn:UpperBound}, we obtain
\[\frac{\mathbb{P}_{m_0}(\xi^{(i)}_1=\tau_{V^{(i)}(m_1)}<\tau_{m_0})}{\mathbb{P}_{m_0}(\zeta_0<\tau_{m_0})}\ 
\leq\ K(\beta)\, |\mathcal{S}|\, e^{-\beta(E(z^*(m_0,m_1))-E(z^*(m_0,x))-7\gamma_{\beta})}\]
and then, upon choosing $x=s_{m_0}$ and noting that $E(z^{*}(m_{0},s_{m_0}))=E(s_{m_0})$,
\[\mathbb{P}_{m_0}(X_{\xi_1}\in V^{(i)}(m_1))\ 
\leq\ K(\beta)\,|\mathcal{S}|\,e^{-\beta\left(E\left(z^*\left(m_0,m_1\right)\right)-E(s_{m_0})-7\gamma_{\beta}\right)}.\]
By combining all previous results, we finally conclude
\[\lim_{\beta\to\infty}\frac{1}{\beta}\ln\mathbb{P}_{m_0}(X_{\xi_1}\in V^{(i)}(m_1))\ 
=\ -\left(E(z^*(m_0,m_1))-E(s_{m_0})\right)\]
as asserted.
\end{bew}

To summarize, which valleys are visited consecutively depends on (a) their spatial arrangement and (b) the energy barriers between them: A transition from one valley $V^{(i)}(m_0)$ to another valley $V^{(i)}(m_1)$ is only possible, if there exists a path from $s_{m_0}$ to $V^{(i)}(m_1)$, not hitting any other valley. This transition is made at small temperatures (i.e. large $\beta$) if the additional energy barrier $E(z^*(s_{m_0},m_1))-E(s_{m_0})$ is sufficiently small or in other words the energy barrier $E(z^*(s_{m_0},m_1))$ is approximately of the same height as all other energy barriers, including the barrier $E(z^*(s_{m_0},m_0))=E(s_{m_0})$. 

\vspace{.2cm}
A result similar to the previous lemma holds true for transitions from $m\in \mathcal{S}^{(i)}\backslash N^{(i)}$ to any $s\in \partial^+V(m)$.

\begin{lem}\label{lem:TransitionProbability2}
Let $m\in \mathcal{S}^{(i)}\backslash N^{(i)}$ and $s\in \partial^+V^{(i)}(m)$. Then
\[\lim_{\beta\to\infty}\frac{1}{\beta}\ln\Prob_m(Y_1=s)=-(E(s)-E(s_m)).\]
\end{lem}

\begin{bew}
For the proof, decompose again the event $\{Y_1=s\}$ with respect to the number of visits to $m$ before $V(m)=V^{(i)}(m)$ is left (or use  Proposition \ref{prop:Cor 1.6 Bovier}), giving
\[\Prob_m(Y_1=s)\ =\ \frac{\Prob_m(\sigma_1=\tau_s<\tau_m)}{\Prob_m(\sigma_1<\tau_m)}.\]
The proof of the lower bound is much more technical. Let $\gamma=(\gamma_1{,}...,\gamma_n)\in \Gamma^*(m,s)$ be a minimal path which leaves $V(m)$ only in the last step and such that for any $\gamma_i,\gamma_j\in \gamma$ both, the subpath from $\gamma_i$ to $\gamma_j$, and the inversed path from $\gamma_j$ to $\gamma_i$, are minimal. Define 
\[r_0:=m \quad \textrm{and}\quad r_1:=\gamma_{i_0}\quad \textrm{with}\quad i_0:=\inf\{0\leq i\leq n-1|E(\gamma_{i+1})\geq E(s_m)\}.\]
In particular, $E(r_1)<E(s_m)$ and $E(z^*(r_1,r_0))<E(s_m)$. Define furthermore the first record by $s_1:=\gamma_{i_1}$ with
\begin{align*}
i_1\ 
:=\ \inf\Big\{i_0< i\leq n\,\Big|\,&E(\gamma_i)\geq E(s_m),\,\inf\{j\geq i|E(\gamma_j)<E(\gamma_{i})\}<\inf\{j\geq i|E(\gamma_j)>E(\gamma_{i})\}\Big\},
\end{align*}
and then successively for $k\geq 1$ with $s_k=\gamma_{i_k}\neq s$ the records $s_{k+1}:=\gamma_{i_{k+1}}$ with
\begin{align*}
i_{k+1}\ 
:=\ \inf\Big\{n\geq i\geq\inf&\{j\geq i_k|E(\gamma_j)<E(s_k)\}\,\Big|\,E(\gamma_i)\geq E(s_k),\\
& \inf\{j\geq i|E(\gamma_j)<E(\gamma_{i})\}<\inf\{j\geq i|E(\gamma_j)>E(\gamma_{i})\}\Big\}.
\end{align*}
Note that the energy of these records is increasing. Let $s_{k-1}$ be the last record defined in this way and $s_k:=s$. Since $E(z^*(m,s))=E(s)$, $s_k$ is as well a record. Given the records $s_1{,}...,s_k$, for $1\leq i\leq k-1$ let $r_{2i}$ be the first minimum along $\gamma$ after $s_i$ and $r_{2i+1}$ the last minimum along $\gamma$ before $s_{i+1}$. Here a minimum along $\gamma$ is some $\gamma_i\in \gamma$ such that it is a minimum of $E$ restricted to $\gamma$. Finally, let $r_{2k}:=s_k=s$. In the following we will proof that
\begin{itemize}
\item[(a)] $\Prob_{r_{2j}}(\tau_{r_{2j+1}}<\zeta_0)\to 1$ as $\beta\to\infty$ for any $0\leq j\leq k-1$,
\item[(b)] $\Prob_{r_1}(\tau_{r_2}<\zeta_0)\geq e^{-\beta(E(z^*(r_1,r_2))-E(s_m)+o(1))}$,
\item[(c)] $\Prob_{r_{2j+1}}(\tau_{r_{2j+2}}<\zeta_0)\geq e^{-\beta(E(z^*(r_{2j+1},r_{2j+2}))-E(z^*(r_{2j-1},r_{2j}))+o(1))}$ for any $1\leq j\leq k-2$, 
\item[(d)] $\Prob_{r_{2k-1}}(\tau_{r_{2k}}=\zeta_0)\geq e^{-\beta(E(z^*(r_{2k-1},r_{2k}))-E(z^*(r_{2k-3},r_{2k-2}))+o(1))}$.
\end{itemize}
This gives for $\beta$ large enough
\begin{align*}
\Prob_m(Y_1=s)\ 
&\geq\ \left(\prod_{j=0}^{2k-2}\Prob_{r_j}(\tau_{r_{j+1}}<\zeta_0)\right)\Prob_{r_{2k-1}}(\tau_{r_{2k}}=\zeta_0)\\
&\geq\ \frac{1}{2}\left(\prod_{j=0}^{k-2}\Prob_{r_{2j+1}}(\tau_{r_{2j+2}}<\zeta_0)\right)\Prob_{r_{2k-1}}(\tau_{r_{2k}}=\zeta_0)\\
&\geq\ \frac{1}{2} e^{-\beta(E(z^*(r_1,r_2))-E(s_m)+o(1))}\cdot \prod_{j=1}^{k-1}e^{-\beta(E(z^*(r_{2j+1},r_{2j+2}))-E(z^*(r_{2j-1},r_{2j}))+o(1))}\\
&=\ \frac{1}{2}e^{-\beta(E(z^*(r_{2k-1},r_{2k}))-E(s_m)+o(1))}\\
&=\ e^{-\beta(E(s)-E(s_m)+o(1))},
\end{align*}
and thus the assertion.

(a) For $j=0$ this is obvious since $E(z^*(r_0,r_1))<E(s_m)$. For $1\leq j\leq k-1$ and any $r'\in \partial^+V(m)$ it holds true that
\begin{align}\label{eqn:LeaveValley1}
E(z^*(r_{2j},r_{2j+1}))\ 
<\ E(z^*(r_{2j},r_{2j-1}))\ 
\leq\ E(z^*(r_{2j},m))\ 
\leq\ E(z^*(r_{2j},r')),
\end{align}
where we make use of the fact that between $r_{2j}$ and $r_{2j+1}$ the energy stays below the last record, and that all subpaths of $\gamma$ are minimal as well.

(b) By the definition of $r_1$ and $r_2$, there is a unimodal path between them so that the cumulative activation energy along this path equals $E(z^*(r_1,r_2))-E(r_1)$. Therefore
\[\Prob_{r_1}(\tau_{r_2}<(\zeta_0\wedge \tau_{r_1}))\ \geq\ e^{-\beta(E(z^*(r_1,r_2))-E(r_1)+o(1))}.\]
Furthermore, for any $r'\in \partial^+V(m)$ we have 
\begin{align*}
E(z^*(r_1,r'))\ 
&\geq\ E(s_m)>E(r_1)
\quad \textrm{ and }\quad
E(z^*(r_1,r_2))\ 
\geq\ E(s_m)>E(r_1).
\end{align*}
Therefore,
\begin{align*}
\Prob_{r_1}((\tau_{r_2}\wedge \zeta_0)<\tau_{r_1})\ 
&\leq\ \sum_{r'\in \partial^+V(m)\cup\{r_2\}}\Prob_{r_1}(\tau_{r'}<\tau_{r_1})\\
&\leq\ \sum_{r'\in \partial^+V(m)\cup\{r_2\}}e^{-\beta(E(z^*(r_1,r'))-E(r_1)+o(1))}\\
&\leq\  e^{-\beta(E(s_m)-E(r_1)+o(1))}.
\end{align*}
Combining the two estimates, we get
\begin{align*}
\Prob_{r_1}(\tau_{r_2}<\zeta_0)\
&=\ \frac{\Prob_{r_1}(\tau_{r_2}<(\zeta_0\wedge \tau_{r_1}))}{\Prob_{r_1}((\tau_{r_2}\wedge \zeta_0)<\tau_{r_1})}\ 
\geq\ e^{-\beta(E(z^*(r_1,r_2))-E(s_m)+o(1))}.
\end{align*}

(c) Let $1\leq j\leq k-2$. We use the same strategy as in the proof of (b). So, again, by the definition of $r_{2j+1}$ and $r_{2j+2}$, there is a unimodal path between them with cumulative activation energy $E(z^*(r_{2j+1},r_{2j+2}))-E(r_{2j+1})$ along this path, and
\[\Prob_{r_{2j+1}}(\tau_{r_{2j+2}}<(\zeta_0\wedge \tau_{r_{2j+1}}))\ 
\geq\ e^{-\beta(E(z^*(r_{2j+1},r_{2j+2}))-E(r_{2j+1})+o(1))}.\]
Furthermore,
\[E(z^*(r_{2j-1},r_{2j}))\ 
\leq\ E(z^*(r_{2j+1},r_{2j+2}))
\quad \textrm{ and }\quad 
E(r_{2j+1})\ <\ E(z^*(r_{2j+1},r_{2j+2})).\]
Finally, for any $r'\in \partial^+V(m)$, by use of Equation \eqref{eqn:LeaveValley1},
\begin{align*}
E(r_{2j+1})\ 
&\leq E(z^*(r_{2j+1},r_{2j}))\\ 
&<\ E(z^*(r_{2j},r'))\\
&\leq\ E(z^*(r_{2j},r_{2j+1}))\vee E(z^*(r_{2j+1},r'))\\
&=\ E(z^*(r_{2j+1},r')).
\end{align*}
Thus,
\begin{align*}
\Prob_{r_{2j+1}}((\tau_{r_{2j+2}}\wedge \zeta_0)<\tau_{r_{2j+1}})\ 
&\leq\ \sum_{r'\in \partial^+V(m)\cup\{r_{2j+2}\}}\Prob_{r_{2j+1}}(\tau_{r'}<\tau_{r_{2j+1}})\\
&\leq\ \sum_{r'\in \partial^+V(m)\cup\{r_{2j+2}\}}e^{-\beta(E(z^*(r_{2j+1},r'))-E(r_{2j+1})+o(1))}\\
&\leq\  e^{-\beta(E(r_{2j-1},r_{2j})-E(r_{2j+1})+o(1))}.
\end{align*}
Combining the two estimates, we get
\begin{align*}
\Prob_{r_{2j+1}}(\tau_{r_{2j+2}}<\zeta_0)\
&=\ \frac{\Prob_{r_{2j+1}}(\tau_{r_{2j+2}}<(\zeta_0\wedge \tau_{r_{2j+1}}))}
{\Prob_{r_{2j+1}}((\tau_{r_{2j+2}}\wedge \zeta_0)<\tau_{r_{2j+1}})}\ 
\geq\ e^{-\beta(E(z^*(r_{2j+1},r_{2j+2}))-E(r_{2j-1},r_{2j})+o(1))}.
\end{align*}

(d) All bounds for energies in (c) can be proved in the very same way (here $r'\in \partial^+V(m)\backslash \{s\}$), so that
\begin{align*}
\Prob_{r_{2k-1}}(\tau_{r_{2k}}=\zeta_0)\
&=\ \frac{\Prob_{r_{2k-1}}(\tau_{r_{2k}}=(\zeta_0\wedge \tau_{r_{2k-1}}))}
{\Prob_{r_{2k-1}}((\tau_{r_{2k}}\wedge \zeta_0)<\tau_{r_{2k-1}})}\ 
\geq\ e^{-\beta(E(z^*(r_{2k-1},r_{2k}))-E(r_{2k-3},r_{2k-2})+o(1))}.
\end{align*}
\end{bew}

Now we see for the reciprocating jumps in the accelerated chain:

\begin{prop}\label{prop:ForthAndBack}
Fix $1\leq i\leq \mathfrak{n}$ and $\varepsilon>0$. Then the AAC at level $i$ exhibits no reciprocating jumps of order $\varepsilon$ if the following three conditions hold true:
\begin{itemize}
\item[(1)] $E(z^*(m_0,m_1))-E(s_{m_0})\leq\varepsilon$ for all distinct $m_0,m_1\in \mathcal{S}^{(i)}\backslash N^{(i)}$.
\item[(2)] For each $m\in\mathcal{S}^{(i)}\backslash N^{(i)}$, there exist at least two distinct $m_1,m_2\in \mathcal{S}^{(i)}\backslash N^{(i)}, m\neq m_1,m_2,$ such that $\Prob_m(X_{\xi_1}\in V^{(i)}(m_j))>0$ for $j=1,2$.
\item[(3)] For each pair $m_0,m_1\in \mathcal{S}^{(i)}\backslash N^{(i)}$ with $\Prob_{m_0}(X_{\xi_1}\in V^{(i)}(m_1))>0$, there exists a minimal uphill-downhill-path from $s_{m_0}$ to $m_1$ not hitting any valley but $V^{(i)}(m_1)$.
\end{itemize}
\end{prop}

The origin of our endeavor to define aggregations with no reciprocating jumps of an order larger than a small $\varepsilon$ is to obtain an associated process with (almost) decorrelated increments (in Euclidean state space), for this and a proper centering causes the variance up to the $n$-th jump to grow with $n$ instead of $n^2$. This is known as diffusive behavior in physics. Without aggregation increments are highly correlated due to the following argument: at any given time, the process is with high probability in a minimum and when leaving it, say by making a positive jump, the next increment is most likely negative because there is a drift back to the minimum. Likewise, the increments of the asymptotic jump chain are neither uncorrelated nor having mean zero since trajectories of $\wh{Y}$ on an irreducibility class are almost surely of the form $m_1\to s\to m_2\to s\to m_3\to...$, where $s\in N^{(i)}$ and $m_1,m_2{,}...\in \mathcal{S}^{(i)}\backslash N^{(i)}$. 
Thus, given the previous increments, it is in general easy to predict the next increment and they do not have mean zero. On the other hand, if we can choose $\beta$ and $i$ such that, for any $m_0,m_1$, we have $\Prob_{m_0}(Y_1=m_1)\in \{p_{m_0}\pm \varepsilon\}\cup [0,\varepsilon]$ for $\varepsilon\ll p_{m_0}$, we obtain an AAC which behaves roughly like a RW on a graph. Such a RW is diffusive if we assume periodic boundary conditions (or sufficiently large state space compared to the observation time $n$) and an energy landscape $E$ which is homogeneous enough to ensure zero-mean increments. In particular $\{m|\Prob_{m_0}(Y_1=m)=p_{m_0}\pm\varepsilon\}$ has to comprise at least two states.

\subsection{Metabasins}\label{subsec:Metabasins}

A path-independent definition of metabasins can now be given on the basis of the previous considerations.

\begin{defn}\label{defn:MB}
A finite Markov chain $X$ driven by an energy function $E$ satisfying the assumptions stated at the beginning of Section \ref{sec:valleys} has \emph{metabasins of order $\varepsilon>0$} if there exists an aggregation level $i<\mathfrak{n}-1$ such that the following conditions are fulfilled for each $m\in \mathcal{S}^{(i)}\backslash N^{(i)}$:
\begin{mydescription}{(MB4)}
\item[(MB1)] $\sup_{m'\in \mathcal{S}^{(i)}\backslash(N^{(i)}\cup\{m\})}E(z^*(m,m'))-E(s_m)\leq \varepsilon$.
\item[(MB2)] There are at least two distinct $m_1,m_2\in \mathcal{S}^{(i)}\backslash(N^{(i)}\cup\{m\})$ with a minimal uphill-downhill-path from $s_m$ to $m_k$ not hitting any other valley but $V^{(i)}(m_{k})$ for $k=1,2$.
\end{mydescription}
In this case, the valleys $(V^{(i)}(m))_{m\in \mathcal{S}^{(i)}}$ are called \textit{metabasins (MB) of order $\varepsilon$}.
\end{defn}

The reader should notice that each singleton set $\{s\}$ consisting of a non-assigned state $s\in\mathcal{N}^{(i)}$ forms a MB. The conditions (MB1) and (MB2) ensure the good nature of (a) the energy barriers and (b) the spatial arrangement of minima. As already pointed out, this determines which valleys are visited consecutively. Properties of MB which can be concluded from the results of the previous sections are summarized in the next theorem. The reader is reminded of Properties 1--4 stated in the Introduction.

\begin{thm}\label{thm:MB}
For MB as defined in Definition \ref{defn:MB} we have
\begin{mydescription}{(4)}
\item[(1)] The transition probabilities for jumps between MB do not depend on the point of entrance as $\beta\to\infty$ (Property 1).
\item[(2)] There are no reciprocating jumps of order $\varepsilon$ (Property 2). 
\item[(3)] The expected residence time in a MB depends on $E$ only via the depth of the MB as $\beta\to\infty$ (Property 3).
\item[(4)] Regarding only MB pertaining to local minima, the system is a trap model (Property 4).
\end{mydescription}
\end{thm}

\begin{bew}
(1) follows from Proposition \ref{prop:Übergänge}, (2) from Proposition \ref{prop:ForthAndBack}, (3) from Theorem \ref{thm:Verlassenszeit}, and (4) directly from the definition.
\end{bew}

It should not be surprising that the path-dependent definition of MB by \textsc{Heuer} \cite{He08} and stated in the Introduction differs from our path-independent one.
\begin{figure}[ht]
\begin{center}
\includegraphics[trim=5.5cm 23cm 5.5cm 3.5cm, clip, width = 10cm]{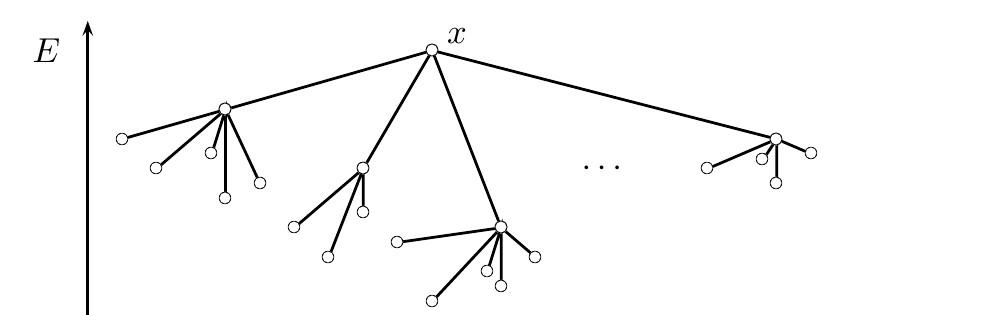}
\end{center}
\caption{Example of an energy landscape with a tree-like structure.}
\label{fig:PEL-Baum}
\end{figure}
For example, the energy landscape in Figure \ref{fig:PEL-Baum} has no reasonable path-dependent MB because every transition between two branches of the shown tree must pass through the state $x$. For a typical trajectory, there will be at most three MB: the states visited before the first occurrence of $x$, the states visited between the first and the last occurrence of $x$, and the states visited after the last occurrence of $x$. The reason for this poor performance is the tree-like structure of the energy landscape or, more generally, the fact that the connectivity between the branches is too small to allow a selfavoiding walk through more than two branches. This results in a small recurrence time for $x$ (compared to the number of states visited in between). However, every branch constitutes a MB when using the path-independent definition for sufficiently small $\varepsilon$, in which case the AAC forms a Markov chain and, given the Metropolis algorithm, even a RW on the graph.

Having thus exemplified that the two definitions of MB do not necessarily coincide, where the path-independent approach applies to a wider class of energy landscapes, we turn to the question about conditions for them to coincide with high probability. As already pointed out, we have to assume a sufficient connectivity between the metastates to ensure the existence of reasonable path-dependent MB. In terms of this connectivity (for a precise definition see Definition \ref{def:ConnectivityParameter}) and the parameter $\beta$ and $\varepsilon$, our last result, Theorem \ref{thm:pd versus pid MB} below, provides lower bounds for the probability that both definitions yield the same partition of the state space.

The first step towards this end is to identify and count, for each $m\in \mathcal{S}^{(i)}$ and a given $\beta$, the states $s\in \mathcal{S}^{(i)}$ for which a transition of $Y$ from $m$ to $s$ is likely. This leads to the announced connectivity parameter.

\begin{defn}\label{def:ConnectivityParameter}
Let $\varepsilon>0$ and suppose that $X$ has MB of order $\varepsilon>0$ at level $i$. Define the \textit{connectivity parameters}
\begin{align}\label{eqn:ConnectivityParameter}
\eta_1\ &=\ \eta_{1,\varepsilon}\ 
:=\ \min_{n\in N^{(i)},r\in\mathcal{S}^{(i)}:V^{(i)}(r)\cap \mathcal{N}(n)\neq \emptyset}
\left|\left\{s\in \mathcal{N}(n)\backslash V^{(i)}(r)\big|\,E(s)\le E(n)+\varepsilon\right\}\right|,\nonumber\\
\eta_2\ &=\ \eta_{2,\varepsilon}\ 
:=\ \min_{m\in \mathcal{S}^{(i)}\backslash N^{(i)}}\big|\{s\in \partial^+V^{(i)}(m)|E(s)\leq E(s_m)+\varepsilon\}\big|,\\
\eta_3\ &=\ \eta_{3,\varepsilon}\ 
:=\ \min_{n\in N^{(i)}}\big|\{s\in \mathcal{S}^{(i)}| E(x)\leq E(n)+\varepsilon 
\textrm{ for some $x\in V^{(i)}(s)\cap\mathcal{N}(n)$}\}\big|.\nonumber
\end{align}
\end{defn}

$\eta_1$ is the minimal number of neighboring sites of a non-assigned state $n$ which do not belong to a particular neighboring valley and whose energy is at most $\varepsilon$ plus the energy of $n$. $\eta_2$ is the  minimal number of neighboring sites/valleys of a valley $V^{(i)}(m)$ whose energy is at most $\varepsilon$ plus the energy of $s_m$ and which can be reached via an uphill-path from $m$. Finally, $\eta_3$ is the minimal number of neighboring valleys of a non-assigned state $n$ which comprise a state with energy of at most $\varepsilon$ plus the energy of $n$. $\eta_1$ and $\eta_3$ are always at least 2 and in fact quite large in the very complex energy landscapes of structural glasses. For very small $\varepsilon$, $\eta_2$ may be 1, but if $X$ has MB of order $\varepsilon$ in a high dimensional energy landscape, then $\eta_2$ can be assumed to be quite large as well.

That transitions to states counted above have reasonable large probabilities is content of the following lemma. Thus, the defined parameters do in fact measure the connectivity of the MB. 

\begin{lem}\label{lem:TransitionProbability3}
Let $\varepsilon>0$ and suppose that $X$ has MB of order $\varepsilon>0$ at level $i$ with connectivity parameters defined in \eqref{eqn:ConnectivityParameter}. Writing $Y_{k}$ for $Y_{k}^{(i)}$ and $V(m)$ for $V^{(i)}(m)$, $m\in\mathcal{S}^{(i)}$, the following assertions hold true for all sufficiently large $\beta$:
\begin{itemize}
\item[(a)] If $m\in\mathcal{S}^{(i)}\backslash N^{(i)}$ and $s\in\partial^{+}V(m)\cap\{x|E(x)-E(s_{m})\le\varepsilon\}$, or $m\in N^{(i)}$ and $s\in\mathcal{S}^{(i)}$ satisfies $V(s)\cap\{x\in\mathcal{N}(m)|E(x)-E(m)\le\varepsilon\}\ne\emptyset$, then
\[\Prob_{m}(Y_{1}=s)\ \ge\ e^{-2\beta\varepsilon}.\]
\item[(b)] For any distinct $m\in N^{(i)}$ and $s\in\mathcal{S}^{(i)}$,
\begin{equation*}
\Prob_{m}(Y_1\ne s)\ \ge\ \eta_1\,e^{-2\beta\varepsilon},
\end{equation*}
\item[(c)] For any distinct $m\in \mathcal{S}^{(i)}\backslash N^{(i)}$ and $s\in \mathcal{S}^{(i)}$,
\begin{equation*}
\Prob_{m}(Y_1\ne s)\ \ge\ (\eta_2-1)\,e^{-2\beta\varepsilon},
\end{equation*}
\end{itemize}
\end{lem}

We see that, for $\varepsilon$ small enough compared to $\beta$, transitions with an energy barrier of at most $\varepsilon$ are still quite likely and thus a jump to a particular valley quite unlikely in the case of high connectivity.

\begin{bew}
(a) Choose $\beta_{0}>0$ so large that, for $\beta\ge\beta_{0}$, $\gamma_{\beta}\leq \varepsilon$ and $\Prob_m(Y_1=s)\geq e^{-2\beta(E(s)-E(s_m))}$ for any $m\in\mathcal{S}^{(i)}\backslash N^{(i)}$ and $s\in\partial^{+}V(m)$, the latter being possible by Lemma \ref{lem:TransitionProbability2}. Then for any such $m$ and $s$, we infer $\Prob_{m}(Y_1=s)\geq e^{-2\varepsilon \beta}$ provided that additionally $E(s)\leq E(s_m)+\varepsilon$ holds true. If $m\in N^{(i)}$, then $\Prob_{m}(Y_1=s)\geq e^{-2\varepsilon\beta}$ for any $s\in \mathcal{S}^{(i)}$ such that $E(x)\leq E(s)+\varepsilon$ for some $x\in V(s)\cap\mathcal{N}(m)$, for
\[\Prob_{m}(Y_1=s)\ \ge\ \Prob_{m}(X_{\sigma_{1}}=x)\ =\ p(m,x)\ \ge\ e^{-\beta((E(x)-E(m))^{+}+\gamma_{\beta})}.\]
(b) Pick again $\beta_{0}$ so large that $\gamma_{\beta}\leq \varepsilon$ for all $\beta\ge\beta_{0}$. Then,
\begin{align*}
\Prob_{m}(Y_1\neq s)\ 
&\ge\ \sum_{x\in\mathcal{N}(m),x\notin V(s)}p(m,x)\\
&\geq\ \sum_{x\in\mathcal{N}(m),x\notin V(s), E(x)\le E(m)+\varepsilon}\exp\big(-\beta((E(x)-E(m))^{+}+\gamma_{\beta})\big)\\
&\geq\ \eta_1\exp(-2\beta\varepsilon),
\end{align*}
by definition of $\eta_1$.

\vspace{.2cm}
(c) Fix $\beta_{0}$ so large that $\Prob_m(Y_1=x)\geq e^{-2\beta(E(x)-E(s_m))}$ for any $x\in\partial^{+}V(m)$ and $\beta\ge\beta_{0}$. In the very same way as in part (b), we then get for all $\beta\geq \beta_0$
\begin{align*}
\Prob_{m}(Y_1\neq s)\
&\ge\ \sum_{x\in\partial^+V(m),x\neq s}\Prob_m(Y_1=x)\\
&\geq\ \sum_{x\in\partial^+V(m),x\neq s, E(x)\le E(s_m)+\varepsilon}\exp\big(-2\beta(E(x)-E(s_m))\big)\\
&\geq\ (\eta_2-1)\exp(-2\beta\varepsilon),
\end{align*}
by definition of $\eta_2$.
\end{bew}

The above result motivates that in the case of high connectivity the probability to revisit a particular valley within a fixed time $T$ is quite small, or in other words, the probability for the AAC to jump along a selfavoiding path is quite high. This is the main step towards the announced theorem and stated below. The observation time $T$ of course has to be small compared to the cover time of the process.

\begin{lem}\label{lem:return prob}
Let $\varepsilon>0$ and suppose that $X$ has MB of order $\varepsilon>0$ at level $i$ with connectivity parameters defined in \eqref{eqn:ConnectivityParameter}. Writing $Y_{k}$ for $Y_{k}^{(i)}$ and $V(m)$ for $V^{(i)}(m)$, $m\in\mathcal{S}^{(i)}$, define
\[\tau_{V(m)}^{(i)}:=\inf\{k\ge 1|Y_{k}=m\}.\]
Then the following assertions hold true for all sufficiently large $\beta$:
\begin{itemize}
\item[(a)] For any $0<\delta<1-\Prob_{m}(Y_{2}=m)$ and $1\leq T\leq T(m,\beta)+1$,
\[\Prob_{m}\left(\tau_{V(m)}^{(i)}>T\right)\ \geq\ \delta,\] 
where 
\[T(m,\beta)\ :=\ \frac{\ln\delta}{\ln\!\big(\min_{m'\neq m}\Prob_{m'}(Y_1\neq m)(1-\mathds{1}_{\{m'\notin N^{(i)}\}}\,\delta(m',\beta))\big)}\]
and
\[\delta(m',\beta)\ :=\ \max_{x\in V(m')}\sum_{z\in \partial^+ V(m')}\tilde{\varepsilon}(x,m',z,\beta).\]
In particular, if $\delta\leq\left((\eta_1\wedge(\eta_2-2))e^{-2\beta\varepsilon}\right)^{T}$ for some $T>0$, then $T(m,\beta)\geq T$.
\item[(b)] For each $k\ge 1$ and $m_{0}\in\mathcal{S}^{(i)}$,
\[\sum_{m_{1}{,}...,m_{k}}\prod_{j=0}^{k-1}\Prob_{m_j}(Y_1=m_{j+1})\ \geq\ [\eta_2\wedge\eta_3]_{k}\,e^{-2k\varepsilon\beta}\]
where summation ranges over all pairwise distinct $m_{1}{,}...,m_{k}\in\mathcal{S}^{(i)}\backslash\{m_0\}$ and for $N\in \mathbb{N}$ we write $[N]_{k}:=N(N-1)\cdot...\cdot (N-k+1)$.
\end{itemize}
\end{lem}

It should be noticed that $\Prob_{m}(\tau_{V(m)}^{(i)}>1)=1$ (the AAC never stays put) and
\begin{equation*}
\Prob_{m}(\tau_{V(m)}^{(i)}>T)\ 
\leq\ \Prob_{m}(\tau_{V(m)}^{(i)}>2)\ 
=\ 1-\Prob_{m}(Y^{(i)}_2=m)
\end{equation*}
for every $T\geq 2$ with equality holding only if $T=2$. We thus see that $\Prob_{m}(\tau_{V(m)}^{(i)}>T)\geq \delta$ entails $\delta<1-\Prob_{m}(Y_2=m)$, the latter being typically large. Furthermore, the bound on the number of self-avoiding path of length $k$ is very crude and can be improved when knowing more about the spatial arrangement of the metastable states. 

\begin{bew}
(a) Recall from the first part of the proof of Proposition \ref{prop:Übergänge} that
\begin{align*}
\Prob_{m}(Y_{n+1}\neq z|Y_n=y,X_{\sigma_n}=x)\ 
&\geq\ \Prob_{y}(Y_1\neq z)\left(1-\mathds{1}_{\{y\notin N^{(i)}\}}\sum_{r\in\partial^{+}V(y)}\Prob_{x}(\tau_{r}<\tau_{y})\right)\\
&\geq\ \Prob_{y}(Y_1\neq z)\left(1-\mathds{1}_{\{y\notin N^{(i)}\}} \sum_{r\in\partial^+V(y)}\tilde{\varepsilon}(x,y,r,\beta)\right)\\
&\geq\ \Prob_{y}(Y_1\neq z)\left(1-\mathds{1}_{\{y\notin N^{(i)}\}}\,\delta(y,\beta)\right)\\
\end{align*}
holds true for all $y,z\in\mathcal{S}^{(i)},\,x\in V(y)$ and $\beta>0$. This will now be used repeatedly to show that
\begin{equation*}
\Prob_{m}\Big(\tau_{V(m)}^{(i)}>T\Big)\ 
\geq\ \left(\min_{m'\neq m}\Prob_{m'}(Y_1\neq m)(1-\mathds{1}_{\{m'\notin N^{(i)}\}}\,\delta(m',\beta))\right)^{T-1}
\end{equation*}
for each $T>2$. Putting $\mathfrak{m}(x):=m'$ if $x\in V(m')$ for $m'\in\mathcal{S}^{(i)}$, we obtain
\begin{align*}
\Prob_{m}&\Big(\tau_{V(m)}^{(i)}>T\Big)\\ 
&=\ \Prob_{m}(Y_{1}\neq m{,}...,Y_{T}\neq m)\\
&=\ \sum_{x_1{,}...,x_{T-1}\notin V^{(i)}(m)}\Prob_{m}(X_{\sigma_1}=x_1)
\prod_{k=1}^{T-2}\Prob_{m}\big(X_{\sigma_{k+1}}=x_{k+1}|X_{\sigma_{k}}=x_{k}, Y_{k}=\mathfrak{m}(x_k)\big)\\
&\hspace*{2.7cm}\times\Prob_{m}(Y_{T}\neq m|X_{\sigma_{T-1}}=x_{T-1}, Y_{T-1}=\mathfrak{m}(x_{T-1}))\\
&\ge\ \sum_{x_1{,}...,x_{T-1}\notin V^{(i)}(m)}\Prob_{m}(X_{\sigma_1}=x_1)
\prod_{k=1}^{T-2}\Prob_{m}\big(X_{\sigma_{k+1}}=x_{k+1}|X_{\sigma_{k}}=x_{k}, Y_{k}=\mathfrak{m}(x_k)\big)\\
&\hspace{1.7cm}\times\left(\min_{m'\neq m}\left(\Prob_{m'}(Y_1\neq m)(1-\mathds{1}_{\{m'\notin N^{(i)}\}}\,\delta(m',\beta))\right)\right)\\
&\hspace{0.2cm}\vdots\\
&\geq\ \min_{m'\neq m}\Big(\Prob_{m'}(Y_1\neq m)(1-\mathds{1}_{\{m'\notin N^{(i)}\}}\,\delta(m',\beta))\Big)^{T-1}.\\
\end{align*}
But this establishes the asserted inequality when finally observing that the last expression is $\ge\delta$ iff $T\le T(m,\beta)+1$.

Having just said that $T(m,\beta)\ge T$ holds iff 
\[\min_{m'\neq m}\Big(\Prob_{m'}(Y_1\neq m)(1-\mathds{1}_{\{m'\notin N^{(i)}\}}\,\delta(m',\beta))\Big)^{T}\ \ge\ \delta,\]
it suffices to note that, as $\beta\to\infty$, $\delta(m',\beta)\to 0$ holds true if $m'\in \mathcal{S}^{(i)}\backslash N^{(i)}$, giving
\[1-\delta(m',\beta)\ \geq\ \frac{\eta_2-2}{\eta_2-1}\]
for sufficiently large $\beta$. Together with Lemma \ref{lem:TransitionProbability3}(b), this further yields
\begin{align*}
\min_{m'\neq m}\Big(\Prob_{m'}(Y_1\neq m)(1-\mathds{1}_{\{m'\notin N^{(i)}\}}\,\delta(m',\beta))\Big)^{T}
&\geq\ \left((\eta_1\wedge(\eta_2-2))e^{-2\beta\varepsilon}\right)^{T}
\end{align*}
and then the assertion.

\vspace{.2cm}
(c) Here it suffices to notice that, by (a), $[\eta_2\wedge\eta_3]_{k}$ forms a lower bound for the number of self-avoiding paths $(m_{0}{,}...,m_{k})$ such that $\Prob_{m_{j}}(Y_{1}=m_{j+1})\ge e^{-2\beta\varepsilon}$ for each $j=0{,}...,k-1$.
\end{bew}

We proceed with the announced result about the relation between path-dependent and path-independent MB. To this end, we fix $T=\sigma_{K}$ for some $K\in\mathbb{N}$. Let $\mathcal{V}_{k}$ for $k=1{,}...,\upsilon$ denote the random number of MB obtained from $X_{0}{,}...,X_{T}$ as defined in the Introduction. For $x\in\mathcal{S}$, we further let $\mathcal{V}(x)$ denote the MB $\mathcal{V}_{k}$ containing $x$ and put $\mathcal{V}(x):=\emptyset$ if no such MB exists which is the case iff $x\notin\{X_{0}{,}...,X_{T}\}$.

\begin{thm}\label{thm:pd versus pid MB}
Let $\varepsilon>0$ and suppose that $X$ has MB of order $\varepsilon>0$ at level $i$ with connectivity parameters defined in \eqref{eqn:ConnectivityParameter}. Fix $K\in\mathbb{N}$, $T=\sigma_{K}$ and $0<\delta \leq \left((\eta_1\wedge(\eta_2-1)-1)e^{-2\beta\varepsilon}\right)^{K}$. Then, for each $0\leq k<K$ and $m_0\in \mathcal{S}^{(i)}$, there exists $\beta_0>0$ such that for all $\beta\geq \beta_0$
\begin{itemize}
\item[(a)] $\Prob_{m_0}\!\left(V^{(i)}_<(Y_{k})\subseteq \mathcal{V}(Y_{k})\right)\geq 1-(\max_{m\in \mathcal{S}^{(i)}\backslash N^{(i)}}|V_<(m)|+2)\max_{m\in \mathcal{S}^{(i)}\backslash N^{(i)}}\delta(m,\beta)$, where $V_{<}^{(i)}(s):=\{s\}$ if $s\in N^{(i)}$.
\item[(b)] $\Prob_{m_0}(\mathcal{V}(Y_{j})\subseteq V^{(i)}(Y_{j}),\,0\le j<k)\geq 1-k(\max_{m\in \mathcal{S}^{(i)}\backslash N^{(i)}}\delta(m,\beta)+(1-\delta))$.
\item[(c)] If $\eta_2\wedge\eta_3>K-1$, then
\[\Prob_{m_0}(\mathcal{V}(Y_j)\subseteq V^{(i)}(Y_j),\,0\leq j\leq K-1)\ \geq\ [\eta_2\wedge\eta_3]_{K}\left(1-\max_{m\in \mathcal{S}^{(i)}\backslash N^{(i)}}\delta(m,\beta)\right)^{K-1}e^{-2K\varepsilon\beta}.\]
\end{itemize}
\end{thm}

For the occurring bounds to be significant, two requirements must be met. First, $K$ must be small compared to the cover time of the AAC and $\varepsilon$ must be small compared to $\beta_0$ to ensure $\exp(-2\beta\varepsilon)\gg 0$. Second, the connectivity must be high to ensure $1-\delta\ll 1$ and $[\eta_2\wedge\eta_3]_{K}e^{-2K\varepsilon\beta}\gg 0$. 

\vspace{.2cm}
Typically, the inclusions in parts (b) and (c) are strict because of high energy states within a valley that will probably be missed during one simulation run and therefore not belong to any path-dependent MB. On the other hand, since our approach strives to cover the state space as completely as possible by valleys the latter comprise such high energy states whenever they are assignable in the sense described in Section \ref{sec:valleys}. 

\begin{bew}
With $i$ being fixed, let us write as earlier $V(m)$ for $V^{(i)}(m)$, and also $V_{<}(m)$ for $V_{<}^{(i)}(m)$.

(a) For a given $0\leq k<K$, define
\begin{align*}
A_k&:=\{\sigma_k\leq \tau_{Y_k}<\sigma_{k+1}\},\\
B_k&:=\{\textrm{ for every $x\in V_{<}(m)$ exists $\tau_{Y_k}\leq l_x<\sigma_k$ such that $X_l=x$}\},\\
C_k&:=\{X_l=Y_k \textrm{ for some $\max_{x\in V_<(Y_k)}\tau_x\leq l<\sigma_k$}\}.
\end{align*}
With $\delta(m,\beta)$ as defined in Lemma \ref{lem:return prob} and using
\begin{align}\label{eqn:tau_m<sigma_1}
\Prob_{x}(\sigma_1<\tau_m)\ 
&\leq\ \sum_{y\in \partial^+V(m)}\tilde{\varepsilon}(x,m,y,\beta)\ 
\leq\ \max_{m\in\mathcal{S}^{(i)}\backslash N^{(i)}}\delta(m,\beta)\ =:\ \delta_{\max}
\end{align}
for $x\in V(m),\,m\in \mathcal{S}^{(i)}\backslash N^{(i)}$, we obtain
\begin{align}\label{eqn:thm:pd versus pid MB(a)}
\Prob_{m_0}\!&\left(V^{(i)}_<(Y_{k})\subseteq \mathcal{V}(Y_{k})\right)\nonumber\\
&\geq\ \Prob_{m_0}(A_k\cap B_k\cap C_k)\nonumber\\
&=\ \sum_{m\in \mathcal{S}^{(i)},r\in V(m)}\Prob_{m_0}(\{X_{\sigma_k}=r\}\cap A_k\cap B_k\cap C_k)\nonumber\\
&=\ \sum_{m\in \mathcal{S}^{(i)},r\in V(m)}\Prob_{m_0}(X_{\sigma_k}=r)\Prob_{m_0}(A_k|X_{\sigma_k}=r)
\Prob_{m_0}(B_k\cap C_k|\{X_{\sigma_k}=r\}\cap A_k)\nonumber\\
&=\ \sum_{m\in \mathcal{S}^{(i)},r\in V(m)}\Prob_{m_0}(X_{\sigma_k}=r)\Prob_r(\tau_m<\sigma_1)\nonumber\\
&\hspace{1cm}\times
\Prob_m(\tau_x<\sigma_1 \textrm{ for every $x\in V_<(m)$, $X_l=m$ for some $\max_{x\in V_<(m)}\tau_x\leq l<\sigma_1$})\nonumber\\
&\geq\ \sum_{m\in \mathcal{S}^{(i)},r\in V(m)}\Prob_{m_0}(X_{\sigma_k}=r)
(1-\delta_{\max})\nonumber\\
&\hspace{1cm}\times
\Prob_m(\tau_x<\sigma_1 \textrm{ for every $x\in V_<(m)$, $X_l=m$ for some $\max_{x\in V_<(m)}\tau_x\leq l<\sigma_1$})\nonumber\\
&=\ (1-\delta_{\max})
\sum_{m\in \mathcal{S}^{(i)}}\Prob_{m_0}(Y_k=m)\nonumber\\
&\hspace{1cm}\times
\Prob_m(\tau_x<\sigma_1 \textrm{ for every $x\in V_<(m)$, $X_l=m$ for some $\max_{x\in V_<(m)}\tau_x\leq l<\sigma_1$}).
\end{align}
Thus, in order to show that with high probability a path-dependent MB comprises the inner part of a valley, we show that with high probability, when starting in its minimum, the whole inner part will be visited and the process will return to the minimum once more before the valley is left. This is trivial if $m\in N^{(i)}$ and thus $V_{<}^{(i)}(m)=\{m\}$, for then
\[\Prob_m(\tau_x<\sigma_1 \textrm{ for every $x\in V_<(m)$, $X_l=m$ for some $\max_{x\in V_<(m)}\tau_x\leq l<\sigma_1$})=1.\]

More needs to be done if $m\in \mathcal{S}^{(i)}\backslash N^{(i)}$, where 
\begin{align*}
\Prob_{m}&\left(\tau_x<\sigma_1 \textrm{ for every $x\in V_<(m)$, $X_l=m$ for some $\max_{x\in V_<(m)}\tau_x\leq l<\sigma_1$}\right)\\
&\geq\ 1-\Prob_{m}\left(\tau_x> \sigma_1 \textrm{ for some $x\in V_<(m)$}\right)
-\Prob_m\big(\textrm{$X_l\neq m$ for each $\max_{x\in V_<(m)}\tau_x\leq l<\sigma_1$}\big).
\end{align*}
The second probability in the preceding line can further be bounded with the help of \eqref{eqn:tau_m<sigma_1}, viz.
\begin{align*}
\Prob_m&(\textrm{$X_l\neq m$ for each $\max_{x\in V_<(m)}\tau_x\leq l<\sigma_1$})\\
&=\ \sum_{y\in V_<(m)}\Prob_m(\max_{x\in V_<(m)}\tau_x=\tau_y,\textrm{ $X_l\neq m$ for every $\max_{x\in V_<(m)}\tau_x\leq l<\sigma_1$})\\
&\leq\ \sum_{y\in V_<(m)}\Prob_m(\max_{x\in V_<(m)}\tau_x=\tau_y)\Prob_y(\tau_m>\sigma_1)\\
&\leq\ \delta_{\max},
\end{align*}
while for the first probability, we obtain with the help of Theorem \ref{thm:DriftZumMinimum}
\begin{align}\label{eqn:return prob}
\Prob_{m}\left(\tau_x> \sigma_1 \textrm{ for some $x\in V_<(m)$}\right)\ 
&\leq \sum_{x\in V_<(m)}\Prob_{m}(\sigma_1<\tau_x)\nonumber\\
&\leq \sum_{x\in V_<(m)}\sum_{y\in \partial^+V(m)}\Prob_m(\tau_y<\tau_x)\nonumber\\
&\leq \sum_{x\in V_<(m)}\sum_{y\in \partial^+V(m)}\varepsilon(m,x,y,\beta),
\end{align}
because $E(z^*(m,y))>E(z^*(m,x))$ for $x\in V_<(m)$ and $y\in \partial^+V(m)$. The latter can be seen as follows: It has been shown in the proof of Theorem \ref{thm:DriftZumMinimumTäler} that $E(z^*(x,m))<E(z^*(x,y))$. Hence,
\[E(z^*(x,m))<E(z^*(x,y))\leq E(z^*(x,m))\vee E(z^*(y,m))=E(z^{*}(y,m))\]
as asserted. Next, we infer
\[E(z^*(x,y))\leq E(z^*(x,m))\vee E(z^*(m,y))\leq E(z^*(x,m))\vee E(z^*(x,y))=E(z^*(x,y)),\]
thus $E(z^*(x,y))=E(z^*(m,y))$. Recalling the definition of $\varepsilon(m,x,y,\beta)$, the last equality implies $\varepsilon(m,x,y,\beta)=\varepsilon(x,m,y,\beta)$ which will now be used to further bound the expression in \eqref{eqn:return prob}, namely
\begin{align*}
\sum_{x\in V_<(m)}\sum_{y\in \partial^+V(m)}\varepsilon(m,x,y,\beta)\ 
&=\ \sum_{x\in V_<(m)}\sum_{y\in \partial^+V(m)}\varepsilon(x,m,y,\beta)\\
&=\ \sum_{x\in V_<(m)}\sum_{y\in \partial^+V(m)}\tilde{\varepsilon}(x,m,y,\beta)\\
&\leq\ |V_<(m)|\,\delta_{\max}\\
&\leq\ \max_{m\in \mathcal{S}^{(i)}\backslash N^{(i)}}|V_<(m)|\,\delta_{\max}.
\end{align*}
Together with \eqref{eqn:thm:pd versus pid MB(a)} this yields as asserted
\begin{align*}
\Prob_{m_0}\!&\left(V^{(i)}_<(Y_{k})\subseteq \mathcal{V}(Y_{k})\right)\\
&\geq\ (1-\delta_{\max})
\sum_{m\in \mathcal{S}^{(i)}}\Prob_{m_0}(Y_k=m)\\
&\hspace{2cm}\times
\Prob_m(\tau_x<\sigma_1 \textrm{ for every $x\in V_<(m)$, $X_l=m$ for some $\max_{x\in V_<(m)}\tau_x\leq l<\sigma_1$})\\
&\geq\ \left(1-\delta_{\max}\right)
\left(1-\left(\max_{m\in \mathcal{S}^{(i)}\backslash N^{(i)}}|V_<(m)|+1\right)\delta_{\max}\right)\\
&\geq\ 1-\left(\max_{m\in \mathcal{S}^{(i)}\backslash N^{(i)}}|V_<(m)|+2\right)\delta_{\max}.
\end{align*}

(b) According to Lemma \ref{lem:return prob}, choose $\beta_{0}>0$ such that
\begin{equation}\label{cond:pd versus pid MB}
\max_{m\in\mathcal{S}^{(i)}}\Prob_{m}\left(\tau_{V(m)}^{(i)}\leq K\right)\ \leq\  1-\delta
\end{equation}
for each $\beta\ge\beta_{0}$. By using \eqref{cond:pd versus pid MB} and \eqref{eqn:tau_m<sigma_1}, we now infer
\begin{align*}
\Prob_{m_0}\!&\left(Y_l= Y_k\textrm{ for some }k+1\leq l\leq K\right)\\
&=\ \sum_{s\in\mathcal{S}}\Prob_{m_0}\!\left(Y_l= Y_k\textrm{ for some }k+1\leq l\leq K,X_{\sigma_{k}}=s\right)\\
&\leq\ \sum_{s\in\mathcal{S}}\Prob_{m_0}(X_{\sigma_{k}}=s)
\Big(\Prob_{s}(\tau_{V(Y_0)}^{(i)}\leq K-k,\,X_{j}=\mathfrak{m}(s) \text{ for some }0\leq j< \sigma_1)\\
&\hspace*{1.7cm}+\mathds{1}_{\{s\notin N^{(i)}\}}\,
\Prob_{s}(\tau_{V(Y_0)}^{(i)}\leq K-k, X_{j}\neq\mathfrak{m}(s) \text{ for all }0\leq j< \sigma_1)\Big)\\
&\leq\ \sum_{s\in\mathcal{S}}\Prob_{m_0}(X_{\sigma_{k}}=s)\,\Prob_{\mathfrak{m}(s)}(\tau_{V(Y_0)}^{(i)}\leq K-k)
+\sum_{s\notin N^{(i)}}\Prob_{m_0}(X_{\sigma_{k}}=s)\,\Prob_{s}(\sigma_1<\tau_{\mathfrak{m}(s)})\\
&\leq\ 1-\delta+\delta_{\max},
\end{align*}
and finally
\begin{align*}
\Prob_{m_0}(\mathcal{V}(Y_j)\subset V^{(i)}(Y_j),\,0\leq j< k)\ 
&\geq\ \Prob_{m_0}\left(\bigcap_{j=0}^{k-1}\left\{Y_l\neq Y_j,j+1\leq l\leq K\right\}\right)\\
&\geq\ 1-\sum_{j=0}^{k-1}\Prob_{m_0}(Y_l= Y_j\textrm{ for some }j+1\leq l\leq K)\\
&\geq\ 1-k(\delta_{\max}+(1-\delta)).
\end{align*}

(c) In the following calculation, let $r_{0}=m_{0}$, $\sum_{m_{j}}$ range over all $K$-vectors $(m_{1},...,m_{K})$ with pairwise distinct components in $\mathcal{S}^{(i)}\backslash\{m_0\}$ and, for each $k<K$, let $\sum_{r_{1}{,}...,r_{k}}$ range over all $k$-vectors $(r_{1}{,}...,r_{k})$ such that $r_{j}\in V(m_{j})$ for each $j=1{,}...,k$. As in part (b), use
\eqref{eqn:tau_m<sigma_1} repeatedly to infer
\begin{align*}
&\Prob_{m_0}(\mathcal{V}(Y_j)\subset V^{(i)}(Y_j),\,0\leq j\leq K-1)\\
&\geq\ \sum_{m_{j}}\sum_{r_{1}{,}...,r_{K-1}}\Prob_{m_0}\left(\bigcap_{j=0}^{K-1}\{Y_j=m_j, X_{\sigma_j}=r_j,\,
\tau_{m_j}<\sigma_{j+1}\}\cap \{Y_K=m_K\}\right)\\
&=\ \sum_{m_{j}}\sum_{r_{1}{,}...,r_{K-1}}\prod_{j=0}^{K-2}
\Prob_{r_{j}}\left(Y_0=m_j,X_{\sigma_{1}}=r_{j+1},\,\tau_{m_{j}}<\sigma_1\right)
\,\Prob_{r_{K-1}}\left(\tau_{m_{K-1}}<\sigma_1\right)\,\Prob_{m_{K-1}}(Y_1=m_K)\\
&\geq\ (1-\delta_{\max})\sum_{m_{j}}\sum_{r_{1}{,}...,r_{K-1}}\prod_{j=0}^{K-2}
\Prob_{r_{j}}\left(Y_0=m_j, X_{\sigma_{1}}=r_{j+1},\,\tau_{m_{j}}<\sigma_1\right)\Prob_{m_{K-1}}(Y_1=m_K)
\end{align*}
\begin{align*}
&=\ (1-\delta_{\max})\sum_{m_{j}}\sum_{r_{1}{,}...,r_{K-2}}\prod_{j=0}^{K-3}
\Prob_{r_{j}}\left(Y_0=m_j, X_{\sigma_{1}}=r_{j+1},\,\tau_{m_{j}}<\sigma_1\right)\\
&\hspace{2.8cm}\times\Prob_{r_{K-2}}\left(Y_0=m_{K-2}, Y_1=m_{K-1},\,\tau_{m_{K-2}}<\sigma_1\right)\,\Prob_{m_{K-1}}(Y_1=m_K)\\
&=\ (1-\delta_{\max})\sum_{m_{j}}\sum_{r_{1}{,}...,r_{K-2}}\prod_{j=0}^{K-3}
\Prob_{r_{j}}\left(Y_0=m_j, X_{\sigma_{1}}=r_{j+1},\,\tau_{m_{j}}<\sigma_1\right)\\
&\hspace{2.8cm}\times\Prob_{r_{K-2}}\left(Y_1=m_{K-1},\,\tau_{m_{K-2}}<\sigma_1\right)\,\Prob_{m_{K-1}}(Y_1=m_K)\\
&\geq\ (1-\delta_{\max})\sum_{m_{j}}\sum_{r_{1}{,}...,r_{K-2}}\prod_{j=0}^{K-3}
\Prob_{r_{j}}\left(Y_0=m_j, X_{\sigma_{1}}=r_{j+1},\,\tau_{m_{j}}<\sigma_1\right)\\
&\hspace{2.8cm}\times\Prob_{m_{K-2}}\left(Y_1=m_{K-1}\right)\,\Prob_{r_{K-2}}(\tau_{m_{K-2}}<\sigma_1)\,\Prob_{m_{K-1}}(Y_1=m_K)\\
&\geq\ (1-\delta_{\max})^{2}\sum_{m_{j}}\sum_{r_{1}{,}...,r_{K-2}}\prod_{j=0}^{K-3}
\Prob_{r_{j}}\left(Y_0=m_j, X_{\sigma_{1}}=r_{j+1},\,\tau_{m_{j}}<\sigma_1\right)\\
&\hspace{5cm}\times\Prob_{m_{K-2}}\left(Y_1=m_{K-1}\right)\,\Prob_{m_{K-1}}(Y_1=m_K)\\
&\hspace{.2cm}\vdots\\
&\geq\ (1-\delta_{\max})^{K-1}\sum_{m_{j}}\prod_{j=0}^{K-1}\Prob_{m_j}(Y_1=m_{j+1})\\
&\geq (1-\delta_{\max})^{K-1}\,[\eta_2\wedge\eta_3]_{K}\,e^{-2K\varepsilon\beta},
\end{align*}
the last line following from Lemma \ref{lem:return prob}.
\end{bew}

\begin{figure}[tb]
\begin{center}
\includegraphics[trim=5.5cm 21.9cm 8cm 3.5cm, clip, width = 7cm]{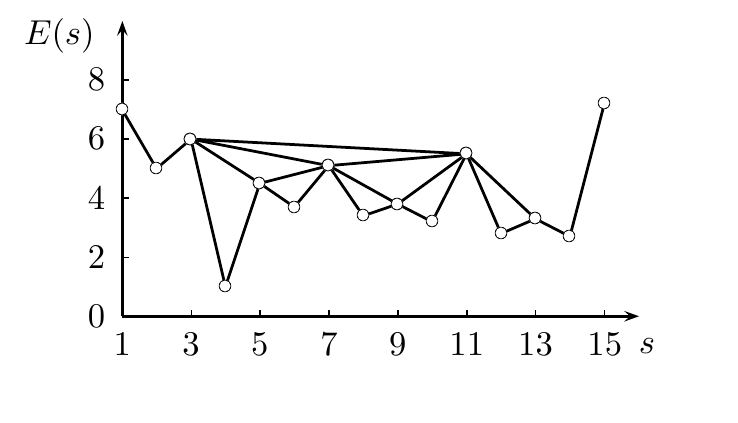}
\includegraphics[trim=6cm 21.9cm 8cm 3.5cm, clip, width = 7cm]{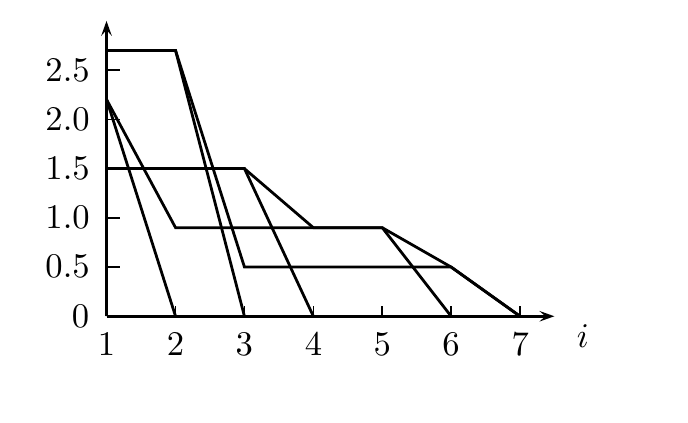}
\end{center}
\caption{(a) 2-dimensional modification of the energy landscape from Example \protect\ref{bsp:Energielandschaft}. \newline
(b) $\sup_{m'\in \mathcal{S}^{(i)}\backslash N^{(i)}}|E(z^*(m,m'))-E(s_m)|$ for the various metastable states in $\mathcal{S}^{(i)}\backslash N^{(i)}$ in dependence of the level $1\leq i\leq \mathfrak{n}$.}
\label{fig:Energiedifferenz}
\end{figure}

\begin{bsp}
We return to Example \ref{bsp:Energielandschaft} given in Section \ref{sec:valleys}, but modify the energy landscape by allowing direct transitions between some saddles (see Figure \ref{fig:Energiedifferenz} (a)) because (MB2) can clearly not be fulfilled in a one-dimensional model.
While having no effect on the metastable states $m\in M^{(i)}$, valleys change in the way that, for levels $i\in \{5,6\}$, the states $\{1,2,3\}$ do no longer belong to the valley around state 4 and $\{1,2\}$ forms its own valley. The energy-differences $\sup_{m'\in \mathcal{S}^{(i)}\backslash N^{(i)}}E(z^*(m,m'))-E(s_m)$ of the various metastable states $m$ at each level $1\leq i\leq \mathfrak{n}=7$ are shown in Figure \ref{fig:Energiedifferenz} (b). The supremum of these energy differences decreases in $i$, and we obtain MB of order 1 for $i\ge 4$, and of order $0.5$ for $i\ge 6$. 

To illustrate the behavior, we have run a Metropolis Algorithm on this energy landscape. For initial state $s=4$ and $\beta=0.75$, the energies of the trajectories of the original chain as well of the aggregated chain at levels $i=3,4,5$ are shown in Figure \ref{fig:Simulation}. The following observations are worth to be pointed out:
\begin{itemize}
\item The number of reciprocating jumps decreases with increasing level of aggregation.
\item The deeper the valley, the longer the residence time.
\item The motion in state space is well described by the aggregated process.
\item Due to the very small size of the state space and a long simulation time, valleys are revisited.
\end{itemize}
\begin{figure}[htb]
\begin{center}
\includegraphics[height = 14cm, width = 14cm]{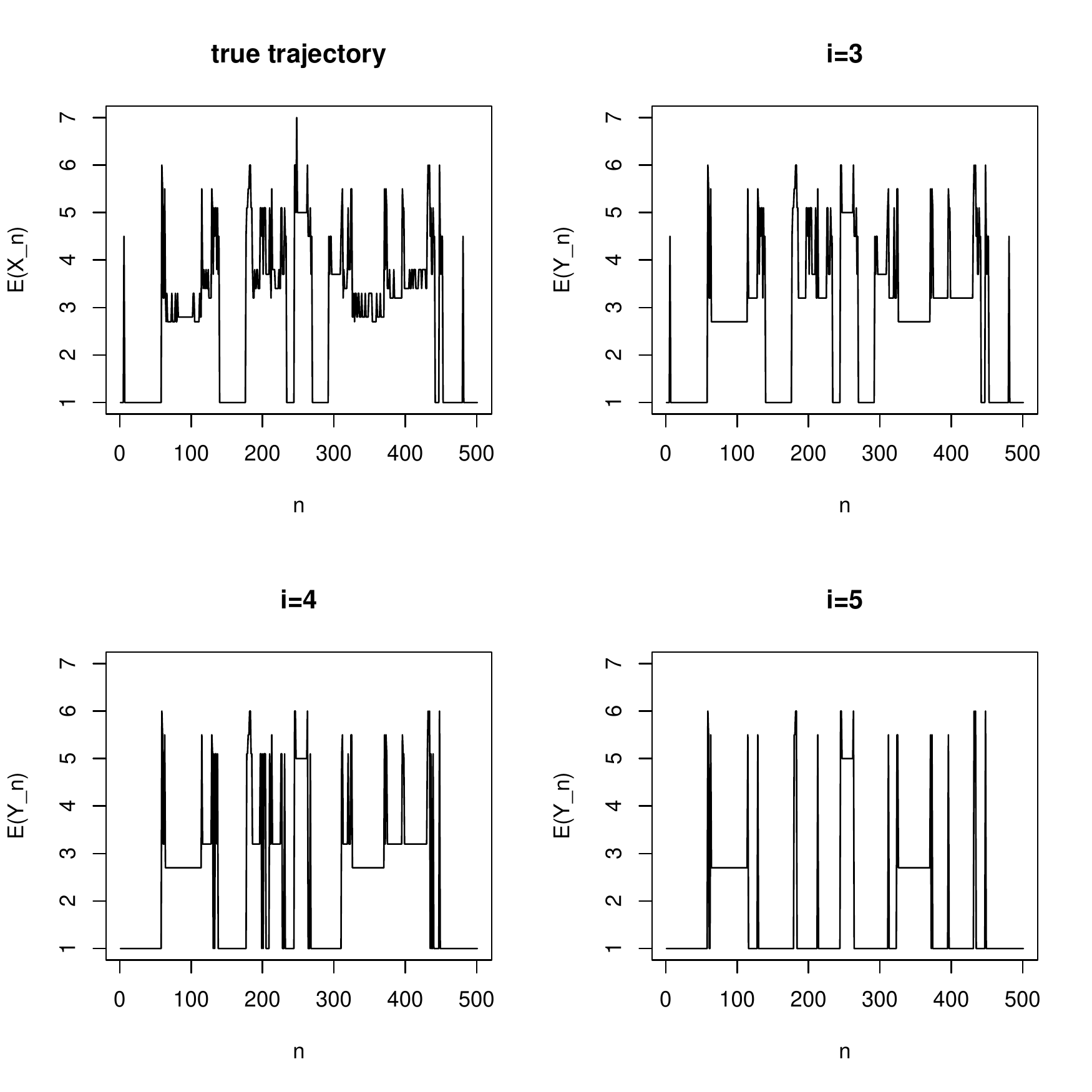}
\end{center}
\vspace{-.5cm}
\caption{Energies of the true trajectory and of the trajectories of the aggregated chain at levels $i=3,4,5$.}
\label{fig:Simulation}
\end{figure}
\end{bsp}

%
%

\section*{Acknowledgment}
\addcontentsline{toc}{section}{Acknowledgment}
\markboth{Acknowledgment}{Acknowledgment}

We are very indebted to Andreas Heuer for sharing his insight about glass forming structures with us and also for his advice and many stimulating discussions that helped to improve the presentation of this article.

\end{document}